\documentclass[a4paper,10pt]{article}
\pdfoutput=1
\usepackage{amsthm}
\usepackage{amsmath}
\usepackage{array}
\usepackage{amssymb}
\usepackage{graphicx}
\usepackage{chngpage}
\usepackage[all]{xy}
\usepackage{tikz}
\usepackage{url}
\usepackage[english]{babel}

\newtheorem{theorem}{Theorem}
\newtheorem{proposition}{Proposition}
\theoremstyle{definition}
\newtheorem{definition}{Definition}	
\newtheorem{example}{Example}
\newtheorem{remark}{Remark}

\title{Algebraic theories, monads, and arities\\\ \\\small{Master thesis\footnote{Under the direction of Paul-Andr\'e Melli\`es (Laboratoire Preuves, Programmes, Syst\`emes - CNRS, Universit\'e Paris Diderot - Paris, France)} - University Paris 6}}
\date{}
\author{Charles Grellois}
\begin{document}
\maketitle

\begin{abstract}
Monads are of interest both in semantics and in higher dimensional algebra. It turns out that the idea behind usual notion finitary monads (whose values on all sets can be computed from their values on finite sets) extends to a more general class of monads called \emph{monads with arities}, so that not only algebraic theories can be computed from a proper set of arities, but also more general structures like $n$-categories, the computing process being realized using Kan extensions. This master thesis compiles the required material in order to understand this question of arities and reconstruction of monads, following mostly \cite{lics2010}, and tries to give some examples of relevant interest from both semantics and higher category theory. A discussion on the promising field of operads is then provided as appendix.
\end{abstract}

\section{Introduction}
In logics or semantics, it is common to think of a theory as generated by elements together with relations between them. Every operation of the theory is then obtained by composition of the generators in a way that respects their \emph{arities}: for example, the usual addition $+$ as arity $2$ since it takes two inputs. In \cite{lawvere}, Lawvere gives a \emph{functorial} presentation of this, where there are no distinguished generators anymore, but only $n$-ary operations whose compositions respect the relations of the theory. The set of all arities is then $\mathbb{N}$: for example, in the case of arithmetics, the composite of $n$ additions gives a $(n+1)$-ary operation, and this for every natural integer $n$. These theories are in turn equivalent to \emph{finitary monads}, that is, monads whose values on finite ordinals determine values on all sets, the computation being realized by a filtered colimit.\\
But monads have an interest outside of semantics: for example, they allow to build mathematical structures such as categories, which are built out of graphs. A natural question is then to determine whether the case of finitary monads can be generalized: is there a finite set of arities, to be thought of as some set of elementary pieces, which can be glued by a proper process to realize all categories ? It happens to be the case: \emph{weighted colimits} of filiform graphs give birth to all categories, as to be shown in this paper. More generally, a notion of \emph{monads with arities} is provided (originally introduced by Weber in \cite{famfun}), together with a notion of Lawvere theories with arities, extending the usual correspondence.\\
This paper is structured as follows: we start at Section 2 by recalling usual Lawvere theories, finitary monads, and the traditional correspondence between them, after what we investigate the structure of categories to determine their set of arities, and then we introduce Kan extensions, which give a proper way of describing the computation of a monad with arities on every value, from its values on arities. An introduction to Yoneda structures is then given.\\
At Section 3, we give the general axiomatization of monads with arities, together with their corresponding Lawvere theories, and investigate the arities of the usual free category and free $2$-category monads.\\
At Section 4, we give several examples of monads encountered in the practice of semantics, together with their corresponding theories; a special treatment of the state monad is made, following \cite{lics2010}.\\
The appendix enlarges the discussion to the promising field that is the one of operads: after a short introduction to them, we explain following \cite{leinster_algeb} that, in spite of their similarities with Lawvere theories, operads are not equivalent to them; we then have a look at how the generalization of arities is treated in the operadic case.

\paragraph{Acknowledgements}
I would like to thank Paul-Andr\'e Melli\`es for his advices and the quality of the discussions we have had, as well as for the freedom I enjoyed while working on this master thesis, Jonas Frey for answering my questions with a lot of patience, the people of the working group "Cat\'egories sup\'erieures, polygraphes et homotopie" of the Laboratoire Programmes, Preuves, Syst\`emes for the interesting talks on the notions of operads and arities they gave, and to the contributors to the nLab for their amazing and very useful work. Special thanks go to Antoine Delignat-Lavaud for his support in improving the language quality of this document.\\
I also would like to apologize by advance to all the people whose paternity in ideas detailed here was not mentioned due to my lack of knowledge of the history of the field.
\newpage
\tableofcontents
\newpage
\section{Preamble}

\subsection{Monads and theories}

\subsubsection{Lawvere theories}
Historically, algebraic theories such as monoids, groups, Lie algebras \ldots\ used to be presented by means of generators and relations between them. In his doctoral dissertation \cite{lawvere}, Lawvere introduced in 1963 an alternative method of specification of algebraic theories using categories. 

\begin{definition}[Lawvere theory]
A Lawvere theory is a category $\mathbb{L}$ with finite products, in which every object is isomorphic to a finite cartesian power $x^n$ of a distinguished object $x$, called the \emph{generic object} of the theory $T$.
\end{definition}

There is a category of Lawvere theories, with morphisms the product-preser-ving functors between these theories such that the distinguished object of the source category is sent to the distinguished object of the target category.\\
The idea behind Lawvere theories is to represent every $n$-ary operation of a theory $\mathbb{L}$ as a morphism from $x^n$ to $x$ - here, the relations of the usual presentation of algebraic structures are encoded in the composition law of the category $\mathbb{L}$. Moreover, there is no notion of "primitive" operation in a Lawvere theory: there are no specified generators, but rather all the operations are given in $\mathbb{L}$.\\
Just like a given group is a model of the algebraic theory of groups, there is a notion of model of a Lawvere theory:

\begin{definition}[Model of a Lawvere theory]
A model of a Lawvere theory $\mathbb{L}$ in a category $\mathcal{C}$ is a finite-products preserving functor $A\,:\, \mathbb{L} \rightarrow \mathcal{C}$ - that is, a functor such that $A[x^{m_1} \times \ldots \times x^{m_k}] \cong A[x^{m_1}] \times \ldots \times A[x^{m_k}]$.
\end{definition}

A straightfoward consequence of this definition is that $A[n] \cong A[1]^{n}$ for every natural number $n$: thus, the functor $A$ sends every $n$-ary operation of the Lawvere theory $\mathbb{L}$ to an operation $A[1]^{n} \rightarrow A[1]$, that is, to a $n$-ary operation on the object $A[1]$, which is called the \emph{underlying object} of the $\mathbb{L}$-model $A$.\\
Again, there is a notion of category of models: given a Lawvere theory $\mathbb{L}$ and a category $\mathcal{C}$, Mod($\mathbb{L},\mathcal{C}$) has the natural transformations between models as morphisms.

\begin{example}[The theory of sets]
The theory of sets is the algebraic theory without operations. Its Lawvere theory $\mathbb{S}$ is just the free category with finite products generated from the category with one object $x$. The free construction of finite products generates $\mathbb{S}$ as the category with objects the (natural) powers $x^n$ of $x$, and morphisms the projections $\pi_i\,:\,x^n \rightarrow x$. $\mathbb{S}$ is thus the opposite of the category $FinSet$ of finite sets $\{1,\ldots,n\}$ and functions between them.\\
Since $\mathbb{S}$ is built freely from a single object $x$, and since the models of Lawvere theories are \emph{product-preserving} functors, its models in $Set$ are entirely characterized by the image of $x$: thus, Mod($\mathbb{S},Set$) is equivalent to $Set$, and $\mathbb{S}$ is the Lawvere theory of sets. More generally, Mod($\mathbb{S},\mathcal{C}$) is equivalent to the category $\mathcal{C}$ itself.
\end{example}

Moreover, this theory is the initial object in the category of Lawvere theories; this induces for every Lawvere theory $\mathbb{L}$ a functor: 
\begin{center}
$U_\mathbb{L}\,:\,$Mod($\mathbb{L},Set$)$ \rightarrow$ Mod($\mathbb{S},Set$)$\,=\,Set$
\end{center}
which sends every model of the theory $\mathbb{L}$ to its \emph{underlying set}. Moreover, this functor has a left adjoint $F_\mathbb{L}\,:\,Set \rightarrow $Mod($\mathbb{L},Set$), called the free $\mathbb{L}$-model functor. As expected, this functor sends a set $S$ to the $\mathbb{L}$-model whose underlying set is the set of formal expressions $\{f(s_1,\cdots ,f_n)/f\in L(n,1), s_i \in S\}$.\label{free}\\
This generalizes to models in other categories than $Set$: since $\mathbb{S}$ is initial, there is a canonical morphism $j\,:\,\mathbb{S} \rightarrow \mathbb{L}$ for every Lawvere theory $\mathbb{L}$, inducing a functor $U_j\,:\,$Mod($\mathbb{L},\mathcal{C}$)$\rightarrow$Mod($\mathbb{S},\mathcal{C}$) (given by the precomposition of the $\mathbb{L}$-model by $j$), which transports a given model to its \emph{underlying object} in $\mathcal{C}$. There is a corresponding notion of free model (when $\mathcal{C}$ is cartesian closed with small colimits), obtained again as left adjoint to $U_j$; we will see a way to compute it after the introduction of Kan extensions.

\begin{example}[The theory of monoids in a category] We now consider the theory $\mathbb{M}$ with objects the natural numbers and freely generated as a category with products whose morphisms from $n$ to $1$ are the finite words over $n$ letters $1,\ldots ,n$. Then Mod($\mathbb{M},\mathcal{C}$) is equivalent to the category of monoids in $\mathcal{C}$ together with their morphisms.\\ 
\end{example}

\begin{example}[The theory of groups]\label{lawgroup} We denote by $F(n)$ the free group on $n$ generators. Let $\mathbb{G}$ be the category opposite to the category of free groups with objects $F(n)$ for every natural integer $n$ and with morphisms the group homomorphisms between them. $F(1)$ then stands as the generic object of $\mathbb{G}$, and $\mathbb{G}$ is a Lawvere theory since it has finite products: this comes from the fact that the category of free groups has finite coproducts, where $G \, {\scriptstyle \coprod}\, H$ is given in the category of groups (and thus in the subcategory of free groups) by the free group built over $G$ and $H$. So $F(m)\, {\scriptstyle \coprod} \, F(n) \cong F(m + n)$ in $FreeGroup$, and $F(n) \cong F(1)^n$ in $\mathbb{G}$.\\
Now that we defined $\mathbb{G}$ and showed that it is a Lawvere theory, we show that its models in $Set$ are precisely the groups. First, any group $G$ defines a model of $\mathbb{G}$, since it defines a product-preserving functor \begin{displaymath} FreeGroup^{op} \hookrightarrow Group^{op} \xrightarrow{Hom(\_,G)} Set \end{displaymath}
Conversely, let $M\,:\,\mathbb{G}\rightarrow Set$ be a model of $\mathbb{G}$, and $G=M(x)$ be the underlying set of $M$. We have to define a group structure on $G$; for this we need a multiplication.  In $FreeGroup$ there is a group homomorphism $F(1)\rightarrow F(2)$, which arises from an element $1\rightarrow F(2)$ by the structural properties of homomorphims between free groups. If we denote by $a$ and $b$ the two generators of $F(2)$, this element is just $*\mapsto ab$, and in $\mathbb{G}$ this gives a multiplication on the generic group $x$, namely $m\,:\,x^2 \rightarrow x$. Since $M$ is a product-preserving functor, we obtain a map $M(x)\times M(x) \cong M(x^2) \xrightarrow{M(m)} M(x)$ which properly defines a multiplication over $G=M(x)$. The group identity and group inversion in $G$ are obtained similarly: models of $\mathbb{G}$ are thus equivalent to groups.
\end{example}

\subsubsection{Monads}
\begin{definition}[Monad]
A monad in a category $\mathcal{C}$ is a triple $(T,\mu,\eta)$ where $T\,:\,\mathcal{C}\rightarrow\mathcal{C}$ is a functor and $\mu\,:\,T \circ T\rightarrow T$ and $\eta\,:\,1\rightarrow T$ are natural transformations satisfying the following diagrammatic conditions:\\
\begin{center}
\begin{tabular}{ccc}
  \xymatrix{
    T\circ T \circ T \ar[r]^{\mu T} \ar[d]_{T \mu} & T \circ T \ar[d]^\mu \\
    T \circ T \ar[r]_\mu & T
  }

& \ \ \ \ &
  \xymatrix{
    T \ar[r]^{T \eta} \ar[rd]_{id} & T \circ T \ar[d]_\mu & T \ar[l]_{\eta T} \ar[ld]^{id} \\
   \ & T &\ \\
  }

\\
\end{tabular}
\end{center}
\end{definition}

The idea of a monad is to endow an endofunctor with a monoid-like structure on its iterated powers (for composition): $\mu$ stands as a multiplication, whose associativity is ensured by the diagram on the left, and $\eta$ as unit for this multiplication, as the diagram on the right shows. Remark that there is no way to get rid of $T$: the notion of monad is an abstract specification, whose models are called algebras.

\begin{definition}[Algebra for a monad]
Given a monad $(T,\mu,\eta)$ in a category $\mathcal{C}$, an algebra for $T$ is a couple $(X,h)$, where $X$ is an object of $\mathcal{C}$ and $h\,:\,TX\rightarrow X$ is an arrow, with the following diagrammatic conditions:\\
\begin{center}
\begin{tabular}{ccc}
  \xymatrix{
    (T\circ T)X \ar[r]^{\mu_X} \ar[d]_{Th} & T(X)\ar[d]^{h} \\
    T(X)\ar[r]_h & X
  }

& \ \ \ \ &
  \xymatrix{
    X \ar[r]^{\eta_X} \ar[rd]_{id} & TX \ar[d]^h \\
   \ & X\\
  }

\\
\end{tabular}
\end{center}
\end{definition}

\begin{example}[The free category monad]
In the category $Graph$ of graphs and morphisms between them, we consider the functor $T$ which sends a graph to its free category (still considered as a graph). This functor takes a graph as input and outputs the category obtained from this graph by adding all required identities and  finite paths obtained associatively by composition as arrows. We need to define a multiplication and unit for it satisfying the relations given in the definition. Let $G$ be a graph, $TG$ is then its free category, and $TTG$ is the free category over the free category of $G$: by construction, its arrows are paths of (composable) paths of $G$. Thus, $\mu\,:\,TTG\rightarrow TG$ can be naturally defined as the natural transformation which composes paths of paths to paths. $\eta$ is defined as the inclusion of the graph $G$ into its free category (considered as a graph): it sends every arrow of the graph $G$ on the corresponding path of length 1 in $TG$. The associativity of $\mu$, corresponding to the diagram on the left in the definition of a monad, is clear: when considering paths of paths of paths, the way to compose them to finally obtain just paths is irrelevant. The diagram on the right, expressing that $\eta$ stands as unit for $\mu$, should be checked more carefully. We consider first the required equality $\mu \circ T\eta\,=\,1$, which instanciates on $TG$ as $\mu \circ T\eta_{G}\,=\,TG$. Since $\eta_G$ simply includes $G$ into $TG$, the heterogeneous composition $T\eta_G$ includes $TG$ in $TTG$, in which the arrows are just paths of arrows of $\eta_G$, that is paths of paths of length one: their concatenation gives back $TG$. The second required equality on every graph $G$ is $\mu \circ \eta_{TG}\,=\,TG$: $\eta_{TG}$ includes $TG$ in $TTG$, in which the arrows are path of length one of paths of $G$. Their concatenation, again, gives the (underlying graph of the) free category $TG$: all the required conditions on $\mu$ and $\eta$ are satisfied, and $(T,\mu,\eta)$ is a monad on $Graph$.\\
One may then wonder which are the algebras for this monad. First of all, every categorical graph $G$ (that is, a category seen as its underlying graph) defines an algebra together with an arrow $h\,:\,TG\rightarrow G$ whose effect is to take an arrow of $TG$, which is seen as a path of a certain length, and to send it to the arrow of $G$ which results from the composition of the arrows of this path according to the categorical structure on $G$. Conversely, if $(G,h)$ is an algebra for the monad $T$, $h$ has to satisfy the equation $h \circ Th\,=\,h \circ \mu_G$. This implies that $h$ gives a notion of composition of paths on $G$, which has to be associative from the associativity in $TG$, and such that $G$ has identities (again, because $TG$ has them). Remark that the algebras for $T$ corresponding to free categories are exactly the ones of the shape $(TG,h)$.
\end{example}

\begin{example}[The free group monad]
In $Set$, we consider the functor $T$ sending a given set $S$ to the free group on it, which is the group generated by the elements of $S$, an unit element $1$, and such that if $E$ and $F$ are in $TS$, $EF$ and $E^{-1}$ also are. Furthermore, there are relations between the elements: if $E,F,G$ are elements of $TS$, we require that:
\begin{itemize}
\item $(EF)G\,=\,E(FG)$,
\item $E1\,=\,E\,=\,1E$,
\item $EE^{-1}\,=\,1\,=\,E^{-1}E$
\end{itemize}
For the sake of clarity, we denote by $[a]$ the element corresponding to $a\in S$ in $TS$, and similarly use brackets to distinguish elements of $S$, $TS$ or even $TTS$. The action of T on arrows is defined inductively: if $f\,:\,X\rightarrow Y$ is a map, $T(f)$ is defined by: 
\begin{itemize}
\item $T(f)([a])\,=\,[f(a)]$ for $a\in S$, 
\item $T(f)(1)\,=\,1$,
\item $T(f)(EF)\,=\,T(f)(E)T(f)(F)$,
\item $T(f)(E^{-1})\,=\,(T(f)(E))^{-1}$
\end{itemize}
$T$ is clearly a functor. We now consider as unit of the monad on $T$ the natural transformation $\eta$ which, on every set $S$, is the application $a\mapsto [a]$, and $\mu_S\,:\,(T\circ T)(S)\rightarrow TS$ sends an element of $(T\circ T)S$ (that is, an element written with two levels of brackets) to an element of $TS$ by removing the external brackets and distributing the inversion if needed. For example: 
\[
\mu_S([[a][b]^{-1}][[c][a]][[d]][[a][b]]]^{-1})\,=\, [a][b]^{-1}[c][a][d][b]^{-1}[a]^{-1}
\]
It should be clear that $\mu$ satisfies the left diagram, since removing the external brackets twice in an expression with three levels of bracketing has the same result than removing the middle ones then the external ones. Now we have to show that $\eta$ stands as unit for $\mu$. On $S$, $T\eta_S$ sends an element $E$ of $TS$ to the element of $TTS$ obtained by puting double brackets instead of simple brackets on each element originally coming from $S$, whereas $\eta_S T$ sends $E$ to $[E]$, simply putting a new bracket on the whole element. For example: 
\[
T\eta_S ([a][c]^{-1})\,=\,[[a]][[c]]^{-1}
\]
and:
\[
\eta_S T([a][c]^{-1})\,=\,[[a][c]^{-1}]
\]
It should be clear that $\mu_S$, acting by removing the external brackets, gives back in both cases $[a][c]^{-1}$. This way, the two required equations $\mu \circ T \eta \, =\, 1$ and $\mu \circ \eta T\,=\, 1$ are satisfied, and we have a monad structure on $T$.\\
Now, what are the algebras for $T$ ? If $(S,h)$ is an algebra for $T$, $h\,:\,TS\rightarrow S$, $h$ sends an element of the free group $TS$ to $S$. In particular, $[a][b]$ is sent to $h([a][b])$, endowing $S$ with a multiplicative structure on $S$ whose unit is $h(1)$: $h$ gives birth to a group structure on $S$. Conversely, every group structure on a set $S$ defines an arrow $h\,:\,TS\rightarrow S$ computing the result of an expression in this group, and this arrow makes the two diagrams defining algebras commute.
\end{example}

We just made the two previous descriptions of monads explicit, but their existence may have been deduced from the fact that the associated free functor is a left adjoint to the forgetful functor ($Cat \rightarrow Graph$ and $Group\rightarrow Set$, respectively), and from the following construction:

\begin{definition}[Monad of an adjunction]
Any adjunction $F\dashv G\,:\,\mathcal{C}\rightarrow \mathcal{D}$, with unit $\eta$ and counit $\epsilon$, induces a monad $T\,=\,G\circ F\,:\,\mathcal{C}\rightarrow \mathcal{C}$, whose unit $\eta\,:\,1\rightarrow T$ is precisely the one of the adjunction, and whose multiplication is $\mu\,=\,G\epsilon F\,:\,GFGF\rightarrow GF$.
\end{definition}


We then have notions of morphisms between monads and between algebras, leading to the definition of the category of monads and to the one of algebras of a monad.

\begin{definition}[Morphism of monads]
A morphism between a monad $(T_1,\mu_1,\eta_1)$ over a category $\mathcal{C}_1$ and a monad $(T_2,\mu_2,\eta_2)$ over a category $\mathcal{C}_2$ is given by a functor $\Phi\,:\,\mathcal{C}_1\rightarrow \mathcal{C}_2$ and a natural transformation $\rho\,:\,T_2 \Phi \Rightarrow \Phi T_1$ making the two following diagrams (of natural transformations) commute:
\begin{center}
\begin{tabular}{ccc}
\xymatrix {
T_2 T_2 \Phi \ar@{=>}[r]^{T_2 \rho} \ar@{=>}[d]_{\mu_2 \Phi} & T_2 \Phi T_1 \ar@{=>}[r]^{\rho T_1} & \Phi T_1 T_1 \ar@{=>}[d]^{\Phi \mu_1}\\
T_2 \Phi \ar@{=>}[rr]_{\rho} & \ & \Phi T_1\\
}
&\ \ \ &
\xymatrix {
\Phi \ar@{=>}[r]^{\Phi \eta_1} \ar@{=>}[d]_{\eta_2 \Phi} & \Phi T_1\\
T_2 \Phi \ar@{=>}[ur]_{\rho} & \ \\
}
\end{tabular}
\end{center}
\end{definition}

\begin{definition}[Morphism of algebras] 
If $(X,h)$ and $(Y,k)$ are two algebras over a monad $T$, a morphism $(X,h)\rightarrow (Y,k)$ is an arrow $f\,:\,X\rightarrow Y$ satisfying:\\
\begin{center}
\begin{tabular}{c}
  \xymatrix{
    TX \ar[r]^{h} \ar[d]_{Tf} & X \ar[d]^f \\
    TY \ar[r]_k & Y\\
  }
\end{tabular}  
\end{center}
\end{definition}

The category of algebras of a monad is often denoted $\mathcal{C}^T$, and is called the Eilenberg-Moore category of the monad. Among the algebras, some will be of special interest for us: these are the free algebras.

\begin{definition}[Free algebras] 
An algebra for a monad $T$ is \emph{free} when it is of the form $(TX,h)$.
\end{definition}

This notion of \emph{freeness} is precised by the fact that these algebras arise from a free functor: the functor $U\,:\,\mathcal{C}^T\rightarrow \mathcal{C}$, sending an algebra  $(X,h)$ to its underlying object $X$, has a left adjoint sending an object $X$ to the free algebra $(TX,\mu_X)$. The restriction of the Eilenberg-Moore category to its subcategory of free algebras gives a full subcategory. Another presentation of this subcategory exists:

\label{kleisli}
\begin{definition}[Kleisli category] 
Given a category $\mathcal{C}$ and a monad $(T,\mu,\eta)$ over it, the Kleisli category of the monad is denoted as $\mathcal{C}_T$, has the same objects as $\mathcal{C}$, and its arrows $A\rightarrow B$ are the arrows $A\rightarrow TB$ in  $\mathcal{C}$. The composition of $f\,:\,A\rightarrow B$ and $g\,:\,B\rightarrow C$ is then defined as $\mu_C \circ Tg \circ f$, that is:\\
\begin{center}
\begin{tabular}{c}
  \xymatrix{
    \ &\ & TTC \ar[d]^{\mu_C} \\
    \ &  TB \ar[ur]^{Tg} & TC \\
    A \ar[ur]^{f} & B \ar[ur]^{g} & \ \\
    }
\end{tabular}  
\end{center}
and unit arrow $A\rightarrow A$ is given by $\eta_A$, which is an arrow $A\rightarrow TA$ in $\mathcal{C}$ and thus an arrow $A\rightarrow A$ in $\mathcal{C}_T$.

\end{definition}

\begin{proposition}
The Kleisli category of a monad $T$ is equivalent to the full subcategory of $\mathcal{C}^T$ whose objects are the free algebras.
\end{proposition}

\begin{definition}[Finitary monad] 
A monad $T\,:\,Set \rightarrow Set$ is \emph{finitary} when it preserves filtered colimits.
\end{definition}

An equivalent definition will be given after the introduction of Kan extensions. The interesting property of finitary monads is that they arise from their values on the finite ordinals  $[n]\,=\,\{0,\ldots,n\}$ (considered as totally ordered sets):

\begin{definition}[Compact object] In a locally small category $\mathcal{C}$ that admits filtered colimits, an object $X \in \mathcal{C}$ is \emph{compact} if the corepresentable functor $\operatorname{Hom}_C (X,\_)\,:\,\mathcal{C}\rightarrow Set$ preserves these filtered colimits.
\end{definition}

The compact objects of $Set$ are precisely the finite sets, so up to isomorphism every set $S$ can be recovered as a filtered colimit of finite ordinals (the idea is to take the filtered colimit over the poset of finite subsets of $S$ with their inclusions). Since filtered colimits commute with finitary monads, and since they are endofunctors of $Set$, the value of a given finitary monad on any set can be computed from its values on finite ordinals.

\subsubsection{Monads and theories} We have seen that groups can be retrieved from free groups as models of a Lawvere theory or equivalently as algebras for a (finitary) monad. We can also see in this example that the Lawvere theory is very close to the data of the monad computed on finite ordinals. In fact, there is a more general result: since a finitary monad is defined by its values on all finite ordinals, $T[n]$ can be thought of as the set of $n$-ary operations of the monad $T$, leading to the following correspondence:

\begin{theorem}[Equivalence between Lawvere theories and monads] The category of Lawvere theories is equivalent to the one of finitary monads.

\end{theorem}

A (complicated) proof of this result may be found in the chapter 3 of \cite{adamek}. The idea is that for $T$ a finitary monad, the corresponding Lawvere theory $\Theta_T$ is given by the subcategory of $Set_T$ whose objects are the finite ordinals.\label{intuitn} The intuition is that the $n$-ary operations of the Lawvere theory correspond to the elements of the free algebra $Tn$ generated by $n$ elements. In addition to this equivalence, it holds that the theory and the monad generate the same models/algebras, generalizing what we observed in the free group example: 

\begin{theorem}[Equivalence between models and algebras] The category $Set^T$ of algebras of a finitary monad $T$ is equivalent to the category Mod($\Theta_T^{op}$,$Set$).
\end{theorem}

Some examples of the equivalence, coming from computer science, may be found at Section \ref{compsci}. 


\begin{remark}
A technical advantage of Lawvere theories is that their models can be chosen in any category with finite products, whereas finitary monads require this category to be $Set$. In the following, we will work on extending this to a more general class of categories, by extending the fundamental property of finitary monads, which is that they can be computed from a given full subcategory of $Set$, in a suitable way. This will also generalize the notion of Lawvere theory, allowing operations with a bigger (or different) class of arities: this allows usual generalizations of Lawvere theories like \emph{countable} theories, with arities the ordinals $[0], \ldots, [\omega ]$ for example.
\end{remark}

\subsection{The nerve functor: a simplicial structure of categories}

Our quest for a generalized notion of arities is guided by an example which arises from algebraic topology, where a space's properties can be investigated through the use of simplicial homology. In this approach, the space is modelized - up to homeomorphism - by the data of simplices of any dimension, of which copies are glued together according to extra informations corresponding to face inclusions and "degeneracies". In category theory, this idea turns out to be a good way to describe higher dimensional structures - by describing how they look like locally and how to paste together the local data to get the whole structure. In some sense, the case of the finitary monads was similar: from the description of the monad's values on finite ordinals (corresponding to \emph{arities} of the theory), we can reconstruct by a given process (filtered colimits) all the values of the monad. One also may see a similarity between this process and countable dense bases in Hilbert spaces: from a given subset of values, a process of limit gives all the data on the space. We should then investigate which may be the proper notions of arities, what the meaning of density is, and find a proper generalization of this colimit/glueing process.

\subsubsection{The simplicial category, simplicial sets} We start by formalizing the idea of simplices in a categorical manner.
\begin{definition}
The simplicial category $\Delta$ has objects all finite ordinals\\\mbox{$n\,=\,\{0,\ldots,n\}$} (considered as totally ordered sets), and arrows all weakly monotone functions: $\Delta ([n],[m])$ is then the set of functions $f\,:\,[n]\rightarrow [m]$ such that for every $(i,j)\in [n]^2$, $i\leq j\Longrightarrow f(i) \leq f(j)$. The \emph{augmented} simplicial category $\Delta_a$ is $\Delta$ enriched with an object $[-1]$, corresponding to the empty set.
\end{definition}

$\Delta_a$ has $[-1]$ as initial element and $[0]$ as final one (so, $\Delta$ only has a final element, but no initial element). There is a bifunctor $+\,:\,\Delta_a \times \Delta_a\rightarrow \Delta_a$ sending the couple of ordinals $([n],[m])$ to the ordinal $[n+m+1]$ and sending two arrows $f\,:\,n\rightarrow n'$ and $g\,:\,m\rightarrow m'$ to the arrow $f+g$ defined as follows:
\begin{center}
$
(f+g)(i) =
\begin{cases}
 f(i), & \text{if }i \in [n] \\
 n'+g(i-n), & \text{if }n\in \{n+1,\ldots,n+m+1\}
\end{cases}
$
\end{center}

This bifunctor clearly gives to $(\Delta_a, +,[-1])$ a strict monoidal structure (note that $\Delta$ doesn't have one for $+$ since it has no unit for it). The interest of this bifunctor $+$ is to give a convenient (and geometrical) way to describe every arrow in $\Delta_a$, and thus on its subcategory $\Delta$. First of all, since $[0]$ is terminal in $\Delta_a$, there is a unique arrow $\mu\,:\,[1]\rightarrow [0]$ and another unique arrow $\eta\,:\,[-1]\rightarrow [0]$\footnote{$([0],\mu,\eta)$ is a monoid in $\Delta_a$, which is universal in a sense made precise in \cite[Section VII.5]{maclane}}. We can use these two arrows, identities, and the monoidal structure to build every arrow of $\Delta_a$ and thus of $\Delta$. For every $n$-simplex, we define:

\begin{itemize}
\item the $i^{th}$ face inclusion $\delta_i^n\,=\,1_{[i-1]}+\eta+1_{[n-i]}$ for $i \in [n+1]$, which is by monoidality an arrow $[n]\rightarrow [n+1]$. It corresponds to the injection whose image leaves out $i \in [n+1]$
\item the $i^{th}$ face degeneracy $\sigma_i^n\,=\,1_{[i-1]}+\mu+1_{[n-i-1]}$ for $i \in [n]$, which is by monoidality an arrow $[n+1]\rightarrow [n]$. It corresponds to the surjection whose image is the same on $i$ and $i+1$ in $[n]$.
\end{itemize}

Moreover, the binary composities of these functions satisfy the following equalities ($n$ being considered fixed and therefore omitted for the sake of readability):

\begin{itemize}
\item $\delta_i \circ \delta_j\,=\,\delta_{j+1} \circ \delta_i$ if $i\leq j$
\item $\sigma_j \circ \sigma_i\,=\, \sigma_i \circ \sigma_{j+1}$ if $i \leq j$
\item $
\sigma_j \circ \delta_i\, =\,
\begin{cases}
 \delta_i \circ \sigma_{j-1} & \text{if } i< j\\
 1 & \text{if } i\in \{j,j+1\}\\
 \delta_{i-1} \circ \sigma_j & \text{if } i>j+1\\
\end{cases}
$
\end{itemize}

An interesting result then gives the shape of every arrow of $\Delta_a$, and thus of $\Delta$:
\begin{proposition}
In $\Delta_a$, every arrow $f\,:\,[n]\rightarrow [m]$ has a unique representation:
\begin{center}
$f\,=\,\delta_{i_1} \circ \ldots \circ \delta{i_k} \circ \sigma_{j_1}\circ \ldots \circ \sigma_{j_h}$
\end{center}
where $h$ and $k$ are such that $n-h+k\,=\,m$ and where the subscripts satisfy:
\begin{center}
$m \geq i_1>\ldots > i_k \geq 0$ and $0 \leq j_1 < \ldots < j_h < n$
\end{center}
\end{proposition}

Informally, the two ordering conditions on the subscripts are due to the commutation equations given just above, while the condition $n-h+k\,=\,m$ expresses the fact that the number of increases and decreases of dimension finally meets the one required to go from the $n^{th}$ standard simplex to the $m^{th}$. We thus obtain $\Delta_a$ as the category whose objects are the finite ordinals together with the empty set $[-1]$, and whose arrows are obtained by composition of face and degeneracy operations, subject to the relations expressed by the equations above\footnote{The proofs for these results may be found in \cite[Section VII.5]{maclane}.}.

\begin{remark}[On other definitions of $\Delta$]
Several definitions of $\Delta$ (and $\Delta_a$) may be found in the litterature - some people call ours the \emph{skeletal version} of $\Delta$. In this case $\Delta$ is defined as having all finite sets (again, totally ordered) as objects and order-preserving functions between them as morphisms, whereas in our definition the objects just are the isomorphism classes of these sets. With our definition, there is only \emph{one} version of a given $n$-simplex, and the numbering of the vertices doesn't necessarily meet the one of subsimplices: e.g., the triangle $[2]$ is mapped in all cases to $[1]$ with our definition by the degeneracy maps, whereas in the other definition it would be mapped to $\{0,1\}$, $\{0,2\}$ or $\{1,2\}$, following the numberings of the vertices. But reasoning with isomorphism classes does not bring any difficulty since the real relationship between simplices is given by the degeneracy functions (and dually by the face inclusion ones), and not by the numbering of the vertices.
\end{remark}

From $\Delta$ we retrieve an usual concept from algebraic topology:

\begin{definition}[Simplicial sets]
A presheaf over $\Delta$ is called a simplicial set.
\end{definition}

The idea of this definition is that such a presheaf $X\,:\,\Delta^{op}\rightarrow Set$ is given by a family $(X_n)_{n \in \mathbb{N}}$ of sets, describing copies of $n$-simplices, and a map $X_f\,:\,X_q \rightarrow X_p$ for every function $f\,:\,[p] \rightarrow [q]$, describing the relations between these simplices - that is, how to glue them to obtain a topological space built from standard simplices. This idea will be detailed after the introduction of the notion of weighted colimits in Section \ref{weight}.

\begin{example}[A simplicial set]
Consider the following triangle:
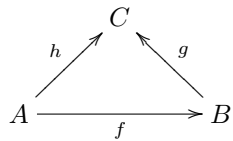
\begin{center}
\begin{tabular}{c}
\xymatrix{
& C &\\
A \ar[rr]_f \ar[ur]^h & & B \ar[ul]_g
}
\end{tabular}
\end{center}
It is described by the simplicial set $X$ such that:
\begin{itemize}
\item $X(0)\,=\,\{A,B,C\}$
\item $X(1)\,=\,\{f,g,h\}$
\item $X(2)$ is the full triangle
\item The images of the functions of $\Delta$ describe how to glue this data to build the triangle: for example, the two applications $[0] \rightarrow [1]$ in $\Delta$:
\begin{center}
\begin{tabular}{ccc}
$0 \mapsto 0$ & and & $0 \mapsto 1$
\end{tabular}
\end{center}
are the following maps (of sets) $X(1) \rightarrow X(0)$:
\begin{center}
\begin{tabular}{ccc}
$\left. \begin{aligned} f \mapsto A\\ g \mapsto B\\ h \mapsto A\\ \end{aligned} \right\}$ & and & $\left\{ \begin{aligned} f \mapsto B\\ g \mapsto C\\ h \mapsto C\\ \end{aligned} \right.$
\end{tabular}
\end{center}
which are just the edge source and edge target maps.
\end{itemize}
\end{example}

\subsubsection{The nerve functor} 
There is an inclusion functor $i\,:\,\Delta \rightarrow Cat$ which embeds $\Delta$ in $Cat$\footnote{We denote by $Cat$ the category of \emph{small} categories and their functors, and by $CAT$ the full one.}; it acts by sending $[n]$ to the free category over the linear quiver $0\rightarrow \cdots \rightarrow n$, which basically is the category induced by the usual order on $\{0,\cdots,n\}$. The action of $i$ on monotone functions derives from its action on the face inclusion and degeneracy maps (due to the factorization lemma in $\Delta_a$), which is natural with the geometrical intuition of $\Delta$ we have developped: the face inclusions are mapped to the functors applying the same inclusion of the category $i[n]$ in $i[n+1]$, and the degeneracies collapse $i[n+1]$ to $i[n]$ by sending the removed vertex to another one.\\

Since $Cat$ is locally small, $i$ induces a functor:
\begin{center} $
\begin{array}{llcc}
Cat(i,1)\,:\, & Cat & \rightarrow & \widehat{\Delta}\\
& \mathcal{C} & \mapsto & 
\bigg\{ \begin{array}{lcr}
[n] & \mapsto & Cat(i[n],\mathcal{C})\\
f & \mapsto & (if)^{*}\\
\end{array}\\ \end{array} $
\end{center}

This is a special case of Yoneda structure, a notion to be introduced at Section \ref{yoneda}. This functor is called the \emph{nerve} of a category. It gives to every category a structure of simplicial set, where the vertices are the objects of the category, the edges ($1$-simplices) its morphims, the triangles ($2$-simplices) are the usual commutative triangles derived from composition:
\begin{center}
\begin{tabular}{c}
  \xymatrix{
    B \ar[r]^g & C \\
    A \ar[ru]_{g\circ f} \ar[u]^{f}&\ \\
  }
\end{tabular}  
\end{center}

The $3$-simplices then are of the following shape:

\begin{center}
\begin{tabular}{c}
 
\xymatrix @!0 @R=3pc @C=4pc {
    \ & D &\  \\
    &&\\
    A \ar[uur]^{h\circ g \circ f} \ar@{.>}[rr]_{\ \ \ \ \ \ \ g\circ f} |!{[ur];[dr]}\hole \ar[rd]_f && C \ar[luu]_h \\
    & B \ar[ru]_g \ar[uuu]^{h \circ g} }

\end{tabular} 
\end{center}

Beware, the identities \emph{are} arrows of a category and thus the following \emph{are} $3$-simplices:\\

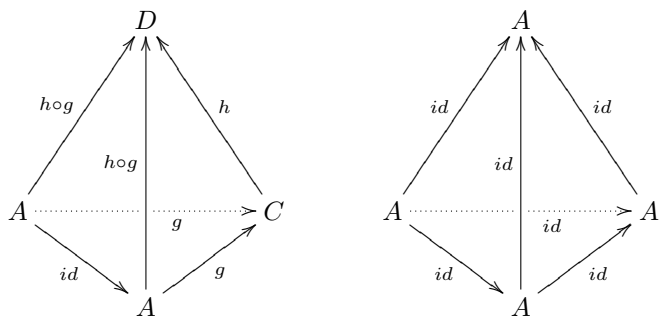
\begin{center}
\begin{tabular}{ccc}
 
\xymatrix @!0 @R=3pc @C=4pc {
    \ & D &\  \\
    &&\\
    A\ar[uur]^{h\circ g} \ar@{.>}[rr]_{\ \ \ \ \ \ \ g} |!{[ur];[dr]}\hole \ar[rd]_{id} && C \ar[luu]_h \\
    & A \ar[ru]_g \ar[uuu]^{h \circ g} }&\ \ \ &
\xymatrix @!0 @R=3pc @C=4pc {
    \ & A &\  \\
    &&\\
    A \ar[uur]^{id} \ar@{.>}[rr]_{\ \ \ \ \ \ \ id} |!{[ur];[dr]}\hole \ar[rd]_{id} && A \ar[luu]_{id} \\
    & A \ar[ru]_{id} \ar[uuu]^{id} }\end{tabular} 
\end{center}

\begin{example}[An example of nerve]
Consider the following category:
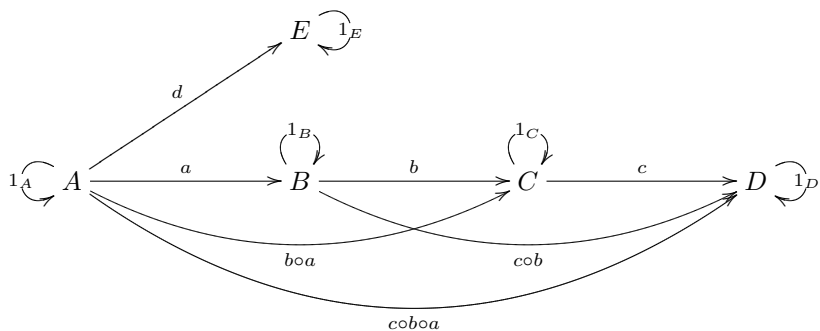
\begin{center}
\begin{tabular}{c}
\xymatrix @!0 @R=20mm @C=30mm {
\ & E \ar@(ur,dr)[]|{1_E} & & \ \\
A \ar[ur]^d \ar@(ul,dl)[]|{1_A} \ar[r]^a \ar@/_2pc/[rr]_{b\circ a} \ar@/_4pc/[rrr]_{c\circ b\circ a} & B \ar@(ul,ur)[]|{1_B} \ar[r]^b \ar@/_2pc/[rr]_{c\circ b} & C \ar[r]^c \ar@(ul,ur)[]|{1_C} & D \ar@(ur,dr)[]|{1_D}\\
}\\
\end{tabular}
\end{center}\ \\\ \\\ \\
Its nerve is the simplical set $X$, which gives on objects:
\begin{itemize}
\item $X[0]$ contains every object together with its identity:
\begin{center} \begin{tabular}{c}\xymatrix{ A \ar@(ur,dr)[]^{1_A}}\end{tabular}, \begin{tabular}{c}\xymatrix{ B \ar@(ur,dr)[]^{1_B}}\end{tabular}, \ldots \end{center}
\item $X[1]$ contains, for every arrow, this arrow and the two objects (with their identities) to which it is linked: 
\begin{center} \begin{tabular}{c}\xymatrix{ A \ar@(ul,dl)[]_{1_A} \ar[r]^a & B \ar@(ur,dr)[]^{1_B}} \end{tabular} \end{center} and the other ones - pay special attention that this is an element of $X(1)$: \begin{center}\begin{tabular}{c}\xymatrix{ C \ar@(ul,dl)[]_{1_C} \ar[r]^{1_C} & C \ar@(ur,dr)[]^{1_C}} \end{tabular} \end{center}
It is therefore tempting to think that $X(0)$ is included in $X(1)$ but it is not strictly true.
\item $X[2]$ contains the two following elements:
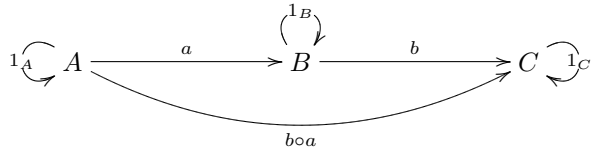
\begin{center}
\begin{tabular}{c}
\xymatrix @!0 @R=20mm @C=30mm {
A \ar@(ul,dl)[]|{1_A} \ar[r]^a \ar@/_2pc/[rr]_{b\circ a} & B \ar@(ul,ur)[]|{1_B} \ar[r]^b & C \ar@(ur,dr)[]|{1_C}\\
}\\
\end{tabular}
\end{center}\ \\
and:
\begin{center}
\begin{tabular}{c}
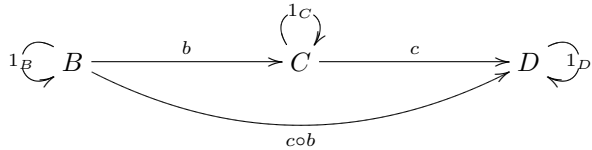
\xymatrix @!0 @R=20mm @C=30mm {
B \ar@(ul,dl)[]|{1_B} \ar[r]^b \ar@/_2pc/[rr]_{c\circ b} & C \ar[r]^c \ar@(ul,ur)[]|{1_C} & D \ar@(ur,dr)[]|{1_D}\\
}\\
\end{tabular}
\end{center}\ \\
and also two sets of elements, corresponding to the pre/postcomposition of these elements by the identity.
\item $X[3]$ contains:
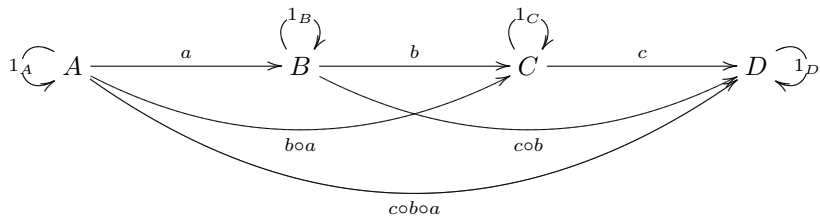
\begin{center}
\begin{tabular}{c}
\xymatrix @!0 @R=20mm @C=30mm {
A \ar@(ul,dl)[]|{1_A} \ar[r]^a \ar@/_2pc/[rr]_{b\circ a} \ar@/_4pc/[rrr]_{c\circ b\circ a} & B \ar@(ul,ur)[]|{1_B} \ar[r]^b \ar@/_2pc/[rr]_{c\circ b} & C \ar[r]^c \ar@(ul,ur)[]|{1_C} & D \ar@(ur,dr)[]|{1_D}\\
}\\
\end{tabular}
\end{center}\ \\\ \\
and degenerated cases where this shape is realized with one or more identities replacing $a$, $b$ or $c$.
\item Values of $X$ on higher arities are compositions of the previous ones with identities.
\end{itemize}
\end{example}

\begin{remark}[Higher nerves]
The definition of the nerve of a category is a bit frustrating from a topological point of view: a category only has vertices and edges, and thus we only have some kind of skeleton of something topologically bigger. In a $2$-category, we also have 2-dimensional transformations, which fill triangles:

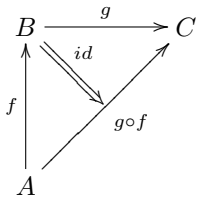
\begin{center}
\begin{tabular}{c}
  \def\enlargexyentry#1{%
    \POS "#1",*{
      \vrule height 1pt depth 1pt width 0pt
      \vrule height 0pt depth 0pt width 2pt
      }="#1",}
  \xymatrix  @!0 @R=5pc @C=5pc { 
    B \ar[r]^g & C \\ A \ar[u]^f 
    \ar[ru]|*{}="A"_{g\circ f} 
    \enlargexyentry A
    \ar@{=>}[u];"A"^{id} }\end{tabular}  
\end{center}

In a $3$-category, we obtain the pyramid as a whole "solid" built out of 4 such $2$-simplices, 6 $1$-simplices and 4 $0$-simplices, and of a $3$-dimensional transformation "filling the pyramid":
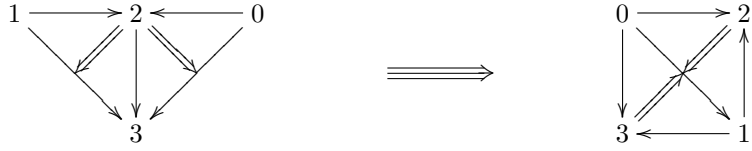
\begin{center}
\begin{tabular}{c}
\xymatrix @!0 @R=8mm @C=16mm {
1 \ar[ddr]|*{}="A" \ar[r] & 2  \ar@2{->}[];"A"^-{} \ar[dd] & 0 \ar[ddl]|*{}="B" \ar[l] \ar@2{->}[l];"B"^-{\ } & & &0 \ar[r] \ar[ddr]|*{}="C" \ar[dd] & 2 \ar@2{->}[];"C"^-{\ }\\
&&& \ar@3{->}[r] & & &\\
& 3 & &\ &\ & 3 \ar@2{->}[];"C"^-{\ } & 1 \ar[l] \ar[uu]\\
}
\end{tabular}
\end{center}
Such simplices are called \emph{orientals}. They may be defined for every dimension, and they are compatible with face and degeneracy maps (they are the translation in a $n$-categorical framework of the notion of simplex in algebraic topology). The category of orientals may then be used instead of $\Delta$ to define a notion of nerve for an $\omega$-category, see the oriental entry at nLab \cite{ncatlab} for details.\\
This higher nerve construction is useful for the construction of weak $n$-categories, see \cite[Section 10.2]{leinster} for details.
\end{remark}

The important following point should seem natural now: every category may be built from its nerve, that is built in some way from vertices (the objects), edges (the arrows), triangles (the compositions), pyramids (ensuring associativity), \ldots This is formalized by saying that $\Delta$ is \emph{dense} in $Cat$, a point we should define at Section \ref{weight}, after the introduction of weighted colimits. 

\subsubsection{The Segal condition on simplicial sets} So, every category gives birth to a simplicial set, but conversely, does simplicial sets give birth to categories ? It evidently isn't the case, think for example of a simplicial set whose image on arities $[i]$ are the empty set for $i\geq  2$: the objects and arrows it describes are not equipped with a notion of composition (moreover, the informations of glueing provided by the maps may be doubtful).\\\ \\

Segal introduced in \cite{segal} a condition named after him (yet attributed to Grothendieck) which characterizes the simplicial sets which realize as a category (that is, the ones isomorphic to the nerve of a category). The idea is quite simple after the following observation: for $([p],[q]) \in \Delta^2$, defining the functions $max$ and $min$ as follows:
\begin{center}
\begin{tabular}{ccc}$
\begin{array}{llcc}
max\,:\, & [0] & \rightarrow & [p]\\
& 0 & \mapsto & p
\end{array}$
& and & $
\begin{array}{llcc}
min\,:\, & [0] & \rightarrow & [q]\\
& 0 & \mapsto & 0
\end{array}$
\end{tabular}
\end{center}
the following square is a pushout in $\Delta$:
\begin{center}
\begin{tabular}{c}
  \def\cocartesien{%
    \ar@{-}[]+L+<-6pt,+1pt>;[]+LU+<-6pt,+6pt>%
    \ar@{-}[]+U+<-1pt,+6pt>;[]+LU+<-6pt,+6pt>%
  }
  \xymatrix{ 
    [0] \ar[r]^{max} \ar[d]_{min} & [p] \ar[d] \\
    [q] \ar[r] & [p+q] \cocartesien }
\end{tabular}
\end{center}

so that $[p+q]$ can be described as $([p] \uplus [q])/\sim$ where the $\sim$ identifies $p \in [p]$ and $0 \in [q]$: geometrically, it means that it glues the linear quiver corresponding to $[p]$, namely $0 \rightarrow \cdots \rightarrow p$, to the (shifted for the sake of comprehension) linear quiver corresponding to $[q]$, namely $p \rightarrow \cdots \rightarrow p+q$, giving $0 \rightarrow \cdots \rightarrow p \rightarrow \cdots \rightarrow p+q$, which corresponds to $[p+q]$.\\
The idea then is that a simplicial set corresponding to the nerve of some category should have a similar property, since it can be thought as a collection of simplices together with information on how to glue them. Because simplicial sets imply the use of $\Delta^{op}$, the arrows are reversed, leading to the following condition:

\begin{theorem}[Segal condition] A simplicial set $X$ is isomorphic to the nerve of a small category $\mathcal{C}$ precisely when this pushout in $\Delta$ is transported to a pullback in $Set$:
\begin{center}
\begin{tabular}{c}
  \def\cartesien{%
    \ar@{-}[]+R+<6pt,-1pt>;[]+RD+<6pt,-6pt>%
    \ar@{-}[]+D+<1pt,-6pt>;[]+RD+<6pt,-6pt>%
  }
  \xymatrix{ 
    X_{p+q} \ar[r] \ar[d] \cartesien & X_p \ar[d]^{X_{max}} \\
    X_q \ar[r]_{X_{min}} & X_0 }\end{tabular}
\end{center}

\end{theorem}

which expresses the fact that the nerves of categories are the simplicial sets whose $(p+q)$-simplices correspond to the pairs of $p$- and $q$-simplices whose extremal vertices coincide.

\begin{remark}[Segal condition and associativity] First, recall that we know explicitely what a pullback in $Set$ is: when the following square is a pullback:
\begin{center}
\begin{tabular}{c}
  \def\cartesien{%
    \ar@{-}[]+R+<6pt,-1pt>;[]+RD+<6pt,-6pt>%
    \ar@{-}[]+D+<1pt,-6pt>;[]+RD+<6pt,-6pt>%
  }
  \xymatrix{ 
    A \ar[r] \ar[d] \cartesien & C \ar[d]^{g} \\
    B \ar[r]_{f} & D}\end{tabular}
\end{center}
we have that $A \cong \{(b,c) \in B \times C / f(b)=g(c)\}$. A categorical structure on a graph requires:
\begin{itemize}
\item a notion of binary composition for two given composable arrows: we have that 
\begin{center}
\begin{tabular}{c}
  \def\cartesien{%
    \ar@{-}[]+R+<6pt,-1pt>;[]+RD+<6pt,-6pt>%
    \ar@{-}[]+D+<1pt,-6pt>;[]+RD+<6pt,-6pt>%
  }
  \xymatrix{ 
    X_2 \ar[r] \ar[d] \cartesien & X_1 \ar[d]^{X_{max}} \\
    X_1 \ar[r]_{X_{min}} & X_0 }\end{tabular}
\end{center}
is a pullback in $Set$, and thus that $X_2 \cong \{(A\rightarrow B,B\rightarrow C)\}$: every $2$-simplex in $X$ is in bijection with a couple of composable arrows, giving a notion of composition on the graph

\item  a property of associativity for composition: we have that
\begin{center}

\begin{tabular}{ccc}
\begin{tabular}{c}
  \def\cartesien{%
    \ar@{-}[]+R+<6pt,-1pt>;[]+RD+<6pt,-6pt>%
    \ar@{-}[]+D+<1pt,-6pt>;[]+RD+<6pt,-6pt>%
  }
  \xymatrix{ 
    X_{3} \ar[r] \ar[d] \cartesien & X_2 \ar[d]^{X_{max}} \\
    X_1 \ar[r]_{X_{min}} & X_0 }\end{tabular}
& \ and \ &
\begin{tabular}{c}
  \def\cartesien{%
    \ar@{-}[]+R+<6pt,-1pt>;[]+RD+<6pt,-6pt>%
    \ar@{-}[]+D+<1pt,-6pt>;[]+RD+<6pt,-6pt>%
  }
  \xymatrix{ 
    X_{3} \ar[r] \ar[d] \cartesien & X_1 \ar[d]^{X_{max}} \\
    X_2 \ar[r]_{X_{min}} & X_0 }\end{tabular} \end{tabular}    
\end{center}
are both pullbacks in $Set$, so:\begin{center}

$X_3 \cong \{(A \rightarrow B,B \rightarrow C \rightarrow D)\}$\end{center}
and: \begin{center}
$X_3 \cong \{(A \rightarrow B \rightarrow C,C \rightarrow D)\}$\end{center}

Transitivity of $\cong$ implies that these sets are isomorphic: to every composition $h \circ (g \circ f)$ corresponds one and only composition $(h \circ g) \circ f$, so that a simplicial set satisfying the Segal condition has an associativity property.

\item an identity over every object of the category: it comes from the very structure of simplicial sets since the degeneracy functions $[1] \rightarrow [0]$ induce on the simplicial set $X$ maps $X_{[0]} \rightarrow X_{[1]}$ which should be thought of as identities.

\end{itemize}
\end{remark}

\subsubsection{Segal condition as a representability property} 
\label{segalrep}
In fact there is another way to see the Segal condition. We first consider the category $\Delta_0$ whose objects are the finite ordinals $[n]$ (so they are exactly the same as in $\Delta$) and whose morphisms are the distance-preserving functions: \mbox{$f\,:\,[m]\rightarrow [n]$} is an arrow in $\Delta_0$ iff $\forall p \in [m-1],\ f(p+1)\,=\,f(p)+1$. $\Delta_0$ is a subcategory of $\Delta$ with the same objects, but also a full subcategory of $Graph$ via the inclusion functor $i_0$. We thus have the following commutative diagram:
\begin{center}
\begin{tabular}{c}
  \xymatrix{
    \Delta \ar[r]^i & Cat \\
    \Delta_0 \ar[r]_{i_0} \ar[u]^{l}& Graph \ar[u]_{Free} \\
  }
\end{tabular}  
\end{center}

where $Free$ takes an oriented graph to its free category and $l$ includes $\Delta_0$ in $\Delta$ (and thus is an identity-on-objects functor). Using $\Delta_0$, we have the following formulation of the Segal condition:

\begin{theorem}[Segal condition as a representability condition]
A simplicial set $X$ is isomorphic to the nerve of a category if and only if there exists a graph $G$ such that the functor:
\begin{center}
$\Delta_0^{op} \xrightarrow{l^{op}} \Delta^{op} \xrightarrow{X} Set$
\end{center} 
is isomorphic to the functor:
\begin{center}
$Graph(i_0,G)\,:\,[n] \mapsto Graph(i_0 [n],G)$
\end{center}
\end{theorem}

It is interesting to see that the construction of $Graph(i_0,G)$ derives from the one of $Graph(i_0,1)$ (mapping a graph $G$ to the corresponding presheaf $Graph(i_0,G)$), which in turn is exactly similar to the construction of the categorical nerve $Cat(i,1)$ mentioned above. In fact, it is also a kind of nerve, as we will see in Section \ref{yoneda}. So, the point is that a simplicial set is isomorphic to the nerve of a category iff its restriction along $l$ is isomorphic to the "nerve" of a given graph: this meaning that the paths of length $n$ described by the simplicial set are isomorphic to the paths of length $n$ of the graph $G$.\\
It is quite clear why such a simplicial set can be realized as a category: as in the former Segal condition, the existence of identities arise from the very structure of the simplicial set, and the two other properties arise from the fact that $n$-simplices have to be isomorphic to the paths of length $n$ of some graph: so, composition of two composable arrows is given by the simplex corresponding to the path of length $2$ generated by composing the two paths of length $1$ corresponding to the two arrows (seen as $1$-simplices), and associativity similarly comes from the associativity of the composition of paths in a graph.

\subsection{Kan extensions}
In our quest for a generalization of the notion of arities, we are not only looking for a formalization of the notion of \emph{density}, that is, of a kind of subcategory which gives all the pieces required to build a bigger one, but also for a proper notion of \emph{glueing} of these pieces. Kan extensions will turn out to be the good concept for such a reconstruction of a category from its subcategory of arities. Before we introduce them, we require the notions of coend and weighted colimit, which extends the idea of geometric realization of a simplicial set (taking copies of simplices and glueing them according to the face inclusions). 

\subsubsection{Coends} When considering functors $S,T\,:\,\mathcal{C}^{op}\times \mathcal{C} \rightarrow \mathcal{D}$, the notion of natural transformation between them is not relevant, since the appearance of $\mathcal{C}$ both in a contravariant and in a covariant position conflicts with the definition of such a transformation. We then have to extend the notion, and this should lead us to the following equivalence appearing in \cite[Section 3.2]{curien}:
\begin{center}
\begin{tabular}{ccc}
natural transformation &\ \ \ & dinatural transformation\\
cone & & wedge\\
end & & limit\\
coend & & colimit\\
\end{tabular}
\end{center}

\begin{example}[Functors mixing covariant and contravariant behaviour]
An example of functor mixing covariant and contravariant behaviour is the one of evaluation when an object $B$ of a cartesian closed category is fixed: in this case, $eval\,:\,A \rightarrow Hom(A,B) \times A  \rightarrow B$ has both behaviours.\\
Another example is the notion of composition in a category: the composition $Hom(B,C) \times Hom(A,B) \rightarrow Hom(A,C)$ is natural (contravariantly) in $a$, natural (covariantly) in $c$, and dinatural in $b$.
\end{example}

\begin{definition}[Dinatural transformation]
Given functors $S,T\,:\,\mathcal{C}^{op}\times \mathcal{C} \rightarrow \mathcal{D}$, a dinatural transformation $\alpha\,:\,S \rightarrow T$ is a function $\alpha$ which assigns to each object $c \in \mathcal{C}$ an arrow $\alpha_C\,:\,S(c,c)\rightarrow T(c,c)$ of $D$ in a way that makes the following diagram commute for each arrow $f\,:\,c\rightarrow c'$ in $\mathcal{C}$:
\begin{center}
\begin{tabular}{c}
\xymatrix{
\ & S(c,c) \ar[r]^{\alpha_c} & T(c,c) \ar[dr]^{T(1,f)} & \ \\
S(c',c) \ar[ur]^{S(f,1)} \ar[dr]_{S(1,f)} & \ & \ & T(c,c')\\
\ & S(c',c') \ar[r]_{\alpha_{c'}} & T(c',c') \ar[ur]_{T(f,1)} & \ \\
}

\end{tabular}
\end{center}
\end{definition}

\begin{remark}[Natural transformations yield dinatural transformations] It is easy to remark that the natural transformations $S\rightarrow T$, with $S,T\,:\,\mathcal{C}^{op}\rightarrow \mathcal{D}$, are in one-to-one correspondence with the dinatural transformations $\lambda(x,y).Fx \rightarrow \lambda(x,y).Gx$ (where the $\lambda$-calculus-like notation denotes exponentiation in the usual categorical sense).
\end{remark}

\begin{remark}[\emph{Splendeur et mis\`{e}re} of the dinatural transformations]
A serious limitation of the dinatural transformations is that they do not compose vertically. However, heterogeneous compositions are possible: if $\mu\,:\,F\rightarrow G$ is a natural transformation and \mbox{$\nu\,:\lambda(x,y).Gx\rightarrow S$} is a dinatural transformation, the one-to-one correspondence gives a transformation $\lambda(x,y).Fx\rightarrow \lambda(x,y).Gx$ whose composition with $\nu$ gives what deserve the name of $\nu \circ \mu\,:\,\lambda(x,y).Fx\rightarrow S$.
\end{remark}

Recall that limits in category theory can be formalized with a notion of universal cone. The generalization in our case is the notion of end, that is, of universal wedge. Similarly to the ordinary case, a wedge is a dinatural transformation between a constant functor and another one. When it is from the constant functor to the general one, it generalizes cones; dually, the other direction generalizes cocones.

\begin{definition}[Wedge]
If $T\,=\,b\,:\,\mathcal{C}^{op}\times \mathcal{C} \rightarrow \mathcal{D}$ is dummy in both variables, a dinatural transformation $S\rightarrow T$ consists of components $\alpha_c\,:\,S(c,c)\rightarrow b$ which makes the following diagram commute for every $f\,:\,c\rightarrow c'$:
\begin{center}
\begin{tabular}{c}
\xymatrix{
S(c',c) \ar[r]^{S(1,f)} \ar[d]_{S(f,1)} & S(c',c') \ar[d]^{\alpha_{c'}} \\
S(c,c) \ar[r]_{\alpha_c} & b\\
}
\end{tabular}
\end{center}
Such a transformation is called a \emph{wedge} from $S$ to $b$. Dually, there is a notion of wedge from $b$ to $S$, satisfying:
\begin{center}
\begin{tabular}{c}
\xymatrix{
b \ar[r]^{\alpha_c} \ar[d]_{\alpha_{c'}} & T(c,c) \ar[d]^{T(1,f)} \\
T(c',c') \ar[r]_{T(f,1)} & T(c,c')\\
}
\end{tabular}
\end{center}
\end{definition}

We can now introduce ends and coends, which are the extension of limits and colimits, that is, \emph{universal} generalized cones/cocones:

\begin{definition}[End]
An end of a functor $S\,:\,\mathcal{C}^{op}\times \mathcal{C}\rightarrow \mathcal{D}$ is a universal dinatural transformation $\omega$ from a constant\footnote{That is, a functor $\mathcal{C}^{op}\times \mathcal{C}\rightarrow \mathcal{D}$ evaluating everywhere to $e$ and sending every arrow to $1_e$.} $e\in D$ to $S$. This means that, for every dinatural transformation $\beta$ from a constant $x\in X$, there is a \emph{unique} arrow $h\,:x\rightarrow e$ in $\mathcal{D}$ satisfying, for all $a \in C$, $\beta_a\,=\,\omega_a h$. Diagrammatically, this means that for each arrow $f\,:b\rightarrow c$ of $\mathcal{C}$ we have:
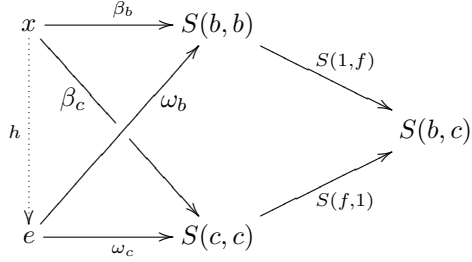
\begin{center}
\begin{tabular}{c}
\xymatrix @R=2pc @C=4pc {
x \ar[ddr]^*++{\ \, \omega_b} |!{[dd];[r]}\hole \ar@{.>}[dd]_h \ar[r]^{\beta_b} & S(b,b) \ar[dr]^{S(1,f)} & \ \\
\ & \ & S(b,c)\\
e \ar[r]_{\omega_c} \ar[uur]^*++{\beta_c \ \,} & S(c,c) \ar[ur]_{S(f,1)} \\
}
\end{tabular}
\end{center}
where the two quadrilaterals commute, and where the universal property of $\omega$ amounts to the existence and uniqueness of such an arrow $h$ making the two triangles on the left commute.

\end{definition}

The constant $e$ is traditionally denoted by $\int_c{S(c,c)}$. A given end (that is, a constant object and an universal dinatural transformation from it) is unique up to isomorphism.

\begin{definition}[Coend]
A coend of a functor $S\,:\,\mathcal{C}^{op}\times \mathcal{C}\rightarrow \mathcal{D}$ is a dinatural transformation $\alpha$ from $S$ to a constant $x \in X$ which is universal among dinatural transformations from $S$ to a constant. $x$ is traditionally denoted by $\int^c{S(c,c)}$.
\end{definition}

\begin{example}[A usual coend] Here we show that the usual tensor product of two modules over a ring is obtained from a coend. We need first to recall the categorical formalization of rings and modules on them:
\begin{itemize}
\item An $Ab$-category\footnote{This is a simple case of \emph{enrichment} of a category. A quick introduction may be found at \cite[Section 1.3]{leinster}, a reference being \cite{kelly}.} is a category whose hom-sets are additive abelian groups, and for which composition is bilinear (that is: if $f,f'\,:\,a\rightarrow b$ and $g,g'\,:\,b\rightarrow c$, $(g+g')\circ (f+f')\,=\,g\circ f + g\circ f' + g' \circ f + g' \circ f'$),
\item A functor $T$ between $Ab$-category is \emph{additive} when it respects the additive structure: if $f$ and $f'$ are parallel arrows from the source category, this means that $T(f+_s f')\,=\,T(f)+_t T(f')$ (where we distinguished $+$ in the source and target category),
\item A ring $R$ is an $Ab$-category with only one object, with arrows its elements, and where composition realizes the product in $R$\footnote{Note that this is a monoid (in the usual categorical sense) in an $Ab$-enriched category.},
\item A left $R$-module $B$ is an additive functor $R \rightarrow Ab$ sending the only object of $R$ to the abelian group $B$, and each arrow in $R$ to the scalar multiplication $r_*\,:\,b\mapsto rb$ in $B$,
\item A right $R$-module is an additive functor $R^{op}\rightarrow Ab$ built in the same way.
\end{itemize}
Given a left $R$-module $A$ and a right one $B$, the usual tensor product $\otimes$ in $Ab$ induces a bifunctor $R^{op}\times R \rightarrow Ab$ (where the image of the single object of $R^{op}\times R$ is $ A \otimes B$). The following coend:
\[
\int^R{A \otimes B}\,=\,A\otimes_R B
\]
turns to be the usual tensor product of $A$ and $B$ over $R$: a wedge \mbox{$\alpha\,:\, A  \otimes B \rightarrow M$} with $M \in Ab$ is precisely a single morphism $\rho\,:\,A\otimes B \rightarrow M$ of abelian groups making the following diagram commute for every arrow $r\in R$:
\begin{center}
\begin{tabular}{c}
\xymatrix @R=4pc @C=4pc {
A \otimes B \ar[r]^{1_A \otimes r_*}  \ar[d]_{r_* \otimes 1_B} & A \otimes B \ar[d]^{\rho}\\
A \otimes B \ar[r]_{\rho} & M\\
}
\end{tabular}
\end{center}
With the interpretation of modules as functors, this means that for every $a\in A$ and $b \in B$ we have $\rho (ar \otimes b)\,=\,\rho (a \otimes rb)$. So, $M$ is a coend iff $M\,=\,A\otimes B /\sim$, where the equivalence relation $\sim$ is the smallest one identifying $ar\otimes b$ and $a\otimes rb$ for every couple $(a,b)$: as claimed, this defines the usual tensor product of $A$ and $B$ over $R$.
\end{example}

\begin{remark}[Monoidal categories have a tensor product]
If $\mathcal{C}$ is a monoidal category, with multiplication $\bullet$, every couple of functors $T\,:\,\mathcal{D}^{op}\rightarrow \mathcal{C}$ and $S\,:\,\mathcal{D}\rightarrow \mathcal{C}$ has a tensor product - which is an element of $\mathcal{C}$ - given by the following coend:
\[
T \bullet_{\mathcal{D}} S\,=\,\int^{\mathcal{D}}{(Td)\bullet(Sd)}
\]
\end{remark}

So, the notion of coend provides some way of $\emph{glueing}$ stuff together. We are thus getting closer to our idea of the realization of a nerve - we now need to mix this glueing operation with one of duplication of elementary simplices; in our case copowers should be enough. This leads to the notion of weighted colimits.

\begin{remark}[Link with the usual notion of limit]
In \cite[Section IX.5]{maclane} and \cite[Exercice IX.6.3]{maclane}, two categories are defined where the concepts of end and coend are expressed using ordinary category theory.
\end{remark}

\subsubsection{Geometric realization of a simplicial set}
The idea of glueing standard simplices together to build a topological space from a simplicial set can be formalized using copowers and coends. Take a simplicial set $S$ - remember it's a presheaf over $\Delta$, that is, a functor $\Delta^{op}\rightarrow Set$ - and consider the functor $\zeta\,:\,\Delta\rightarrow Top$ which sends every ordinal $[n]$ to the standard $n$-simplex and arrows to the obvious ones (face inclusions to topological face inclusions, degeneracies to topological degeneracies, and every arrow to the composite of the images of its decomposition by the factorization property of arrows in $\Delta$). Now we introduce a shortcut\footnote{This is the notion of \emph{copower} - but here this intuition should be enough.} for the iterated coproduct in $Top$: for $S'$ a set and $X$ a topological set, we denote by $S'\cdot X$ the coproduct in $Top$ of as many copies of $X$ as there are elements in $S'$. It should now be clear that to realize the simplicial set in $Top$ we need, for every arity $[n]$, as many $n$-simplices as there are elements in $S[n]$. We also need to paste them according to the face and degeneracy operations, which are induced in $S$ by the arrows of $\Delta$: $S$ has the following realization in $Top$ (where the use of a coend is legitimated by the fact that $(n,m)\rightarrow S[n] \cdot \zeta [m]$ is a functor $\Delta^{op}\times \Delta \rightarrow Top$):
\begin{center}
$|S|\,=\,\int^{n \in \mathbb{N}}{S[n]\cdot \zeta [n]}$
\end{center}
\ \\
This idea of taking some kind of colimit of an usual diagram but with some entries duplicated is formalized by the notion of weighted colimit.

\subsubsection{Weighted colimits}
\label{weight}
Here we follow \cite{baezcolim}\footnote{It seemed us interesting to see the concept of weighted colimit as a \emph{categorification} of the usual notion of integral instead of just giving its definition.}, in the case where the category is Set-enriched (\textit{i.e.} the usual case). The usual Lebesgue integral associated to a measure $w\,:\,X\rightarrow \mathbb{C}$ with $X$ a finite set is basically a map giving a complex "weight" to each of the elements of $X$. Then any function from $X$ to a vector space $L$ can be integrated as usual:\\
\begin{center}
$\int_X{fdw}\,=\, \sum_{x\in X}{w(x)f(x)}$
\end{center}
Weighted colimits can be seen as a categorification of this idea. One may have remarked that in our purpose of generalizing the nerve construction and its dual which is realization we precisely need to be able to relate in some way potentially \emph{many} copies of a given standard simplex: here is what we need weights for. Now, $X$ is not a set anymore but a category. But we take a restriction here to the usual notion of weighted colimit: we only consider weights with values in $Set$, so that we take for weight function $w\,;\,X^{op}\rightarrow Set$, and the function we "measure" is now a functor $F\,:\,X\rightarrow L$. Since the image of $w$ is in $Set$, for a given object $x \in X$, $w(x)\cdot F(x)$ has a sense as the coproduct of $\operatorname{Card}(w(x))$ copies of $F(x)$ (if it exists\footnote{It especially is the case when $L$ is locally small and has all coproducts.}), generalizing what we did before. So, $X^{op}\times X\xrightarrow{w\times F} Set \times L \xrightarrow{\cdot} L$ is a functor which can be understood as a weighted version of $F$, and whose colimit can be taken with a coend, since coends provide the notion of colimit for such "mixed behaviour" functors.\\\ \\
We can now formalize the idea of density we previously sketched.

\subsubsection{Density}
\label{density}
\begin{definition}[Density]
\label{defdensity}
A functor $i\,:\,\mathcal{C}\rightarrow \mathcal{D}$ (with $\mathcal{D}$ locally small) is \emph{dense} when it satisfies one of the two following equivalent conditions:
\begin{itemize}
\item The induced nerve functor (see Section \ref{yoneda}):
\begin{center} $
\begin{array}{llcc}
\mathcal{D}(i,1)\,:\, & \mathcal{D} & \rightarrow & \widehat{\mathcal{C}}\\
& D & \mapsto & 
\bigg\{ \begin{array}{lcr}
C & \mapsto & \mathcal{D}(iC,D)\\
f & \mapsto & (if)^{*}\\
\end{array}\\ \end{array} $
\end{center}
is fully faithful
\item Every object $d \in \mathcal{D}$ is obtained as the colimit of $i$ weighted by the functor $c\mapsto \mathcal{D}(ic,d))$
\end{itemize}
\end{definition}

Density formalizes the condition under which every object of a category may be obtained as a weighted colimit of objects of a given subcategory - this generalizes the idea of the reconstruction of categories from elementary simplices, of which we take several copies properly glued together. The following fact should not be surprising since it has been our guideline: $\Delta$ is dense in $Cat$\footnote{As expected, $\mathcal{C}$ is said to be dense in $\mathcal{D}$ when the inclusion functor $i\,:\,\mathcal{C}\rightarrow \mathcal{D}$ is dense.}. 

Next step\footnote{According to \cite[Historical note, p. 6]{diwo}, this is how Kan extensions were historically introduced.} towards the generalization of the realization construction comes with the introduction of the Kan extension.

\subsubsection{Kan extensions} The notion of Kan extension is purely 2-categorical, and we assume some familiarity reader with $2$-categories - if needed, one may consult \cite[Section XII.3]{maclane} (for example) for an introduction. Kan extensions are closely related to the change-of-base operation, as explained in a remark to follow: the idea is to find the best approximation of a functor when it is transported to another basis by precomposition. There is a strong link with what we introduced before in the section, as we shall see soon. We start defining Kan extensions by their universal property:

\begin{definition}[Kan extensions]
Let $F\,:\,\mathcal{C}\rightarrow \mathcal{D}$ and $i\,:\,\mathcal{C}\rightarrow \mathcal{C'}$ be functors. The left Kan extension of $F$ along $i$ is given by a functor $\operatorname{Lan}_i F$ and a natural transformation $\eta_F\,:\,F\Rightarrow \operatorname{Lan}_i F \circ i$:

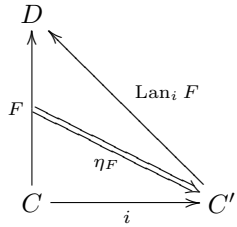
\begin{center}
\begin{tabular}{c}
  \def\enlargexyentry#1{%
    \POS "#1",*{
      \vrule height 1pt depth 1pt width 0pt
      \vrule height 0pt depth 0pt width 2pt
      }="#1",}
  \xymatrix @!0 @R=25mm @C=25mm  { 
    D & \  \\
    C \ar[u]|*{}="A"^F \ar[r]_i  
     \enlargexyentry A
    \ar@2{<-}[r];"A"^-{\eta_F} 
    & C'     \ar[lu]_{\operatorname{Lan}_i F}  
   \\}
\end{tabular}
\end{center}
satisfying the following universal property: every 2-cell:

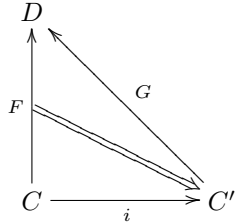
\begin{center}
\begin{tabular}{c}
  \def\enlargexyentry#1{%
    \POS "#1",*{
      \vrule height 1pt depth 1pt width 0pt
      \vrule height 0pt depth 0pt width 2pt
      }="#1",}
  \xymatrix @!0 @R=25mm @C=25mm  { 
    D & \  \\
    C \ar[u]|*{}="A"^F \ar[r]_i  
     \enlargexyentry A
    \ar@2{<-}[r];"A"^-{\ } 
    & C'     \ar[lu]_{G}  
   \\}
\end{tabular}
\end{center}

factors through the Kan extension:

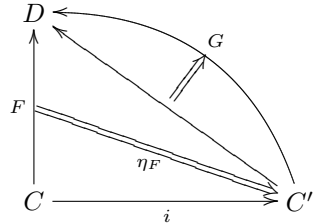
\begin{center}
\begin{tabular}{c}
  \def\enlargexyentry#1{%
    \POS "#1",*{
      \vrule height 1pt depth 1pt width 0pt
      \vrule height 0pt depth 0pt width 2pt
      }="#1",}
  \xymatrix @!0 @R=25mm @C=35mm  { 
    D & \  \\
    C \ar[u]|*{}="A"^F \ar[r]_i  
     \enlargexyentry A
    \ar@2{<-}[r];"A"^-{\eta_F} 
    & C'   \ar[lu]_{}="B"     \ar@/_2pc/[lu]_{G}="C"
    \ar@{=>}"B";"C"  
   \\}
\end{tabular}
\end{center}

which means in a more traditional language that:
\begin{center} 
$\operatorname{Nat}(\operatorname{Lan}_i F,G)\,\cong\,\operatorname{Nat}(F,G\circ i)$.
\end{center}

Dually, the right Kan extension of $F$ along $i$ is given by a functor $\operatorname{Ran}_i F$ and a natural transformation $\epsilon_F\,:\,\operatorname{Ran}_i F \circ i\Rightarrow F$:

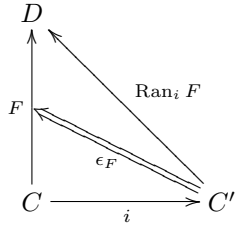
\begin{center}
\begin{tabular}{c}
  \def\enlargexyentry#1{%
    \POS "#1",*{
      \vrule height 1pt depth 1pt width 0pt
      \vrule height 0pt depth 0pt width 2pt
      }="#1",}
  \xymatrix @!0 @R=25mm @C=25mm  { 
    D & \  \\
    C \ar[u]|*{}="A"^F \ar[r]_i  
     \enlargexyentry A
    \ar@2{->}[r];"A"^-{\epsilon_F} 
    & C'     \ar[lu]_{\operatorname{Ran}_i F}  
   \\}
\end{tabular}
\end{center}
satisfying the following universal property: every 2-cell:

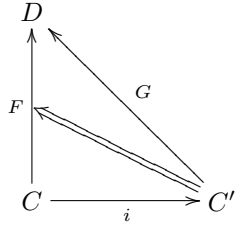
\begin{center}
\begin{tabular}{c}
  \def\enlargexyentry#1{%
    \POS "#1",*{
      \vrule height 1pt depth 1pt width 0pt
      \vrule height 0pt depth 0pt width 2pt
      }="#1",}
  \xymatrix @!0 @R=25mm @C=25mm  { 
    D & \  \\
    C \ar[u]|*{}="A"^F \ar[r]_i  
     \enlargexyentry A
    \ar@2{->}[r];"A"^-{\ } 
    & C'     \ar[lu]_{G}  
   \\}
\end{tabular}
\end{center}

factors through the Kan extension:

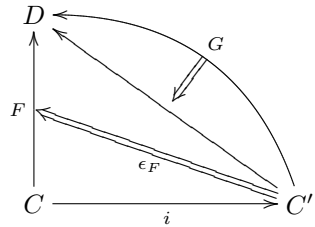
\begin{center}
\begin{tabular}{c}
  \def\enlargexyentry#1{%
    \POS "#1",*{
      \vrule height 1pt depth 1pt width 0pt
      \vrule height 0pt depth 0pt width 2pt
      }="#1",}
  \xymatrix @!0 @R=25mm @C=35mm  { 
    D & \  \\
    C \ar[u]|*{}="A"^F \ar[r]_i  
     \enlargexyentry A
    \ar@2{->}[r];"A"^-{\epsilon_F} 
    & C'   \ar[lu]_{}="B"     \ar@/_2pc/[lu]_{G}="C"
    \ar@2{<-}"B";"C"  
   \\}
\end{tabular}
\end{center}

which means in a more traditional language that:
\begin{center} 
$\operatorname{Nat}(G,\operatorname{Ran}_i F)\,\cong\,\operatorname{Nat}(G\circ i,F)$.
\end{center}
\end{definition}

The reader unaware of $2$-categories may feel a bit confused: these diagrams are \emph{not} commutative diagrams but are $2$-cells, that is, representation of natural transformations between functors.\\ 
A first interesting theorem is the following:
\begin{theorem}[Existence of the left Kan extension]
Let $F\,:\,\mathcal{C}\rightarrow \mathcal{D}$ and \mbox{$i\,:\,\mathcal{C}\rightarrow \mathcal{C'}$} be functors. If $\mathcal{C}$ is small and $\mathcal{D}$ is cocomplete, a left Kan extension of $F$ along $i$ exists. Moreover, if $i$ is fully faithful, $F\,\cong\,\operatorname{Lan}_i F \circ i$ - equivalently, the natural transformation giving the universal $2$-cell is an isomorphism.
\end{theorem}

The connection with nerves and realization arises from the following theorem\footnote{See \cite[Section X.4]{maclane}}:
\begin{theorem}[The coend formula for the left Kan extension\footnote{The reader may deduce from our previous discussion that this formula is more about a weighted colimit than about a coend.}]
Given functors $F\,:\,\mathcal{C}\rightarrow \mathcal{D}$ and $i\,:\,\mathcal{C}\rightarrow \mathcal{E}$ such that for all $(c,c',e)\in \mathcal{C}^2\times \mathcal{E}$ the copowers $\mathcal{E}(ic',e)\cdot Fc$ (that is, the coproduct of as many copies of $Fc$ as there are elements in $\mathcal{E}(ic',e)$) exist in $\mathcal{D}$, $\operatorname{Lan}_i F$ exists if (and only if) the following coends, giving then its value on objects, exist for every $e\in \mathcal{E}$:
\begin{center}
$(\operatorname{Lan}_i F)(e)\,=\,\int^c{\mathcal{E}(ic,e)\cdot Fc}$
\end{center}
\end{theorem}

Specifically, the condition on copowers is fulfilled when $\mathcal{D}$ is locally small and has all coproducts. An important point is that when dealing with $Set$ the copower meets the usual notion of cartesian product of sets: if $X$ and $Y$ are sets, $X \cdot Y\,\cong \,X\times Y$.

\paragraph{Realization as a Kan extension}\label{realkan} If we denote by $Y_{\Delta}$ the usual Yoneda embedding $\Delta\rightarrow \widehat{\Delta}$ and by $i\,:\,\Delta\rightarrow Cat$ the inclusion we introduced before, we have that:
\begin{center}
\begin{tabular}{c}
  \def\enlargexyentry#1{%
    \POS "#1",*{
      \vrule height 1pt depth 1pt width 0pt
      \vrule height 0pt depth 0pt width 2pt
      }="#1",}
  \xymatrix @!0 @R=25mm @C=25mm  { 
    Cat & \  \\
    \Delta \ar[u]|*{}="A"^i \ar[r]_{Y_{\Delta}}  
     \enlargexyentry A
    \ar@2{<-}[r];"A"^-{\eta_{Y_{\Delta}}} 
    & \widehat{\Delta}     \ar[lu]_{\operatorname{Lan}_{Y_{\Delta}} i}  
   \\}
\end{tabular}
\end{center}
is a left Kan extension. The coend formula gives on every simplical set $S$:
\begin{center}
\begin{tabular}{rcl}
$(\operatorname{Lan}_{Y_{\Delta}} i) S$ & $\cong$ & $\int^{[n]}{\widehat{\Delta}(Y_{\Delta}[n],S)\cdot i [n]}$\\
& $\cong$ & $\int^{[n]}{S[n]\cdot i [n]}$\\
\end{tabular}
\end{center}
which is the weighted colimit corresponding to realization\footnote{Here in $Cat$ instead of $Top$.}: $\operatorname{Lan}_{Y_{\Delta}} i$ realizes a simplicial set in $Cat$, and corresponds to the left Kan extension of the embedding functor in the category of realization $i$ along the Yoneda embedding. Note that, $i$ being fully faithful, $\eta_{Y_{\Delta}}$ is an isomorphism.

\paragraph{Nerve as a Kan extension} With the same notations, we have that:
\begin{center}
\begin{tabular}{c}
  \def\enlargexyentry#1{%
    \POS "#1",*{
      \vrule height 1pt depth 1pt width 0pt
      \vrule height 0pt depth 0pt width 2pt
      }="#1",}
  \xymatrix @!0 @R=25mm @C=25mm  { 
    \widehat{\Delta} & \  \\
    \Delta \ar[u]|*{}="A"^{Y_{\Delta}} \ar[r]_i  
     \enlargexyentry A
    \ar@2{<-}[r];"A"^-{\eta_i} 
    & Cat     \ar[lu]_{\operatorname{Lan}_{i} Y_{\Delta}\,=\,Cat(i,1)}  
   \\}
\end{tabular}
\end{center}
is a left Kan extension  (so that the nerve operation can be defined as the left Kan extension of the Yoneda embedding along the inclusion functor $i$). As the reader may expect it, there is a link between  the nerve and realization operations: the realization operation is the left adjoint of the associated nerve. $2$-categorically, this gives a bijection between $2$-cells:

\begin{center}
\begin{tabular}{ccc}
\begin{tabular}{c}
  \def\enlargexyentry#1{%
    \POS "#1",*{
      \vrule height 1pt depth 1pt width 0pt
      \vrule height 0pt depth 0pt width 2pt
      }="#1",}
  \xymatrix @!0 @R=25mm @C=25mm  { 
    \widehat{\Delta} & \  \\
    \Delta \ar[u]|*{}="A"^{Y_{\Delta}} \ar[r]_i  
     \enlargexyentry A
    \ar@2{<-}[r];"A"^-{\eta_i} 
    & Cat     \ar[lu]_{Cat(i,1)}  
   \\}
\end{tabular}
&\ $\cong$ \   &
\begin{tabular}{c}
  \def\enlargexyentry#1{%
    \POS "#1",*{
      \vrule height 1pt depth 1pt width 0pt
      \vrule height 0pt depth 0pt width 2pt
      }="#1",}
  \xymatrix @!0 @R=25mm @C=25mm  { 
    \widehat{\Delta} & \  \\
    \Delta \ar[u]|*{}="A"^{Y_{\Delta}} \ar[r]_i  
     \enlargexyentry A
    \ar@2{<-}[r];"A"^-{\eta} 
    & Cat     \ar@{<-}[lu]_{nerve}  
   \\}
\end{tabular}\\
\end{tabular}
\end{center}

where the natural transformation in the $2$-cell is an  isomorphism in the left diagram (since $i$ is fully faithful) and thus also in the right one: so we can write it in the other direction:

\begin{center}
\begin{tabular}{ccc}
\begin{tabular}{c}
  \def\enlargexyentry#1{%
    \POS "#1",*{
      \vrule height 1pt depth 1pt width 0pt
      \vrule height 0pt depth 0pt width 2pt
      }="#1",}
  \xymatrix @!0 @R=25mm @C=25mm  { 
    \widehat{\Delta} & \  \\
    \Delta \ar[u]|*{}="A"^{Y_{\Delta}} \ar[r]_i  
     \enlargexyentry A
    \ar@2{<-}[r];"A"^-{\eta_i} 
    & Cat     \ar[lu]_{Cat(i,1)}  
   \\}
\end{tabular}
&\ $\cong$ \   &
\begin{tabular}{c}
  \def\enlargexyentry#1{%
    \POS "#1",*{
      \vrule height 1pt depth 1pt width 0pt
      \vrule height 0pt depth 0pt width 2pt
      }="#1",}
  \xymatrix @!0 @R=25mm @C=25mm  { 
    \widehat{\Delta} & \  \\
    \Delta \ar[u]|*{}="A"^{Y_{\Delta}} \ar[r]_i  
     \enlargexyentry A
    \ar@2[r];"A"^-{\eta^{-1}} 
    & Cat     \ar@{<-}[lu]_{nerve}  
   \\}
\end{tabular}\\
\end{tabular}
\end{center}

and we now have on the right, even if the diagram is rotated in a way that doesn't meet our usual convention, a $2$-cell looking like a left Kan extension. It happens to be one, since the universal property is just about cuting and pasting $2$-cells: our claim that the nerve is a left Kan extension was righteous\footnote{Actually, it is an axiom in Yoneda structures (see Section \ref{yoneda}), and the fact that the realization is a left Kan extension arises from the converse reasonment.}.

\paragraph{Kan extensions and (co)limits}
The reader should be quickly convinced from the previous discussion that in the case where the functor the Kan extension is taken along is $i\,:\,\mathcal{C}\rightarrow \{*\}$, we have that
\begin{itemize}
\item a functor $F$ has a limit iff its right extension along $i$ exists, and then the value of the limit is the value of $\operatorname{Ran}_i F$ on the only object $*$ of the target category for $i$,
\item a functor $F$ has a colimit iff its left extension along $i$ exists, and then the value of the colimit is the value of $\operatorname{Lan}_i F$ on the only object $*$ of the target category for $i$
\end{itemize}
It is quite obvious since in the second case the coend formula exhibits the left Kan extension as a weighted colimit with uniform weight $1$ on all objects: the dinatural transformation is then dummy in its contravariant parameter and basically is a natural transformation, giving an usual cocone which is universal. The first case can be treated dually.\\
In fact all fundamental concepts of category theory may be expressed with Kan extensions, as detailed in \cite[Section X.7]{maclane}. Other interesting examples of reformulations of usual concepts with Kan extensions follow. First, we say that a given right\footnote{As usual with Kan extensions, we often only treat the case immediatly useful for our purpose, the other one being dual.} Kan extension $(\operatorname{Ran}_i F, \epsilon)$ is preserved by the functor $G\,:\mathcal{D}\rightarrow \mathcal{D'}$ if $(G\operatorname{Ran}_i F, G\epsilon)$ is a right Kan extension of $G\circ F$ along $i$. For $i\,:\,\mathcal{C}\rightarrow \mathcal{D}$ and $j\,:\,\mathcal{D}\rightarrow \mathcal{C}$, the following properties are equivalent:
\begin{itemize}
\item $j \dashv i$, with $\epsilon$ as counit of the adjunction,
\item $(j,\epsilon)$ is a right Kan extension of $id$ along $i$ that is preserved by all functors,
\item $(j,\epsilon)$ is a right Kan extension of $id$ along $i$ that is preserved by $i$.
\end{itemize}
Two other interesting properties deal with the notion of density:
\begin{itemize}
\item A functor $i\,:\,\mathcal{A} \rightarrow \mathcal{C}$ is \emph{dense} iff:
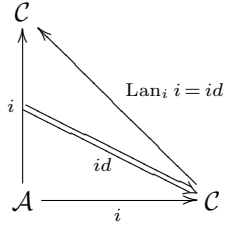
\begin{center}
\begin{tabular}{c}
  \def\enlargexyentry#1{%
    \POS "#1",*{
      \vrule height 1pt depth 1pt width 0pt
      \vrule height 0pt depth 0pt width 2pt
      }="#1",}
  \xymatrix @!0 @R=25mm @C=25mm  { 
    \mathcal{C} & \  \\
    \mathcal{A} \ar[u]|*{}="A"^i \ar[r]_{i}  
     \enlargexyentry A
    \ar@2{<-}[r];"A"^-{id} 
    & \mathcal{C}     \ar[lu]_{\operatorname{Lan}_{i} i\,=\,id}  
   \\}
\end{tabular}
\end{center}
is a left Kan extension of $i$ along $i$. This expresses, through the coend formula for Kan extensions, exactly the idea that a dense subcategory is like a box of elementary pieces of which we can glue together as many copies as needed to built the whole category.

\item When considering a functor $T\,:\,\mathcal{C}\rightarrow \mathcal{D}$, it has a left Kan extension along the Yoneda embedding $Y_{\mathcal{C}}\,:\,\mathcal{C}\rightarrow \widehat{\mathcal{C}}$ iff the functor $\mathcal{D}(F,1)\,:\,D \mapsto \mathcal{D}(F\cdot,D)$ has a left adjoint; and then this left adjoint is precisely $\operatorname{Lan}_{Y_{\mathcal{C}}} F$.
\end{itemize}

\begin{remark}[Kan extensions and (pre)sheaves] Historically, the idea behind Kan extensions is the one of change of basis, a very common operation in algebraic geometry which is due to Alexander Grothendieck. An easy and common example of is the one of complexification of a $\mathbb{R}$-module: tensoring it with $\mathbb{C}$ over $\mathbb{R}$, one obtains a ring generated by the same elements but with relations over $\mathbb{C}$ instead of $\mathbb{R}$. The quest for Kan extensions may then be seen as starting from the fact that every functor $F\,:\,\mathcal{C}\rightarrow \mathcal{E}$ induces a functor $F^*\,:\,\widehat{\mathcal{E}}\rightarrow \widehat{\mathcal{C}}$, sending a presheaf $\phi\,:\,\mathcal{E}^{op}\rightarrow Set$ to the presheaf $\phi\circ F^{op}\,:\,\mathcal{C}^{op}\rightarrow Set$. When $\mathcal{C}$ is a small category, $F^*$ has both a right and a left adjoint, transporting presheaves $\phi \in \mathcal{C}$ respectively to their right and left Kan extensions along $F^{op}$.\\
This change-of-basis operation has a well known effect on sheaves (see \cite{schapira1} or \cite{schapira2}): the adjoints to precomposition then give the operation of sheaf extension (right adjoint) and of sheaf inverse image (which is not always a sheaf, usually sheafification is performed on the result). This can be intuited in a quite graphical way: when the extension is a left Kan, we have as universal $2$-cell a natural transformation from a functor to the composite of $F^{op}$, which can be thought of $F^{-1}$ as in the case of $Set$, and of another functor provided by the Kan extension. So, we have an universal way of sending a sheaf to a sheaf composed with "the inverse of F": it is the sheaf inverse image. Dually, the right Kan extension "removes $F^{op}$", extending the sheaf.\\
The problem of finding adjoints to precomposition generalizes to the case of categories of functors. Denote by $\mathcal{D}^{\mathcal{C}}$ the category whose objects are the functors $\mathcal{C} \rightarrow \mathcal{D}$ and arrows natural transformations between them: $F\,:\,\mathcal{C}\rightarrow \mathcal{C'}$ induces $F^*\,:\,\mathcal{D}^{\mathcal{C'}}\rightarrow \mathcal{D}^{\mathcal{C}}$, again by precomposition. The existence of ajoints then gives a notion of Kan extension.\\
But we have to distinguish several "levels" of Kan extensions: when such an ajunction property exists, \emph{all} functors with domain $\mathcal{C}$ have a Kan extension along $F^{op}$.The Kan extensions we used so far are in contrast called "pointwise". 
\end{remark}

\paragraph{Finitary monads, a second formulation}
An important point in our quest for a generalization of arities is the fact that a \emph{finitary} monad can be defined differently that we did: a monad $T\,:\,Set\rightarrow Set$ is finitary when the following diagram:
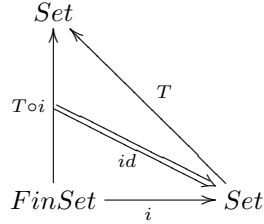
\begin{center}
\begin{tabular}{c}
  \def\enlargexyentry#1{%
    \POS "#1",*{
      \vrule height 1pt depth 1pt width 0pt
      \vrule height 0pt depth 0pt width 2pt
      }="#1",}
  \xymatrix @!0 @R=25mm @C=25mm  { 
    Set & \  \\
    FinSet \ar[u]|*{}="A"^{T\circ i} \ar[r]_{i}  
     \enlargexyentry A
    \ar@2{<-}[r];"A"^-{id} 
    & Set     \ar[lu]_{T}  
   \\}
\end{tabular}
\end{center}
exhibits the functor $T$ together with the identity natural transformation as a left Kan extension of $T\circ i$ along the functor $i$ (which here is just the canonical inclusion $FinSet \rightarrow Set$). This means, by the coend formula, and since the natural transformation associated to this left Kan extension is the identity, that giving the arities (that are: the finite sets) and the values of $T$ on arities is enough to compute $T$ on any set. This was exactly the point of our previous definition of finitary monads, in a case where the coend formula for the left Kan extension gives exactly an inductive colimit: as many copies of each finite ordinal as necessary to represent the set we want to compute $T$ on are taken, and then they are glued together by taking their coend, reconstructing in this way the set according to the total order on it (provided by the arrows between the copies of sets). So this was just a degenerated case of a nerve/realization construction !

\begin{remark}[Free models of an algebraic theory]
The reconstruction \emph{via} left Kan extensions also applies to free models of algebraic theories, as described by Lawvere in \cite{lawvere}\footnote{One also may consult \cite{freemodels}, where the notion of free model is generalized.}. Recall from Section \ref{free} that the forgetful functor $U_j$ sending a model of a Lawvere theory $\mathbb{L}$ to its underlying object in $\mathcal{C}$ has a left adjoint when $\mathcal{C}$ is cartesian closed with small colimits. Moreover, the following $2$-cell where $F$ is a $\mathbb{S}$-model in $\mathcal{C}$ and $j$ the canonical morphism (since $\mathbb{S}$ is initial) is a Kan extension (under these conditions on $\mathcal{C}$):
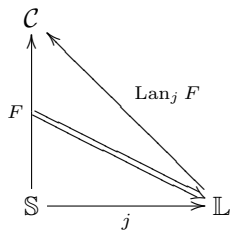
\begin{center}
\begin{tabular}{c}
  \def\enlargexyentry#1{%
    \POS "#1",*{
      \vrule height 1pt depth 1pt width 0pt
      \vrule height 0pt depth 0pt width 2pt
      }="#1",}
  \xymatrix @!0 @R=25mm @C=25mm  { 
   \mathcal{C} & \  \\
   \mathbb{S} \ar[u]|*{}="A"^{F} \ar[r]_{j}  
     \enlargexyentry A
    \ar@2{<-}[r];"A"^-{\,} 
    & \mathbb{L}     \ar[lu]_{\operatorname{Lan}_{j} F}  
   \\}
\end{tabular}
\end{center}
so that there is a bijection $[\operatorname{Lan}_{j} F,G]\cong [F,G\circ j]$ between sets of natural transformations, and this for every functor $G\,:\,\mathbb{L}\rightarrow \mathcal{C}$. This induces a bijection:
\begin{center}
Mod($\mathbb{L},\mathcal{C}$)($\operatorname{Lan}_{j} F,G)\cong\,$Mod($\mathbb{S},\mathcal{C})(F,U_j G)$
\end{center}
which is natural in $G$, so that $\operatorname{Lan}_{j} F$ should be thought of as the \emph{free} $\mathbb{L}$-model generated by the model $F$ in the Lawvere theory $\mathbb{L}$. Moreover, this construction is functorial, so that $\operatorname{Lan}_{j}$ gives the expected left adjoint of $U_j$: a free model of a given theory in a given category is just obtained by the change-of-base operation provided by the left Kan extension !\\
If we take the theory $\mathbb{G}$ of groups as an example, with models in $Set$, we have the following situation:
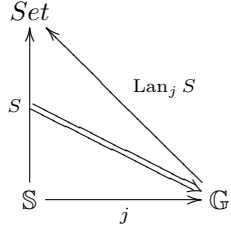
\begin{center}
\begin{tabular}{c}
  \def\enlargexyentry#1{%
    \POS "#1",*{
      \vrule height 1pt depth 1pt width 0pt
      \vrule height 0pt depth 0pt width 2pt
      }="#1",}
  \xymatrix @!0 @R=25mm @C=25mm  { 
   Set & \  \\
   \mathbb{S} \ar[u]|*{}="A"^{S} \ar[r]_{j}  
     \enlargexyentry A
    \ar@2{<-}[r];"A"^-{\,} 
    & \mathbb{G}     \ar[lu]_{\operatorname{Lan}_{j} S}  
   \\}
\end{tabular}
\end{center}
is a left Kan extension, and $S$ is (represents) a set. So, $\operatorname{Lan}_{j} S$ is a model of the theory of groups in $Set$, whose underlying set is $(\operatorname{Lan}_{j} S)(F(1))$, and the coend formula states that:
\begin{center}
\begin{tabular}{rcl}
$(\operatorname{Lan}_{j} S)(F(1))$&$=$&$\int^{n\in \mathbb{N}}{\mathbb{G}(j[n],F(1))\cdot \mathbb{S}[n]}$\\
&$\cong$&$\int^{n\in \mathbb{N}}{\mathbb{G}(F(n),F(1))\cdot S^n}$\\
\end{tabular}
\end{center}
Since $\mathbb{G}(F(n),F(1))$ is by definition the set of group homomorphisms from the generic free group on $n$ elements to the generic free group on one element, its cardinality is just $1$, and we have:
\begin{center}
\begin{tabular}{rcl}
$(\operatorname{Lan}_{j} S)(F(1))$&$\cong$&$\int^{n\in \mathbb{N}}{\mathbb{S}^n}$\\
&$\cong$&$\coprod_{n\in \mathbb{N}}{\mathbb{S}^n}$\\
\end{tabular}
\end{center}
where in the last step we just have an usual colimit, since the coend is taken over a functor whose contravariant behaviour appears to be dummy, so that it only is an usual covariant functor, a case on which coends meet the usual colimits. So, the Kan extension builds out of a set $S$ the free group on it (since it is a model of the theory of groups, and since its underlying set is the one of a free group).

\end{remark}

\subsection{Yoneda structures}
\label{yoneda}
Yoneda structures were introduced in 1978 by Street and Walters in \cite{yonedastruct}, and reexplored recently by Weber in \cite{yonedatopos} - where most proofs of the results presented in this section can be found. What follows is strongly inspired by Melli\`es' slides \cite{yonedamellies}. The point of these structures is to axiomatise the idea of a nerve in a general way. The notion of left Kan extension will play a central role here, but the notion of left Kan lifting will be required as well.

\begin{definition}[Left Kan lifting]

One says that the following $2$-cell:

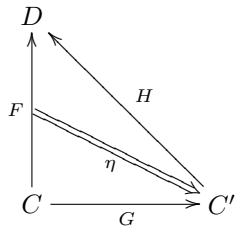
\begin{center}
\begin{tabular}{c}
  \def\enlargexyentry#1{%
    \POS "#1",*{
      \vrule height 1pt depth 1pt width 0pt
      \vrule height 0pt depth 0pt width 2pt
      }="#1",}
  \xymatrix @!0 @R=25mm @C=25mm  { 
    D & \  \\
    C \ar[u]|*{}="A"^F \ar[r]_G  
     \enlargexyentry A
    \ar@2{<-}[r];"A"^-{\eta} 
    & C'     \ar[lu]_{H}  
   \\}
\end{tabular}
\end{center}
exhibits the arrow $G$ as an absolute left lifting of $F$ through $H$ when it has the following universal property: every 2-cell from $F$ to a functor composed with $H$ factors through the previously mentionned $2$-cell:\\
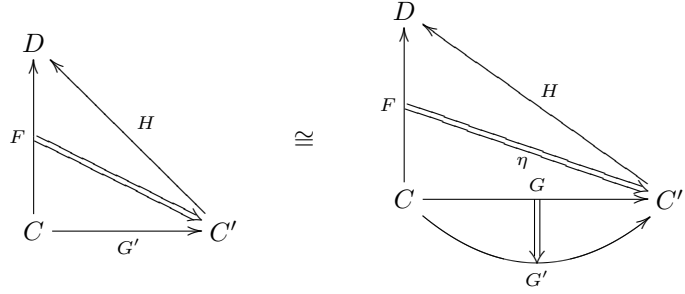
\begin{center}
\begin{tabular}{ccc}
\begin{tabular}{c}
 \def\enlargexyentry#1{%
    \POS "#1",*{
      \vrule height 1pt depth 1pt width 0pt
      \vrule height 0pt depth 0pt width 2pt
      }="#1",}
  \xymatrix @!0 @R=25mm @C=25mm  { 
    D & \  \\
    C \ar[u]|*{}="A"^F \ar[r]_{G'}  
     \enlargexyentry A
    \ar@2{<-}[r];"A"^-{\ } 
    & C'     \ar[lu]_{H}  
   \\}
   \end{tabular}
&  $\cong$ &
\begin{tabular}{c}
  \def\enlargexyentry#1{%
    \POS "#1",*{
      \vrule height 1pt depth 1pt width 0pt
      \vrule height 0pt depth 0pt width 2pt
      }="#1",}
  \xymatrix @!0 @R=25mm @C=35mm  { 
    D & \  \\
    C \ar[u]|*{}="A"^F \ar[r]^G="B" \ar@/_2pc/[r]_{G'}="C"
     \enlargexyentry A 
    \ar@2{<-}[r];"A"^-{\eta} 
    & C'   \ar[lu]_{H}     
    \ar@{=>}"B";"C"  
   \\}
\end{tabular}   
\\
\end{tabular}
\end{center}

this providing a bijection between $2$-cells $G\Rightarrow G'$ and $2$-cells $F \Rightarrow H \circ G'$. The lifting is \emph{absolute} when every arrow out of $D$ preserves it.

\end{definition}

We can then define the notion of Yoneda structure.

\begin{definition}[Yoneda structure]
A \emph{good Yoneda structure} on a $2$-category $\mathcal{K}$ consists of:
\begin{enumerate}
\item a right ideal of admissible arrows (this meaning that if an arrow $g$ is admissible, for every precomposable arrow $f$, $g \circ f$ is admissible)
\item a collection of admissible objects, which are the ones whose identity is admissible
\item an admissible object $\widehat{A}$ and an admissible map $y_A\,:\,A \rightarrow \widehat{A}$ for each admissible object $A$
\item an arrow $B(f,1)$ and a $2$-cell:

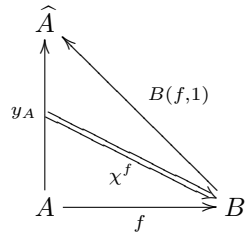
\begin{center}
\begin{tabular}{c}
  \def\enlargexyentry#1{%
    \POS "#1",*{
      \vrule height 1pt depth 1pt width 0pt
      \vrule height 0pt depth 0pt width 2pt
      }="#1",}
  \xymatrix @!0 @R=25mm @C=25mm  { 
    \widehat{A} & \  \\
    A \ar[u]|*{}="A"^{y_A} \ar[r]_{f}  
     \enlargexyentry A
    \ar@2{<-}[r];"A"^-{\chi^f} 
    & B   \ar[lu]_{B(f,1)}  
   \\}
\end{tabular}
\end{center}

for every admissible object $A$ and every admissible arrow $f\,:\,A \rightarrow B$
\end{enumerate}
which satisfy the following axioms:
\begin{enumerate}
\item the provided $2$-cells exhibit $B(f,1)$ as a left extension of $y_A$ along $f$
\item the provided $2$-cells exhibit $f$ as an absolute left lifting of $y_A$ through $B(f,1)$
\item for every admissible object $A$, the following $2$-cell:
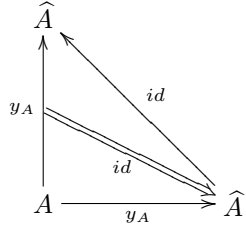
\begin{center}
\begin{tabular}{c}
  \def\enlargexyentry#1{%
    \POS "#1",*{
      \vrule height 1pt depth 1pt width 0pt
      \vrule height 0pt depth 0pt width 2pt
      }="#1",}
  \xymatrix @!0 @R=25mm @C=25mm  { 
    \widehat{A} & \  \\
    A \ar[u]|*{}="A"^{y_A} \ar[r]_{y_A}  
     \enlargexyentry A
    \ar@2{<-}[r];"A"^-{id} 
    & \widehat{A}   \ar[lu]_{id}  
   \\}
\end{tabular}
\end{center}
exhibits the identity on $A$ as a left extension of $y_A$ along $y_A$
\item for every couple $(A,B)$ of admissible objects and every couple\\$(f\,:\,A\rightarrow B,g\,:\,B \rightarrow C)$ of accessible maps, the following $2$-cell (where $B(1,f)$ obviously denotes $y_B \circ f$):

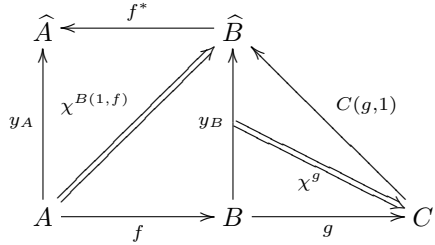
\begin{center}
\begin{tabular}{c}
  \def\enlargexyentry#1{%
    \POS "#1",*{
      \vrule height 1pt depth 1pt width 0pt
      \vrule height 0pt depth 0pt width 2pt
      }="#1",}
  \xymatrix @!0 @R=25mm @C=25mm  { 
    \widehat{A} & \widehat{B} \ar[l]_{f^*} & \  \\
    A \ar[r]_f \ar[u]^{y_A} \ar@{=>}[ru]^{\chi^{B(1,f)}} 
    & B \ar[u]|*{}="A"^{y_B} \ar[r]_{g}  
     \enlargexyentry A
    \ar@2{<-}[r];"A"^-{\chi^g} 
    & C     \ar[lu]_{C(g,1)}  
   \\}
\end{tabular}
\end{center}

exhibits the arrow $f^* \circ C(g,1)$ - where $f^*$ shall be introduced quickly - as a left Kan extension of $y_A$ along $g\circ f$.

\end{enumerate}
\end{definition}

\begin{remark}
Here we follow the axioms of the original paper by Street and Walters \cite{yonedastruct}. In \cite{yonedatopos}, Weber adopts a different axiomatization which exhibits some of our axioms as a consequence of his.
\end{remark}

\begin{example}[A good Yoneda structure on $CAT$] \label{yoncat} The $2$-category $CAT$ of all categories and functors between them has the following obvious Yoneda structure:
\begin{itemize}
\item admissible functors $f\,:\,A \rightarrow B$ are the functors such that for every $a\in A$ and $b\in B$ the homset $B(fa,b)$ is a set
\item admissible objects are the locally small categories
\item for every category $A$, $\widehat{A}$ is defined as the usual presheaf category on $A$, and $y_A$ is the usual Yoneda embedding
\item for every functor $f\,:\,A\rightarrow B$ such that $A$ and $f$ are admissible, $B(f,1)$ is defined as in Section \ref{density}: $B(f,1)(b)(a)\,=\,B(fa,b)$, and $\chi^f$ evaluates on $a$ to $(\chi^f)a\,:\,A(\cdot,a) \mapsto B(f \cdot,fa)$.
\end{itemize}
\end{example}
Yoneda structures can be understood as a formalization and a generalization of this example, relaxing the requirement of hom-sets - somewhere, it is similar to the generalization topoi provide, in the sense that it abstracts some of the properties of sets outside their world. There also are non-trivial examples\footnote{See \cite[Section  7]{yonedastruct} for Yoneda structures on Hom-enriched categories, internal categories and finetely complete categories.}:

\begin{example}[A good Yoneda structure for compact $2$-categories]
A \emph{compact bicategory}\footnote{An interesting - yet degenerated - case of compact bicategories is the one of groupoids, which can be seen as such categories with only identity $2$-cells.} is a bicategory\footnote{A definition may be found at \cite[Section 1.5]{leinster}.} $\mathcal{C}$ in which every arrow has a right adjoint. A Yoneda structure is obtained on $\mathcal{C}$ by taking all arrows to be admissible (and thus also all objects), by defining $\widehat{A}\,=\,A$ and $y_A\,=\,1_A$ for every object $A$ of $\mathcal{C}$, and by taking for every arrow $f\,:\,A\rightarrow B$ its right adjoint as the arrow $B(f,1)$.
\end{example}

The inverse arrow $f^*$\footnote{The operation $(\cdot)^*$ is functorial, see \cite[Section 2]{yonedastruct} for a proof.}, which has been mentionned when we introduced the Kan extension, can be defined in an axiomatized way in Yoneda structures: for every admissible arrow $f\,:\,A \rightarrow B$ with both $A$ and $B$ admissible, the arrow $f^*\,:\,\widehat{B} \rightarrow \widehat{A}$ is defined as follows:
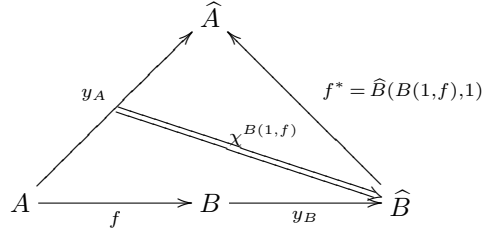
\begin{center}
\begin{tabular}{c}
  \def\enlargexyentry#1{%
    \POS "#1",*{
      \vrule height 1pt depth 1pt width 0pt
      \vrule height 0pt depth 0pt width 2pt
      }="#1",}
  \xymatrix @!0 @R=25mm @C=25mm  { 
    \ & \widehat{A}  &\   \\
    A \ar[ru]|*{}="A"^{y_A} \ar[r]_{f}  
     \enlargexyentry A
    \ar@2{<-}[rr];"A"_-{\chi^{B(1,f)}} 
    & B \ar[r]_{y_B} & \widehat{B} \ar[lu]_{\ f^*\,=\,\widehat{B}(B(1,f),1)}  
   \\}
\end{tabular}

\end{center}and this $2$-cell factors as follows:

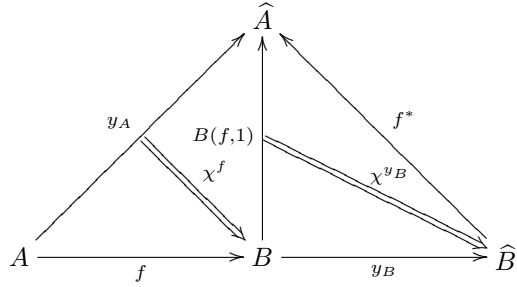
\begin{center}
\begin{tabular}{c}
  \def\enlargexyentry#1{%
    \POS "#1",*{
      \vrule height 1pt depth 1pt width 0pt
      \vrule height 0pt depth 0pt width 2pt
      }="#1",}
  \xymatrix @!0 @R=32mm @C=32mm  { 
    \ & \widehat{A}  &\   \\
    A \ar[ru]|*{}="A"^{y_A} \ar[r]_{f}  
     \enlargexyentry A
    \ar@2{<-}[r];"A"_-{\chi^{f}} 
    & B \ar[r]_{y_B} \ar[u]|*{}="B"^{B(f,1)}
         \enlargexyentry B
    \ar@2{<-}[r];"B"_-{\chi^{y_B}} 
     & \widehat{B} \ar[lu]_{f^*}  
   \\}
\end{tabular}
\end{center}

where both $2$-cells are Kan extensions.\\\ \\

This arrow, in the case of Example \ref{yoncat}, meets the one we introduced before with the same notation. To this generalization of the operation of precomposition corresponds a generalization of its left and right adjoints - remember used to correspond to left and right Kan extensions. Here we use the logical notations\footnote{They are related to the notion of quantification in a topos.} for these adjoints: if $f\,:\,A \rightarrow B$ is an arrow, the adjoints to $f^*\,:\,\widehat{B} \rightarrow \widehat{A}$ are denoted $\exists_f$ (left adjoint to $f^*$) and $\forall_f$ (right one). These arrows are defined as follows:

\begin{itemize}
\item Given an admissible arrow $f\,:\,A\rightarrow B$ with $A$, $\widehat{A}$ and $B$ admissible, $f^*$ has a right adjoint, defined as follows:
\begin{displaymath}
\forall_f\,=\,\widehat{A}(B(f,1),1))
\end{displaymath}
or equivalently by the following commutative square:
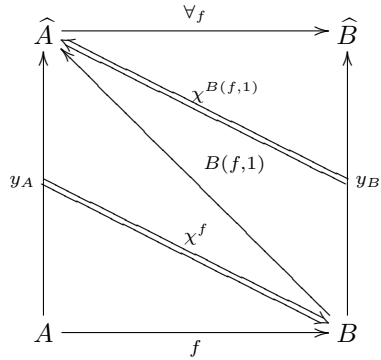
\begin{center}
\begin{tabular}{c}
\xymatrix @!0 @R=40mm @C=40mm {
\widehat{A} \ar[r]^{\forall_f} & \widehat{B} \\
A \ar[r]_f \ar[u]|*{}="A"^{y_A}  \ar@2{<-}[r];"A"_-{\chi^f} & B \ar[u]|*{}="B"_{y_B}  \ar@2{<-}[ul];"B"^-{\chi^{B(f,1)}} \ar[ul]_{B(f,1)\ \,} \\
}
\end{tabular}
\end{center}

\item Given an admissible arrow $f\,:\,A\rightarrow B$ with $A$, $\widehat{A}$ and $B$ admissible, $f^*$ has a left adjoint $\exists_f$, defined as the left Kan extension of $y_B \circ f$ along $y_A$:

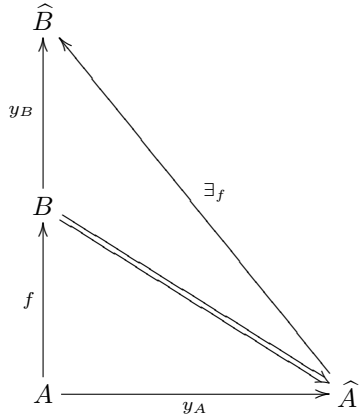
\begin{center}
\begin{tabular}{c}
\xymatrix @!0 @R=25mm @C=40mm  {
\widehat{B} & \ \\
B \ar[u]^{y_B} & \ \\
A \ar[u]^f \ar[r]_{y_A} & \widehat{A} \ar[uul]_{\exists_f} \ar@2{<-}[ul]\\
}
\end{tabular}
\end{center}
\end{itemize}

Again, in $CAT$ together with its canonical Yoneda structure, these definitions meet the usual ones.\\\ \\
The interest of Yoneda structures is to provide a general setting for nerves-like constructions, since it gives a general framework for the change-of-base operation on presheaves - which, remember, were used as models of Lawvere theories and, not surprisingly, will be reused in this purpose when it comes to generalize these theories for different notions of arities.

\section{Monads with arities}

\subsection{Monads with arities}

We want to generalize the fact that the values of finitary monads on any set can be computed from their values on finite ordinals using a left Kan extension.

\subsubsection{A definition via Kan extensions}
\begin{definition}[Monad with arities]
A monad with arities consists of a monad  $(T,\mu,\eta)$ on a category $\mathcal{A}$ and of a fully faithful and dense functor $i_0\,:\,\Theta_0\rightarrow \mathcal{A}$, with $\Theta_0$ as small category, satisfying the following conditions:
\begin{enumerate}
\item The following $2$-cell is a left Kan extension:
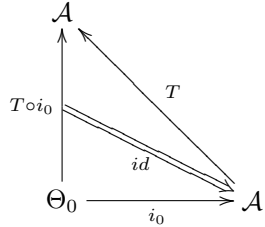
\begin{center}
\begin{tabular}{c}
  \def\enlargexyentry#1{%
    \POS "#1",*{
      \vrule height 1pt depth 1pt width 0pt
      \vrule height 0pt depth 0pt width 2pt
      }="#1",}
  \xymatrix @!0 @R=25mm @C=25mm  { 
   \mathcal{A} & \  \\
   \Theta_0 \ar[u]|*{}="A"^{T\circ i_0} \ar[r]_{i_0}  
     \enlargexyentry A
    \ar@2{<-}[r];"A"^-{id} 
    & \mathcal{A}     \ar[lu]_{T}  
   \\}
\end{tabular}
\end{center}

\item This Kan extension is preserved by the nerve functor \mbox{$\mathcal{A}(i_0,1)\,:\,\mathcal{A}\rightarrow \widehat{\Theta}_0$}
\end{enumerate}

Since $\mathcal{A}(i_0,1)$ is fully faithful, it reflects left Kan extensions; thus, both conditions happen to be equivalent to the following single requirement:
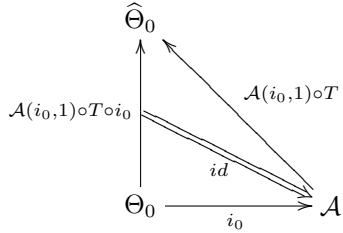
\begin{center}
\begin{tabular}{c}
  \def\enlargexyentry#1{%
    \POS "#1",*{
      \vrule height 1pt depth 1pt width 0pt
      \vrule height 0pt depth 0pt width 2pt
      }="#1",}
  \xymatrix @!0 @R=25mm @C=25mm  { 
   \widehat{\Theta}_0 & \  \\
   \Theta_0 \ar[u]|*{}="A"^{\mathcal{A}(i_0,1) \circ T\circ i_0} \ar[r]_{i_0}  
     \enlargexyentry A
    \ar@2{<-}[r];"A"^-{id} 
    & \mathcal{A}     \ar[lu]_{\mathcal{A}(i_0,1) \circ T}  
   \\}
\end{tabular}
\end{center} is a left Kan extension.
\end{definition}

The first condition should seem natural, by its proximity with the idea of realization: it just means that the values of $T$ on the subcategory of arities are enough to compute all the values of $T$; moreover, the requirement that the universal natural transformation in this $2$-cell is $id$ enforces the exactness of this computation. The second condition means that the computation of the values of $T$ on $\mathcal{A}$ should also be able to be computed in $Set$ using presheaves.\\
The "unified" condition will be the one we use in what follows.

\begin{remark}
According to Melli\`es \cite{lics2010}, this second arity condition was introduced by Weber - it turns out that it didn't appear before since it happens to be satisfied when $\Theta_0$ is the full subcategory of compact\footnote{Also called \emph{finitely presentable}.} objects of a locally finitely presentable category $\mathcal{A}$ (that is, a category  where every object can be computed as the filtered colimit of a diagram of compact objects), since $\mathcal{A}(i_0,1)$ preserves filtered colimits. 
\end{remark}

\begin{remark}
A finitary monad is a monad whose arity functor is just:
\[
i_0\,:\,FinSet\rightarrow Set
\]
\end{remark}

\subsubsection{Monads with arities and factorizations of arrows}
\label{monadfact}
Another formulation of the conditions on arities we just gave can be deduced from the coend formula for left Kan extensions. This formula gives, for our "unified" condition, the value on $A$ of the functor $\mathcal{A}(i_0,1)\circ T$:
\begin{center}
$\mathcal{A}(i_0 \cdot,TA)\,=\,\int^{p\in \Theta_0}{\mathcal{A}(i_0 p,A)\times \mathcal{A}(i_0 \cdot, Ti_0 p)}$
\end{center}
(recall that when we deal with sets the copower happens to meet the usual notion of cartesian product). Since the evaluation of $\mathcal{A}(i_0,1)\circ T$ on an object $A \in \mathcal{A}$ gives a presheaf, there is a $\cdot$ at every place where the evaluation of this presheaf should lead to a substitution: the evaluation of this presheaf on any arity $n \in \Theta_0$ gives:
\begin{center}
$\mathcal{A}(i_0 n,TA)\,=\,\int^{p\in \Theta_0}{\mathcal{A}(i_0 p,A)\times \mathcal{A}(i_0 n, Ti_0 p)}$
\end{center}

So, every element of $\mathcal{A}(i_0 n,TA)$ corresponds to one only element of\\\mbox{$\int^{p\in \Theta_0}{\mathcal{A}(i_0 p,A)\times \mathcal{A}(i_0 n, Ti_0 p)}$}. Since this coend is obtained as a quotient of \mbox{$\mathcal{A}(i_0 p,A)\times \mathcal{A}(i_0 n, Ti_0 p)$}, every element of $\mathcal{A}(i_0 n,TA)$ - that is every arrow \mbox{$i_0 n \rightarrow TA$} in $\mathcal{A}$ - corresponds to an equivalence class of elements of \\\mbox{$\mathcal{A}(i_0 p,A)\times \mathcal{A}(i_0 n, Ti_0 p)$} for an equivalence relation given by the coend - that is, an equivalence class of couples of arrows $(e,f)$ where $e\,:\,i_0 n\rightarrow Ti_0 p$ and $f\,:\,i_0 p \rightarrow A$. This just means that every arrow $i_0 n \rightarrow TA$ in $\mathcal{A}$ factorizes\footnote{Interestingly, there seems to be a similar result on the dual case of Lawvere theories in \cite{linton2}.} as follows:
\begin{center}
$i_0 n \xrightarrow{e} Ti_0 p  \xrightarrow{Tf} TA$
\end{center}
Now we have to find which equivalence relation between these factorizations is induced by the coend. In our case, the universal wedge for every arrow  $u\,:\,p \rightarrow q$ in the category of arities $\Delta_0$ is the following commutative square:

\begin{center}
\begin{tabular}{c}
\xymatrix{
\mathcal{A}(i_0 q,A)\times \mathcal{A}(i_0 n, Ti_0 p) \ar[r]^{f_*} \ar[d]_{f^*} & \mathcal{A}(i_0 q,A)\times \mathcal{A}(i_0 n, Ti_0 q) \ar[d] \\
\mathcal{A}(i_0 p,A)\times \mathcal{A}(i_0 n, Ti_0 p) \ar[r] & \mathcal{A}(i_0 n,TA)\\
}
\end{tabular}
\end{center}

where there is an abuse of notation for the sake of readability: we denote here by $u_*$ the composition with $f$ on the covariant side of the functor and by $u^*$ its composition on the contravariant side.\\
So, every such arrow $p\rightarrow q$ induces an identification in the coend between factorizations via $Ti_0 p$ and via $Ti_0 q$, precisely, the coend identifies two factorizations $(e_1,f_1)$ and $(e_2,f_2)$ when they are equivalent modulo the transitive, symmetric and reflexive closure of the binary relation: $(e_1,f_1) \rightarrow (e_2,f_2)$ iff there exists $u\,:\,p\rightarrow q$ in $\Theta_0$ making the following diagram commute in $\mathcal{A}$:

\begin{center}
\begin{tabular}{c}
\xymatrix @!0 @R=25mm @C=25mm  {
i_0 n \ar[r]^{e_1} \ar[d]_{e_2} & Ti_0 p \ar[d]^{Tf_1} \ar[ld]_{Ti_0 u}\\
Ti_0 q \ar[r]_{Tf_2} & TA\\
}
\end{tabular}
\end{center}

\subsubsection{Finitary monads, a third formulation} 
\label{finfact}
In the particular case of finitary monads, the functor $i$ including the category of arities (finite ordinals) into the category of sets is $i\,:\,Finset \rightarrow Set$. The associated nerve functor is $FinSet(i,1)$. Our second formulation for finitary monads stated that a monad is finitary iff:

\begin{center}
\begin{tabular}{c}
  \def\enlargexyentry#1{%
    \POS "#1",*{
      \vrule height 1pt depth 1pt width 0pt
      \vrule height 0pt depth 0pt width 2pt
      }="#1",}
  \xymatrix @!0 @R=25mm @C=25mm  { 
    Set & \  \\
    FinSet \ar[u]|*{}="A"^{T\circ i} \ar[r]_{i}  
     \enlargexyentry A
    \ar@2{<-}[r];"A"^-{id} 
    & Set     \ar[lu]_{T}  
   \\}
\end{tabular}
\end{center}
is a left Kan extension. Since $Set$ is a locally finitely presentable category (recall that every set may be computed as a filtered colimit of ordinals, so that $FinSet$ can be used as a category of arities for $Set$), the second condition on monad with arities is automatically satisfied, and $i\,:FinSet \rightarrow Set$ being a fully faithful and dense functor, all what we did previously applies. So, if $T$ is a finitary monad, the following $2$-cell is a left Kan extension:

\begin{center}
\begin{tabular}{c}
  \def\enlargexyentry#1{%
    \POS "#1",*{
      \vrule height 1pt depth 1pt width 0pt
      \vrule height 0pt depth 0pt width 2pt
      }="#1",}
  \xymatrix @!0 @R=25mm @C=25mm  { 
    \widehat{FinSet} & \  \\
    FinSet \ar[u]|*{}="A"^{Set(i,1) \circ T\circ i} \ar[r]_{i}  
     \enlargexyentry A
    \ar@2{<-}[r];"A"^-{id} 
    & Set     \ar[lu]_{Set(i,1) \circ T}  
   \\}
\end{tabular}
\end{center}

and the coend formula gives (again, the copowers are just cartesian products when we deal with sets):

\begin{center}
\begin{tabular}{rcl}
$Set(i \cdot,TA)$&$\,=\,$&$\int^{[n]\in FinSet}{Set(i [n],A)\times Set(i \cdot, Ti [n])}$\\
&$\,=\,$&$\int^{[n]\in FinSet}{Set([n],A)\times Set(\cdot, [n])}$\\
\end{tabular}
\end{center}

This is a presheaf on $FinSet$, which evaluates on $[p]$ as follows:
\begin{center}
$Set([p],TA)\,=\,\int^{[n]\in FinSet}{Set([n],A)\times Set([p], [n])}$\\
\end{center}

so that, as previously, we have the factorization property that every arrow $[p]\rightarrow TA$ should factorize as $[p] \xrightarrow{e} T[q] \xrightarrow{Tf} TA$, two factorizations being identified when there exists $u\,:\,[q] \rightarrow [r]$ making the following square commutes:

\begin{center}
\begin{tabular}{c}
\xymatrix @!0 @R=25mm @C=25mm  {
[p] \ar[r]^{e_1} \ar[d]_{e_2} & T[q] \ar[d]^{Tf_1} \ar[ld]_{Tu}\\
T[r] \ar[r]_{Tf_2} & TA\\
}
\end{tabular}
\end{center}

\begin{remark}
According to Melli\`es (private communication), it is enough in the finitary case to check this factorization property for all arrows with domain  $[1]$.
\end{remark}

\subsection{Categories with arities}

After what we have done so far, it is natural to define a notion of category with arities, that is, of a category built out of subpieces:

\begin{definition}[Category with arities]
A category with arities is a couple $(\mathcal{A},i_0)$, with $\mathcal{A}$ a category and $i_0\,:\,\Theta_0 \rightarrow \mathcal{A}$ a fully faithful and dense functor whose domain $\Theta_0$ is a small category.
\end{definition}

The notion of morphism between categories with arities yet is not as simple as we might expect:

\begin{definition}[Morphism of categories with arities]
A morphism between categories with arities $(F,l)\,:\,(\mathcal{A},i_0) \rightarrow (\mathcal{B},i_1)$ is defined as a pair of functors $(F,l)$ making the following diagram commute:
\begin{center}
\begin{tabular}{c}
\xymatrix @!0 @R=17mm @C=17mm  {
\Theta_1 \ar[r]^{i_1} & \mathcal{B}\\
\Theta_0 \ar[u]^{l} \ar[r]_{i_0} & \mathcal{A} \ar[u]_F \\
}
\end{tabular}
\end{center}

Moreover, the induced diagram on presheaves:

\begin{center}
\begin{tabular}{c}
\xymatrix @!0 @R=17mm @C=17mm  {
\widehat{\Theta_1} \ar[d]_{l^*} & \widehat{\mathcal{B}} \ar[d]^{F^*} \ar[l]_{i_1^*}\\
\widehat{\Theta_0} & \widehat{\mathcal{A}} \ar[l]^{i_0^*}  \ar@{=>}[lu]^{id} \\
}
\end{tabular}
\end{center}

is required to be an exact square in the sense of Guitart; equivalently, the Beck-Chevalley condition has to be satisfied - meaning that the natural transformation occuring in the following $2$-cell:
\begin{center}
\begin{tabular}{c}
\xymatrix @!0 @R=17mm @C=17mm  {
\widehat{\Theta_1} \ar[d]_{l^*} \ar[r]^{\forall_{i_1}} & \widehat{\mathcal{B}} \ar[d]^{F^*} \ar@{=>}[ld]\\
\widehat{\Theta_0}  \ar[r]_{\forall_{i_0}} & \widehat{\mathcal{A}}  \\
}
\end{tabular}
\end{center}

has to be reversible (where $\forall_{i_0}$ denotes the right adjoint to $i_0^*$, that is, the right Kan extension along $i_0^{op}$ - the logical notation being related to the notion of quantification in a topos).
\end{definition}

Such functors happen to compose: in fact, this $2$-cell is defined as the \emph{mate} under adjunction of the previous $2$-cell (the identity one $id\,:\,i_0^* \circ F^* \Rightarrow l^* \circ i_1^*$). Mateship is defined in \cite{review}\footnote{The reader willing to consult this article should be aware that the notion of algebra of a monad is different there, yet related: the algebras are defined as $2$-cells acting over $1$-cells, the usual case being obtained by reducing to the case where both the domain of both cells is the unit category $1$.}, where the functoriality of mateship is detailed (cf \cite[Proposition 2.2]{review}): the two mates under adjunction are in bijection as $2$-cells, and this notion of bijection respects compositions and identities, allowing the composition of morphisms between categories with arities. This naturally leads to a notion of category of categories with arities. 

\begin{remark}
According to Melli\`es (private communication), there is another definition of category of categories with arities, which has a $2$-categorical structure:
\begin{itemize}
\item The $0$-cells are the categories with arities - that is, categories $\mathcal{C}$ equipped with a dense generator $A$, giving a dense functor $A \rightarrow \mathcal{C}$,
\item A $1$-cell from a category $\mathcal{C}$ with dense object $A$ to a category $\mathcal{D}$ with dense object $B$ is a functor $F\,:\,\mathcal{C} \rightarrow \mathcal{D}$ such that:
\begin{itemize}
\item $id\,:\,F \circ i_A \Rightarrow F \circ i_A$ exhibits $F$ as a left Kan extension of $F \circ i_A$ along $i_A$
\item This left Kan extension is preserved by the nerve functor $\mathcal{D}(i_B,1)$
\end{itemize}
\item $2$-cells are the usual natural transformations
\end{itemize}
The interesting point with this definition is that its monads precisely are monads with arities. The link between this definition and the previous is still unclear.
\end{remark}

\begin{remark}
The $2$-cell which is the mate of the identity $id\,:\,i_0^* \circ F^* \Rightarrow l^* \circ i_1^*$ factors in fact as the following one:
\begin{center}
\begin{tabular}{c}
\xymatrix  @!0 @R=25mm @C=32mm  {
\widehat{\Theta_1} \ar[r]|*{}="A"^1 \ar@2{->}[d];"A"^-{\epsilon_1} \ar[d]_{\forall_{i_1}} & \widehat{\Theta_1} \ar@2{<-}[d]^{id} \ar[r]^{l^*} & \widehat{\Theta_0} \ar[d]^{\forall_{i_0}}\\
\widehat{B} \ar[ur]_{i_1^*} \ar[r]_{F^*} & \widehat{A} \ar[r]|*{}="B"_{1} \ar[ur]^{i_0^*} \ar@2{<-}[ur];"B"^-{\eta_0} & \widehat{A}\\
}
\end{tabular}
\end{center}
where $\epsilon_1$ denotes the counit of the adjunction $i_1^* \vdash \forall_{i_1}$ and $\eta_0$ the unit of the adjunction $i_0^* \vdash \forall_{i_0}$. So, the requirement is that this composite is reversible.
\end{remark}

\subsection{An abstract Segal condition}
\subsubsection{The Linton condition}
Every monad $T\,:\,\mathcal{C}\rightarrow \mathcal{C}$ induces a commutative diagram:
\begin{center}
\begin{tabular}{c}
\xymatrix {
\mathcal{C}_T \ar[r]^i & \mathcal{C}^T \\
\mathcal{C} \ar[r]_{id} \ar[u]^F & \mathcal{C} \ar[u]_{Free}\\
}
\end{tabular}
\end{center}
where:
\begin{itemize}
\item $F$ is the identity-on-objects functor sending $\mathcal{C}$ in $\mathcal{C}_T$: recall from Section \ref{kleisli} that the category of free algebras for a monad is equivalent to its Kleisli category, whose objects are the same - and, as expected, every arrow $A\rightarrow B$ in $\mathcal{C}$ is sent by $F$ on its corresponding arrow $A\rightarrow TB$ in $\mathcal{C}_T$.
\item $i$ is the comparison functor, sending an object $A$ of the Kleisli category to the canonical free algebra over it $(TA,\mu_A)$.
\end{itemize}

Linton observe in \cite{linton} that the following square is a pullback under certain conditions - Street and Walters showed that this especially was the case in any $2$-category endowed with a Yoneda structure:
\begin{center}
\begin{tabular}{c}
\xymatrix {
\mathcal{C}^T \ar[r] \ar[d] & \widehat{\mathcal{C}_T} \ar[d]^{F^*} \\
\mathcal{C} \ar[r]_{Y} & \widehat{\mathcal{C}} \\
}
\end{tabular}
\end{center}
so that all the algebras of a monad $T$ on a category $\mathcal{C}$ may be computed from the presheaves over its free algebras: a $T$-algebra is just a presheaf on $\mathcal{C}_T$ whose restriction to $\mathcal{C}$ is representable. In \cite{lics2010}, Melli\`es highlights the fact that the functor $\mathcal{C}^T \rightarrow \widehat{\mathcal{C}_T }$ appearing in the previous square is the nerve-alike functor $\mathcal{C}^T(i,1)$ (where $i$ is the comparison functor defined in the first square of this subsection). This functor transports every algebra $(A,h)$ to a presheaf over $\mathcal{C}_T$ which we will call its \emph{nerve} by analogy with the construction previously seen. Moreover $i$ is dense, thus fully faithful\footnote{Remember that this is a very important point since it is what allows the proper glueing of a nerve: "it embeds enough information" or "the Legos set has enough pieces".}, and we have:

\begin{theorem}[Linton condition]
A presheaf $\phi$ on the Kleisli category $\mathcal{C}_T$ is isomorphic to the nerve of an algebra if and only if the presheaf:
\[
\mathcal{C}^{op} \xrightarrow{F^{op}} \mathcal{C}_T^{op} \xrightarrow{\phi} Set
\]
is representable in the category $\mathcal{C}$, i.e. if and only if $F^* \phi$ is isomorphic to $y_A$ for a certain object $A$ of $\mathcal{C}$.
\end{theorem}

\subsubsection{The Segal condition}
Recall that we have seen the Segal condition as a representability property at Section \ref{segalrep}: here we have the same situation, with $F$ playing the r\^ole of $l$ and just $id$ playing the one of $i_0$: the Linton condition\footnote{The connection between Segal and Linton conditions was first established by Melli\`es in \cite{lics2010}.} is in fact some kind of a Segal condition, where the arity functor is trivial. Nevertheless this connection between the two conditions is crucial, since the Linton condition opens this work on arities to the whole world of monads - and among them the free category monad, which gives the following commutative square:
\begin{center}
\begin{tabular}{c}
\xymatrix {
FreeCat \ar[r]^i & Cat \\
Graph \ar[r]_{id} \ar[u]^F & Graph \ar[u]_{Free}\\
}
\end{tabular}
\end{center}
which is almost the square we had when exhibiting the Segal condition as a representability property except that we have $id$ instead of $i_0$: the notion of monad with arities precisely provides the convenient framework for replacing $id$ by $i_0$, as we shall see now (following Melli\`es \cite{lics2010}, where the demonstrations may be found). Recall first that a functor is \emph{essentially surjective} when every object of the target category is isomorphic to the image by this functor of an object of the source category, and that a presheaf over a category $\Theta$, that is, an element of $\widehat{\Theta}$, is called \emph{representable along the functor $i_0\,:\,\Theta \rightarrow \mathcal{C}$} when it is isomorphic to the restriction along $i_0$ of a presheaf representable in $\mathcal{C}$ - or, equivalently, when it is isomorphic to a presheaf $\mathcal{C}(i_0,A)$ for a given object $A$. 
\begin{proposition}
For every morphism between categories with arities:
\[
(F,l)\,:\,(\mathcal{A},i_0) \rightarrow (\mathcal{B},i_1)
\]
the adjunction $i_1^* \vdash \forall i_1$ induces an adjunction between the full subcategory of $\widehat{\mathcal{B}}$ of presheaves over $\mathcal{B}$ whose restriction along $F$ is representable in $\mathcal{A}$ and the full subcategory of $\widehat{\Theta_1}$ of presheaves over $\Theta_1$ whose restriction along $l$ is representable along $i_0$. Moreover, when $F$ is essentially surjective, these two subcategories are equivalent.
\end{proposition}

\begin{proposition}
Every monad with arities $(T,i_0)$ induces a commutative diagram:
\begin{center}
\begin{tabular}{c}
\xymatrix {
\Theta_T \ar[r]^{i_T} & \mathcal{A}_T \\
\Theta_0 \ar[r]_{i_0} \ar[u]^l & \mathcal{A} \ar[u]_{F}\\
}
\end{tabular}
\end{center}
where $(F,l)\,:\,(\mathcal{A},i_0)\rightarrow (\mathcal{A}_T,i_T)$ is a morphism of categories with arities.
\end{proposition}

This corresponds to what we had in our previous version of the Segal condition: recall we had the following commutative diagram:

\begin{center}
\begin{tabular}{c}
  \xymatrix{
    \Delta \ar[r]^i & Cat \\
    \Delta_0 \ar[r]_{i_0} \ar[u]^{l}& Graph \ar[u]_{Free} \\
  }
\end{tabular}  
\end{center}

which illustrates the second proposition: for every monad with arities, its free algebras may either be computed as the free algebras over arities injected in the target category, or injected first in the category where the free algebras are then computed.\\

These propositions combine to give the following abstract Segal condition:

\begin{theorem}[An abstract Segal condition]
The canonical functor:
\begin{center}
$H\,:\,\Theta_T \xrightarrow{i_T} \mathcal{A}_T \rightarrow \mathcal{A}^T$
\end{center}
induces an equivalence $\mathcal{A}^T \xrightarrow{\mathcal{A}^T(H,1)} \widehat{\Theta_T}$ between $\mathcal{A}^T$ and the full subcategory of $\widehat{\Theta_T}$ of presheaves whose restriction along $l$ is representable along $i_0$.
\end{theorem}

Note that this is exactly what we had in the case of both Linton and Segal condition, which are properly generalized here. What differs from the Linton condition is exactly what we claimed would: the representability along $id$, which induced the use of an ordinary Yoneda embedding, has been changed to the representability along the arity inclusion functor $i_0$ - as in the Segal condition.

\subsection{Lawvere theories with arities}

To the notion of monad with arities correspond an equivalent notion of Lawvere theory with arities, generalizing what we saw at Section 1.1.

\begin{definition}[Lawvere theory with arities] A Lawvere theory $\mathbb{L}$ on a category $\mathcal{C}$ with arities functor $i_0\,:\,\Theta_0 \rightarrow \mathcal{C}$ is defined as an identity-on-objects functor $\mathbb{L}\,:\,\Theta_0 \rightarrow \Theta_{\mathbb{L}}$ such that the endofunctor $\widehat{\Theta_0} \xrightarrow{\exists_{\mathbb{L}}} \widehat{\Theta_{\mathbb{L}}} \xrightarrow{\mathbb{L}^*} \widehat{\Theta_0}$ preserves representability along $i_0$.

\end{definition}

The condition of representability preservation may seem a bit strange; but it is just a requirement making the free construction work as we shall see now.

\begin{definition}[Models of a theory with arities]
A model of the Lawvere theory $\mathbb{L}$ with arity functor $i_0$ is a presheaf $\phi \in \widehat{\Theta_\mathbb{L}}$ whose restriction along $\mathbb{L}$:
\[
\Theta_0^{op} \xrightarrow{\mathbb{L}^{op}} \Theta_{\mathbb{L}}^{op} \xrightarrow{\phi} Set
\]
is representable along $i_0$.
\end{definition}
This defines a category $Mod(\mathbb{L})$, which is a full subcategory of $\widehat{\Theta_{\mathbb{L}}}$. There is then a forgetful functor $U\,:\,Mod(\mathbb{L})\rightarrow \mathcal{C}$, extending the usual one, defined as the unique functor (up to natural isomorphism) making the following diagram commute:
\begin{center}
\begin{tabular}{c}
\xymatrix {
Mod(\mathbb{L}) \ar[d] \ar[rr]^U & & \mathcal{C} \ar[d]^y\\
\widehat{\Theta_{\mathbb{L}}} \ar[r]_{\mathbb{L}^*} & \widehat{\Theta_0} \ar[r]_{\forall{i_0}} & \widehat{C}\\
}
\end{tabular}
\end{center}
The preservation-of-representability property of Lawvere theories with arities then ensures that $U$ has a left adjoint $Free$, making the following diagram commute (again, up to natural isomorphism):
\begin{center}
\begin{tabular}{c}
\xymatrix {
Mod(\mathbb{L}) \ar[d] & & \mathcal{C} \ar[ll]_{Free} \ar[d]^y \\
\widehat{\Theta_{\mathbb{L}}} & \widehat{\Theta_0} \ar[l]^{\exists_{\mathbb{L}^*}}  & \widehat{C} \ar[l]^{i_0^*}\\
}
\end{tabular}
\end{center}

We then have the two following results:

\begin{proposition}
Every monad $T$ with arity functor $i_0$ induces a Lawvere theory $\mathbb{L}_T\,:\,\Theta_0 \rightarrow \Theta_{\mathbb{L}}$.
\end{proposition}

\begin{proposition}
When $T$ is the monad induced by the adjunction $Free \vdash U$, it has arity functor $i_0$, and its induced Lawvere theory $\mathbb{L}_T\,:\,\Theta_0 \rightarrow \Theta_T$ coincides with the theory $\mathbb{L}\,:\,\Theta_0 \rightarrow \Theta_{\mathbb{L}}$ (up to isomorphism between $\Theta_T$ and $\Theta_{\mathbb{L}}$).
\end{proposition}

One may feel that we are getting closer to a generalization of the traditional equivalence of Lawvere theories and finitary monads; this is indeed the case. We have a notion of category of Lawvere theories for a given arity over a given category, whose morphisms are:

\begin{definition}[Morphisms of Lawvere theories with arities]
Given an arity functor $i_0$, a morphism $\mathbb{L}_1 \rightarrow \mathbb{L}_2$ between two Lawvere theories with arities $i_0$ is defined as an identity-on-objects functor $\theta\,:\,\Theta_{\mathbb{L}_1} \rightarrow \Theta_{\mathbb{L}_2}$ making the following diagram commute:
\begin{center}
\begin{tabular}{c}
\xymatrix {
\Theta_{\mathbb{L}_1} \ar[rr]^{\theta} & \ & \Theta_{\mathbb{L}_2}\\
\ & \Theta_0 \ar[ul]^{\mathbb{L}_1} \ar[ur]_{\mathbb{L}_2} & \ \\
}
\end{tabular}
\end{center}
\end{definition}

This provides a notion of category of Lawvere theories with arity functor $i_0\,:\,\Theta_0 \rightarrow \mathcal{C}$, denoted $Law(\mathcal{C},i_0)$. There also is a category of monads with arity functor $i_0$, denoted $Mon(\mathcal{C},i_0)$ - the morphisms between monads corresponding to the usual definition, since the arity functor is the same for every monad in this category. The traditional correspondence then extends in the following way:

\begin{theorem}[Equivalence between Lawvere theories with arities and monads with arities]
The categories $Mon(\mathcal{C},i_0)$ and $Law(\mathcal{C},i_0)$ are equivalent.
\end{theorem}

\subsection{Arities for higher categories}
\label{higherarities}
\subsubsection{Arities for categories}
We start here by showing formally that the free category monad is in fact a monad with arity functor $i_0\,:\,\Delta_0 \rightarrow Graph$. First of all, we need to recall the structure of $Graph$:

\begin{definition}[Graphs]
A \emph{graph} is a presheaf over the following category:
\begin{center}
\begin{tabular}{ccc}
$\mathbb{G}$ & $\ =\ $ &
\xymatrix {
0 \ar@/^/[r]^s \ar@/_/[r]_t & 1\\
}\\
\end{tabular}
\end{center}
\end{definition}

Since a presheaf over $\mathbb{G}$ is a functor $\mathbb{G}^{op} \rightarrow Set$, a morphism of graphs $G_1 \rightarrow G_2$ is a natural transformation $\theta\,:\,G_1 \Rightarrow G_2$; so that it has to make the two following diagrams commute:
\begin{center}
\begin{tabular}{ccc}
\xymatrix {
G_1 1 \ar[r]^{\theta_{1}} \ar[d]_{G_1 (s^{op})} & G_2 1 \ar[d]^{G_2 (s^{op})}\\
G_1 0 \ar[r]_{\theta_{0}} & G_2 0\\
}
& $\ \ \ $ &
\xymatrix {
G_1 1 \ar[r]^{\theta_{1}} \ar[d]_{G_1 (t^{op})} & G_2 1 \ar[d]^{G_2 (t^{op})}\\
G_1 0 \ar[r]_{\theta_{0}} & G_2 0\\
}\\
\end{tabular}
\end{center}

This just means that an edge of the initial graph has to be sent to the target graph "glued with its source and target vertices": their images by $\theta$ have to be the source and target of the image of the edge. Basically, this is a functoriality-like property, over graphs instead of over categories. There is a category of graphs $Graph$, with these objects and morphisms. Recall the interpretation of $\Delta_0$ from Section \ref{segalrep}: $i_0$ takes an object - a finite ordinal $[n]$ - to the path of length $n$ seen as a filiform graph. Since the morphisms of $\Delta_0$ correspond to distance-preserving functions over the ordinals, their interpretation once sent to $Graph$ by $i_0$ is that they just are path inclusions in wider paths (the contraction of an edge would violate the requirement of distance preservation). We are now ready to show that the free category monad has arity functor $i_0$ thanks to the factorization property from Section \ref{monadfact}:

\begin{itemize}
\item $i_0$ is obviously fully faithful, and is dense due to the second item of Definition \ref{defdensity}: it should be clear that every graph $G$ can be reconstructed by glueing as many finite-length paths as can be included in $G$ (this is exactly the meaning of weighting the colimit by the functor $[n] \mapsto Graph(i_0[n],G)$).
\item Now take a morphism $g\,:\,i_0 [n] \rightarrow TG$ in Graph. It maps the object $i$ to an object $gi$, and the edge $e_i\,:\,i-1 \rightarrow i$ ($i\in\{1,\cdots,n\}$) to a path $(f_{i,1},\cdots,f_{i,q_i})$ of length $q_i$ from $g(i-1)$ to $gi$ in $G$. We take for convenience $q_0\,=\,0$. Now $g$ factors as follows:
\[
i_0[n] \xrightarrow{e} Ti_0 [q_1 + \cdots + q_n] \xrightarrow{f} TG
\]
where:
\begin{itemize}
\item $e$ maps the object $i$ of $i_0[n]$ to the object $q_0 + \cdots + q_i$ of $Ti_0[q_1+\cdots+q_n]$, and the edge $e_i\,:\,i-1 \rightarrow i$ to the obvious path of length $q_i$ from $e_{i-1}\,=\,q_0 +  \cdots + q_{i-1}$ to $e_i\,=\,q_0 +  \cdots + q_{i}$, corresponding to the following path over the filiform graph $i_0[q_1+\cdots+q_n]$:
\[
q_0 + \cdots + q_{i-1} \rightarrow q_0 +  \cdots + q_{i-1}+1 \rightarrow \cdots \rightarrow q_0 + \cdots + q_{i}
\]
\item $f$ includes the free category over $i_0[q_1+\cdots+q_n]$ in $TG$:
\begin{itemize}
\item $f$ maps the arrow of $Ti_0[q_1+\cdots+q_n]$ corresponding to the path of length $1$ over the graph $i_0[q_1+\cdots+q_n]$ from $i-1$ to $i$ to $f_{j,k}$, where: $0\leq i \leq q_1+\cdots+q_n$, $j$ is defined as the unique integer such that $i\in  \{q_{j-1} +1, q_{j}\}$, and  $k\,=\,\sum^{j-1}_{l=0}{q_l}$,
\item $f$ maps the arrow of $Ti_0[q_1+\cdots+q_n]$ corresponding to the path $(h_1,\cdots,h_l)$ of length $l$ over the graph $i_0[q_1+\cdots+q_n]$ to the arrow of $TG$ corresponding to the concatenation of the path of paths over $G$ $(f(h_1),\cdots,f(h_l))$,
\item $f$ maps the objects of $Ti_0[q_1+\cdots+q_n]$ in the only obvious way which respects the sources and targets of arrows.
\end{itemize}
\end{itemize}
An example makes things easier to understand. Take as $G$ the following graph:
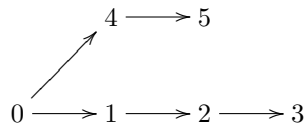
\begin{center}
\begin{tabular}{c}
\xymatrix {
\ & 4 \ar[r]& 5 & \ \\
0 \ar[ur] \ar[r] & 1 \ar[r] & 2 \ar[r] & 3\\
}\\
\end{tabular}
\end{center}
and take the following morphism $i_0[2]\rightarrow TG$, whose action on arrows, depicted by the dotted lines, is to take $0\rightarrow 1$ (in $i_0[2]$) to the path $0\rightarrow 1 \rightarrow 2$ of $TG$ (depicted here as an arrow $0\rightarrow 2$), and $1\rightarrow 2$ (in $i_0[2]$) to $2\rightarrow 3$ (seen as a path of length 1 in $TG$):
\ \\\ \\
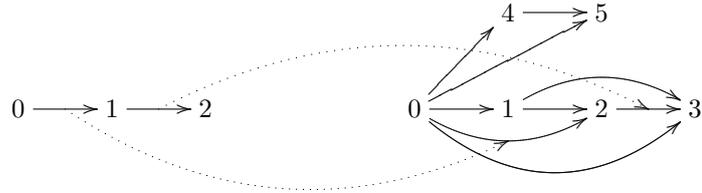
\begin{center}
\begin{tabular}{c}
\xymatrix {
\ & \ & \ & \ & \ & 4 \ar[r]& 5 & \ \\
0 \ar[r]|*{}="A" & 1 \ar[r]|*{}="C"  & 2  &\ \ \ \  & 0 \ar[r] \ar[urr] \ar[ur] \ar@/_1pc/[rr]|*{}="B" \ar@/_2pc/[rrr] & 1 \ar[r] \ar@/^1pc/[rr] & 2 \ar[r]|*{}="D" & 3 \ar @/_2pc/ @{.>} "A";"B" \ar @/^2pc/ @{.>} "C";"D"\\
}\\
\end{tabular}
\end{center}\ \\ \ \\ \ \\
(in $TG$, the identities are not depicted but should be). The first morphism of the factorization is the following morphism $i_0[2] \rightarrow i_0[2+1]$, whose action on arrows, depicted by the dotted lines, is to take $0\rightarrow 1$ to the path $0\rightarrow 1 \rightarrow 2$ (depicted here as an arrow $0\rightarrow 2$) and $1\rightarrow 2$ to $2\rightarrow 3$:
\ \\\ \\
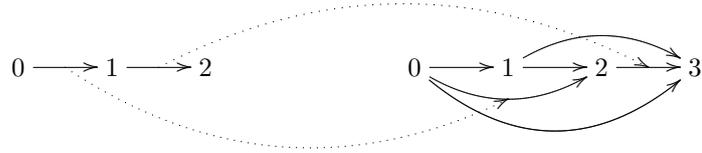
\begin{center}
\begin{tabular}{c}
\xymatrix {
0 \ar[r]|*{}="A" & 1 \ar[r]|*{}="C"  & 2  &\ \ \ \  & 0 \ar[r] \ar@/_1pc/[rr]|*{}="B" \ar@/_2pc/[rrr] & 1 \ar[r] \ar@/^1pc/[rr] & 2 \ar[r]|*{}="D" & 3 \ar @/_2pc/ @{.>} "A";"B" \ar @/^2pc/ @{.>} "C";"D"\\
}\\
\end{tabular}
\end{center}
\ \\ \ \\ \ \\
(where the identities of $Ti_0[3]$ are, again, not depicted) - and then $f$ is just the obvious inclusion of $i_0[2+1]$ in $TG$.

\item This factorization is unique modulo zig-zag: it is clear that there is no such factorization by $Ti_0[p]$ for $p<q_1+ \cdots + q_n$, and any factorization via a greater ordinal $[p]$ amounts to send a wider free category over a filiform graph in $TG$, which leads to a zig-zag commuting diagram:
\begin{center}
\begin{tabular}{c}
\xymatrix @!0 @R=25mm @C=45mm  {
i_0 [n] \ar[r]^{e} \ar[d]_{e'} & Ti_0 [q_1+\cdots+q_n] \ar[d]^{Tf} \ar[ld]_{Ti_0 u}\\
Ti_0 [p] \ar[r]_{Tf'} & TG\\
}
\end{tabular}
\end{center}
since all the arrows of $\Delta_0$ are strict inclusions of paths in wider ones.
\end{itemize}

\subsubsection{Arities for $2$-categories}
First of all, we need to generalize graphs to higher dimensions: this is the purpose of $n$-globular sets.
\begin{definition}[$n$-globular set]
Define $\mathbb{G}_n$ as the category:
\begin{center}
\begin{tabular}{ccc}
$\mathbb{G}_n$ & $\ =\ $ &
\xymatrix {
0 \ar@/^/[r]^{s_1} \ar@/_/[r]_{t_1} & 1 \ar@/^/[r]^{s_2} \ar@/_/[r]_{t_2} & \cdots \ar@/^/[r]^{s_n} \ar@/_/[r]_{t_n} & n \\
}\\
\end{tabular}
\end{center}
such that, for $i\in \{2,\cdots,n\}$, we have $s_i \circ s_{i-1}\,=\,t_i \circ s_{i-1}$ and $t_i \circ t_{i-1}\,=\,s_i \circ t_{i-1}$. A \emph{$n$-globular set} is then a presheaf over $\mathbb{G}_n$\footnote{When $n \rightarrow \infty$, this defines a \emph{globular set} and allows the definition of $\omega$-categories: for details see \cite{leinster}.}.
\end{definition}

To understand these two equalities: when it comes to presheaves, the order of composition is reversed; the idea is that for $i\geq 2$, given a $i$-cell (that is, an element of the image of $i$ by the presheaf), its source and target are $(i-1)$-cells, which should have the same sources and targets (which are $(i-2)$-cells).

\begin{example}[A $2$-globular set]
A very simple example of a $2$-globular set is the presheaf $X\,:\,\mathbb{G}_2^{op}\rightarrow Set$, such that:
\begin{itemize}
\item $X[0]\,=\,\{a,b\}$
\item $X[1]\,=\,\{a\xrightarrow{f} b,a\xrightarrow{g} b\}$
  \newcommand{\UN}[4][r]{%
    \ar@/^1pc/[#1]^{#2}_*=<0.3pt>{}="HAUT"
    \ar@/_1pc/[#1]_{#3}^*=<0.3pt>{}="BAS"
    \ar @{=>} "HAUT";"BAS" ^{#4}
  }
\item $X[2]$ is the singleton \begin{tabular}{c} \xymatrix{
    a \UN[r]{f}{g}{\alpha} & b } \end{tabular}
\end{itemize}
\end{example}

Since globular sets are functors, their morphisms are natural transformations and thus for every integer $1\leq i\leq n$ a morphism $\theta\,:\,\mathbb{G}_1 \rightarrow \mathbb{G}_2$ of $n$-globular sets makes the following diagrams commute:

\begin{center}
\begin{tabular}{ccc}
\xymatrix {
G_1 i \ar[r]^{\theta_{i}} \ar[d]_{G_1 (s_i^{op})} & G_2 i \ar[d]^{G_2 (s_i^{op})}\\
G_1 (i-1) \ar[r]_{\theta_{i-1}} & G_2 (i-1)\\
}
& $\ \ \ $ &
\xymatrix {
G_1 i \ar[r]^{\theta_{i}} \ar[d]_{G_1 (t_i^{op})} & G_2 i \ar[d]^{G_2 (t_i^{op})}\\
G_1 (i-1) \ar[r]_{\theta_{i-1}} & G_2 (i-1)\\
}\\
\end{tabular}
\end{center}

So, generalizing what we had in $Graph$, a morphism of $n$-globular sets has to take the source and target of a $i$-cell to the source and target of its image.\\\ \\

Now we need the notion of \emph{pasting diagram}, described by Leinster in \cite[Section 8.1]{leinster} - we will only use an informal version of them, the reader willing to learn more about them being strongly encouraged to have a look at the original source. This notion is dual to the one of \emph{level tree} by Batanin (see \cite[p. 268]{leinster} for details), and we prefer it just because it is geometrically explicit.
\label{pd}
\begin{definition}[Pasting diagram]\ 
\begin{itemize}
\item A $0$-pasting diagram is an element of $pd(0)\,=\,\{\bullet\}$,
\item A $n$-pasting diagram, whose set is denoted by $pd(n)$, is a finite sequence of $(n-1)$-pasting diagrams.
\end{itemize}
\end{definition}

This has to be understood pictorially:
\begin{itemize}
\item A $1$-pasting diagram is a sequence of dots: if it has $n+1$ dots, we can write it $(\bullet,\cdots,\bullet)$, or even $\bullet \rightarrow \cdots \rightarrow \bullet$: this is exactly the filiform graph shapes we had in $[n] \in \Delta_0$ !
\item A $2$-pasting diagram can thus be understood as a sequence of elements of $\Delta_0$: on an example, the idea is to think of:
\begin{center}
$(\bullet \rightarrow \bullet \rightarrow \bullet,\bullet,\bullet \rightarrow \bullet)$
\end{center}
as:
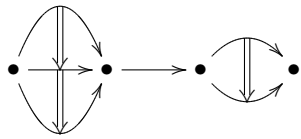
\begin{center}
\begin{tabular}{c} 

  \newcommand{\UN}[4][r]{%
    \ar@/^1pc/[#1]^{#2}_*=<0.3pt>{}="HAUT"
    \ar@/_1pc/[#1]_{#3}^*=<0.3pt>{}="BAS"
    \ar @{=>} "HAUT";"BAS" ^{#4}
  }
  \newcommand{\DEUX}[6][r]{
    \ar@/^2pc/[#1]^{#2}_*=<0.3pt>{}="HAUT"
    \ar@{}    [#1]     ^*=<0.3pt>{}="MILIEUHAUT"
                       _*=<0.3pt>{}="MILIEUBAS"
    \ar[#1]_(0.3){#3}                  
    \ar@/_2pc/[#1]_{#4}^*=<0.3pt>{}="BAS"
    \ar @{=>} "HAUT";"MILIEUHAUT" ^{#5}
    \ar @{=>} "MILIEUBAS";"BAS" ^{#6}}
    
\xymatrix{
    \bullet \DEUX[r]{\ }{\ }{\ }{\ }{\ } & \bullet \ar[r]& \bullet \UN[r]{\ }{\ }{\ } & \bullet } \end{tabular}
\end{center}
\end{itemize} 

Basically, a $2$-pasting diagram realizes as a $2$-globular set (so, a generalization in dimension $2$ of a graph), of finite size both horizontally and vertically, and without branching in the compositions: it is still in a sense filiform, but with one more dimension. 

\begin{remark}[Batanin's level trees]
The link of our notion  of pasting diagram with the one of level tree by Batanin can be intuited on an example: the $2$-pasting diagram we just took as an example can be equivalently portrayed by the following tree:
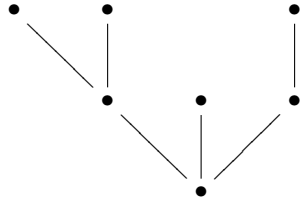
\begin{center}
\begin{tabular}{c}
\xymatrix {
\bullet & \bullet &\ & \bullet \\
\ & \bullet \ar@{-}[ul] \ar@{-}[u] & \bullet & \bullet  \ar@{-}[u]\\
\ & \ & \bullet \ar@{-}[ur] \ar@{-}[u] \ar@{-}[ul] & \ \\
}
\end{tabular}
\end{center}
The idea is that the $2$-pasting diagram it represents is three $1$-cells long: so the tree starts growing by branching in three directions. Every of these directions represent a region of the pasting diagram: the left node which just grew represents the left part of the diagram, which has two $2$-cells: so two more branches start from this node. The middle of the diagram has no $2$-cell, so the middle node of the tree stops growing, and the right node, corresponding to the right part of the diagram, has one $2$-cell, so one branch grows from the node on the right. There are no $3$-cells in our diagram, so the tree stops growing.
\end{remark}

Now, we should introduce the notion of free $2$-category over a $2$-pasting diagram, and then show that the corresponding monad has arities functor $i_1\,:\,pd(2) \rightarrow \widehat{\mathbb{G}_2}$ - where this functor just takes a $2$-pasting diagram to its interpretation as a $2$-globular set ("a $2$-graph"), realizing the graphical intuition on pasting diagrams we have given just before.

\begin{definition}[Free $2$-category over a $2$-globular set]
Given a $2$-globular set $X$, its \emph{free $2$-category} is the $2$-category with objects the elements of $X(0)$, arrows the paths of edges of $X(1)$ (a nullary path starting on a given element of $X(0)$ being thought of as the $1$-dimensional identity over the corresponding object in the free  $2$-category), and $2$-cells the $2$-dimensional paths over $X$, that is paths starting either:
\begin{itemize}
\item from an element $f \in X(1)$ and stopping there: it is just what should be thought of as the identity $2$-cell $f\Rightarrow f$,
\item from a path of elements of $X(1)$ ending on an element of $X(2)$ whose source's source is the target of the path,
\item or from an element of $X(2)$
\end{itemize}
and obtained from a finite number of such moves over the $2$-globular set:
\begin{itemize}
\item $2$-vertical move: from a $2$-cell whose target is $f \in X(1)$ it is licit to move to a $2$-cell whose source is $f$,
\item $2$-horizontal move: from a $2$-cell whose target's target is $a \in X(0)$ it is licit to move to a $2$-cell whose source's source is $a$,
\item $1$-move: from a $2$-cell whose target's target is $a \in X(0)$ it is licit to move along a finite path of elements of $X(1)$ and then to a $2$-cell whose source's source is the target of the last arrow of the $1$-dimensional path.  
\end{itemize}
where in the $1$-move there may be no target $2$-cell iff the path stops there.
\end{definition}

The idea is that these paths just simulate the compositions in a $2$-category, just as the 1-dimensional paths simulate composition in the free category monad: vertical moves are vertical compositions of $2$-cells, horizontal moves are horizontal compositions of $2$-cells, and the $1$-move gives the notion of composition with $1$-cells (for convenience here it is more a notion of \emph{iterated} composition). From this point of view it should be quite evident to understand which the compositions in the free $2$-category should be.\\\ \\
The interesting point with this definition is that it makes clear that all the $2$-cells of the free $2$-category are just $2$-pasting diagrams. To a $2$-cell of the free $2$-category corresponds (a labelling of) the pasting diagram obtained by the following algorithm:
\begin{itemize}
\item Initialization: 
\begin{itemize}
\item if the path is only a $1$-cell, the $2$-pasting diagram is just a $1$-pasting diagram thought of as a $2$-pasting diagram - it is just $\bullet \rightarrow \bullet$.
\item if the path starts by a path of $1$-cells and then meets a $2$-cell, we start building the pasting diagram as expected:
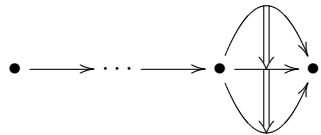
\begin{center}
\begin{tabular}{c} 

  \newcommand{\UN}[4][r]{%
    \ar@/^1pc/[#1]^{#2}_*=<0.3pt>{}="HAUT"
    \ar@/_1pc/[#1]_{#3}^*=<0.3pt>{}="BAS"
    \ar @{=>} "HAUT";"BAS" ^{#4}
  }
  \newcommand{\DEUX}[6][r]{
    \ar@/^2pc/[#1]^{#2}_*=<0.3pt>{}="HAUT"
    \ar@{}    [#1]     ^*=<0.3pt>{}="MILIEUHAUT"
                       _*=<0.3pt>{}="MILIEUBAS"
    \ar[#1]_(0.3){#3}                  
    \ar@/_2pc/[#1]_{#4}^*=<0.3pt>{}="BAS"
    \ar @{=>} "HAUT";"MILIEUHAUT" ^{#5}
    \ar @{=>} "MILIEUBAS";"BAS" ^{#6}}
    
\xymatrix{
    \bullet \ar[r] & \cdots \ar[r] & \bullet \DEUX[r]{\ }{\ }{\ }{\ }{\ } & \bullet } \end{tabular}
\end{center}
\item if the path starts by a $2$-cell, so does the $2$-pasting diagram.
\end{itemize}
\item Iteration:
\begin{itemize}
\item if the next move is a $2$-vertical one: we append to the last added $2$-cell in the $2$-pasting diagram a new $2$-cell just below it,
\item if the next move is a $2$-horizontal one: we append to the last added $2$-cell in the $2$-pasting diagram a new $2$-cell just at its right,
\item if the next move is a $1$-move: we append to the last added $2$-cell a finite path of $1$-cells whose length corresponds to the length of the $1$-path taken in the path describing the $2$-cell, and then a $2$-cell at its right if the path doesn't stop here - in this case, the iteration loop starts again.
\end{itemize}
\end{itemize}

We need now to describe the notion of \emph{composition} of pasting diagrams. It is very easy to understand it pictorially - details may be found in Leinster's book \cite[Section 8.1]{leinster}. Informally, the compositions of pasting diagrams may be described by pasting diagrams: for instance, take the following pasting diagram:

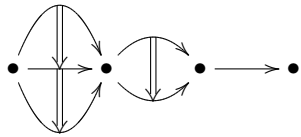
\begin{center}
\begin{tabular}{c} 

  \newcommand{\UN}[4][r]{%
    \ar@/^1pc/[#1]^{#2}_*=<0.3pt>{}="HAUT"
    \ar@/_1pc/[#1]_{#3}^*=<0.3pt>{}="BAS"
    \ar @{=>} "HAUT";"BAS" ^{#4}
  }
  \newcommand{\DEUX}[6][r]{
    \ar@/^2pc/[#1]^{#2}_*=<0.3pt>{}="HAUT"
    \ar@{}    [#1]     ^*=<0.3pt>{}="MILIEUHAUT"
                       _*=<0.3pt>{}="MILIEUBAS"
    \ar[#1]_(0.3){#3}                  
    \ar@/_2pc/[#1]_{#4}^*=<0.3pt>{}="BAS"
    \ar @{=>} "HAUT";"MILIEUHAUT" ^{#5}
    \ar @{=>} "MILIEUBAS";"BAS" ^{#6}}
    
\xymatrix{
    \bullet \DEUX[r]{\ }{\ }{\ }{\ }{\ } & \bullet \UN{\ }{\ }{\ } & \bullet \ar[r] & \bullet } \end{tabular}
\end{center}

and use it as a description of the composition of the following pasting diagrams:

\begin{itemize}
\item The top-left $2$-cell of the pasting diagram should be replaced by:
\begin{center}
\begin{tabular}{c} 
    \xymatrix{
    \bullet \ar[r] & \bullet \ar[r] & \bullet} \end{tabular}
\end{center}

\item The bottom-left $2$-cell of the pasting diagram should be replaced by:

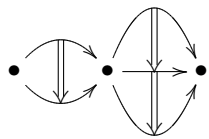
\begin{center}
\begin{tabular}{c} 

  \newcommand{\UN}[4][r]{%
    \ar@/^1pc/[#1]^{#2}_*=<0.3pt>{}="HAUT"
    \ar@/_1pc/[#1]_{#3}^*=<0.3pt>{}="BAS"
    \ar @{=>} "HAUT";"BAS" ^{#4}
  }
  \newcommand{\DEUX}[6][r]{
    \ar@/^2pc/[#1]^{#2}_*=<0.3pt>{}="HAUT"
    \ar@{}    [#1]     ^*=<0.3pt>{}="MILIEUHAUT"
                       _*=<0.3pt>{}="MILIEUBAS"
    \ar[#1]_(0.3){#3}                  
    \ar@/_2pc/[#1]_{#4}^*=<0.3pt>{}="BAS"
    \ar @{=>} "HAUT";"MILIEUHAUT" ^{#5}
    \ar @{=>} "MILIEUBAS";"BAS" ^{#6}}
    
\xymatrix{
    \bullet \UN[r]{\ }{\ }{\ } & \bullet \DEUX{\ }{\ }{\ }{\ }{\ } & \bullet} \end{tabular}
\end{center}

\item The middle $2$-cell of the pasting diagram should be replaced by:

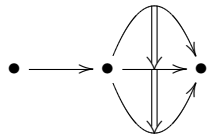
\begin{center}
\begin{tabular}{c} 

  \newcommand{\UN}[4][r]{%
    \ar@/^1pc/[#1]^{#2}_*=<0.3pt>{}="HAUT"
    \ar@/_1pc/[#1]_{#3}^*=<0.3pt>{}="BAS"
    \ar @{=>} "HAUT";"BAS" ^{#4}
  }
  \newcommand{\DEUX}[6][r]{
    \ar@/^2pc/[#1]^{#2}_*=<0.3pt>{}="HAUT"
    \ar@{}    [#1]     ^*=<0.3pt>{}="MILIEUHAUT"
                       _*=<0.3pt>{}="MILIEUBAS"
    \ar[#1]_(0.3){#3}                  
    \ar@/_2pc/[#1]_{#4}^*=<0.3pt>{}="BAS"
    \ar @{=>} "HAUT";"MILIEUHAUT" ^{#5}
    \ar @{=>} "MILIEUBAS";"BAS" ^{#6}}
    
\xymatrix{
    \bullet \ar[r] & \bullet \DEUX[r]{\ }{\ }{\ }{\ }{\ } & \bullet } \end{tabular}
\end{center}

\item The $1$-cell on the right should be replaced by:

\begin{center}
\begin{tabular}{c} 
\xymatrix{
    \bullet \ar[r] & \bullet \ar[r] & \bullet } \end{tabular}
\end{center}

\end{itemize}

this composition giving (as one may imagine, we do vertical composition of the $2$-cells which are stacked, and then horizontal composition of what is left):

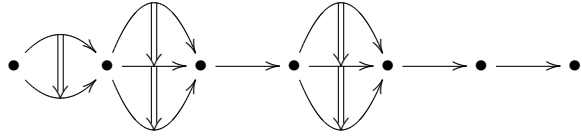
\begin{center}
\begin{tabular}{c} 

  \newcommand{\UN}[4][r]{%
    \ar@/^1pc/[#1]^{#2}_*=<0.3pt>{}="HAUT"
    \ar@/_1pc/[#1]_{#3}^*=<0.3pt>{}="BAS"
    \ar @{=>} "HAUT";"BAS" ^{#4}
  }
  \newcommand{\DEUX}[6][r]{
    \ar@/^2pc/[#1]^{#2}_*=<0.3pt>{}="HAUT"
    \ar@{}    [#1]     ^*=<0.3pt>{}="MILIEUHAUT"
                       _*=<0.3pt>{}="MILIEUBAS"
    \ar[#1]_(0.3){#3}                  
    \ar@/_2pc/[#1]_{#4}^*=<0.3pt>{}="BAS"
    \ar @{=>} "HAUT";"MILIEUHAUT" ^{#5}
    \ar @{=>} "MILIEUBAS";"BAS" ^{#6}}
    
\xymatrix{
    \bullet \UN{\ }{\ }{\ }& \bullet \DEUX{\ }{\ }{\ }{\ }{\ }& \bullet \ar[r] & \bullet \DEUX[r]{\ }{\ }{\ }{\ }{\ } & \bullet \ar[r] &\bullet \ar[r] & \bullet } \end{tabular}
\end{center}

It should be quite clear that there is a monad $T$ over $2$-globular set sending them to the free $2$-category over them. Its algebras are the $2$-categories. We now show that this is a monad with arity functor $i_1\,:\,pd(2) \rightarrow \widehat{\mathbb{G}_2}$ - where this functor just takes a $2$-pasting diagram to its interpretation as a $2$-globular set ("a $2$-graph"), realizing the graphical intuition on pasting diagrams we have had so far. As in the $1$-dimensional case, the multiplication of the monad is the composition of elements of the free $2$-category, that is, composition $2$-dimensional paths whose shape is the one of $2$-pasting diagrams - so it works just as the composition of $2$-pasting diagrams, and unit is just the inclusion of a $2$-globular set in its free $2$-category.

\begin{itemize}
\item It is clear that $i_1$ is fully faithful, and it is dense since a $2$-globular set may be obtained by glueing of as many copies of $2$-pasting diagrams as can be embedded into it.
\item Now take a $2$-pasting diagram $P$, a $2$-globular set $X$ and a morphism $g\,:\,P \rightarrow TX$. We define as $F$ the $2$-pasting diagram which is obtained as the result of the following process:
\begin{itemize}
\item For each $2$-cell of $P$, take its image by $g$.
\item For each $1$-cell of $P$ with no vertically attached $2$-cell - for example, the middle one in this diagram:
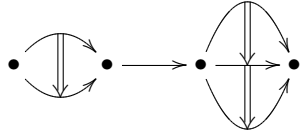
\begin{center}
\begin{tabular}{c} 

  \newcommand{\UN}[4][r]{%
    \ar@/^1pc/[#1]^{#2}_*=<0.3pt>{}="HAUT"
    \ar@/_1pc/[#1]_{#3}^*=<0.3pt>{}="BAS"
    \ar @{=>} "HAUT";"BAS" ^{#4}
  }
  \newcommand{\DEUX}[6][r]{
    \ar@/^2pc/[#1]^{#2}_*=<0.3pt>{}="HAUT"
    \ar@{}    [#1]     ^*=<0.3pt>{}="MILIEUHAUT"
                       _*=<0.3pt>{}="MILIEUBAS"
    \ar[#1]_(0.3){#3}                  
    \ar@/_2pc/[#1]_{#4}^*=<0.3pt>{}="BAS"
    \ar @{=>} "HAUT";"MILIEUHAUT" ^{#5}
    \ar @{=>} "MILIEUBAS";"BAS" ^{#6}}
    
\xymatrix{
    \bullet \UN{\ }{\ }{\ }& \bullet \ar[r] & \bullet \DEUX{\ }{\ }{\ }{\ }{\ }& \bullet } \end{tabular}
\end{center}
take its image by $g$.
\item "Glue the images together": that is, compose all the images - which are elements of the free $2$-category and can thus be thought of as $2$-pasting diagrams - according to the composition information depicted by $P$, inserting the image by $g$ of a given cell at the position of this cell in $P$. This works by structure of the globular set morphims: the images of two parallel $2$-cells by $g$ have to be parallel since their image's sources and targets are the images of their sources and targets.
\end{itemize}
Now there is a factorization of $g$ as follows:
\[
P  \xrightarrow{e} TF \xrightarrow{f} TX
\]
where $e$ maps a cell of $P$ to the cell on $TF$ corresponding to its image by $g$ (that is, the pasting diagram which was inserted in $P$ during the process of building of $F$), and where $f$ sends $TF$ in $TX$ by sending the unary paths (corresponding to $F$ in $TF$) to the cell they correspond to in $TX$ via $g$, and whose image on other paths is given by respecting the free structure of $TF$.
\item A little thought shows that the factorization is unique modulo zig-zag.
\end{itemize}

\paragraph{A Segal condition for $2$-categories} We introduce the following notations:
\begin{itemize}
\item $T$ is, as before, the free $2$-category monad over a $2$-globular set
\item $fpd(2)$ is the category of free $2$-categories over $2$-pasting diagrams
\item $l$ is the map sending $pd(2)$ to $fpd(2)$
\item $2$-$Glob$ is the category of $2$-globular sets
\item $2$-$Cat$ is the category of $2$-categories
\item $i$ is the inclusion of $fpd(2)$ in $2$-$Cat$
\item $i_1$ is the inclusion of $pd(2)$ in $2$-$Glob$
\end{itemize}

We then have that the following square is commutative, and is the one induced by the monad with arities $T$ (whose arity functor is $i_1$):
\begin{center}
\begin{tabular}{c} 
\xymatrix{
    fpd(2) \ar[r]^{i} & 2-Cat\\
    pd(2) \ar[u]^{l} \ar[r]_{i_1} & 2-Glob \ar[u]_{T}\\
    }
\end{tabular}
\end{center}

Since $i$ is fully faithful and dense, there is a generalized nerve functor for $2$-categories $2$-$Cat(i,1)$, sending a $2$-category $\mathcal{C}$ to the presheaf over $fpd(2)$ which sends an element of $fpd(2)$ to the set of $2$-categorical morphisms from it to $\mathcal{C}$.\\
The generalized Seagal condition which was previously introduced then characterizes the presheaves over $fpd(2)$ which are such nerves: such a presheaf $X \in \widehat{fpd(2)}$ is such a nerve if and only if there exists a $2$-globular set $G$ such that $l^* X$ is isomorphic to $2$-$Glob(i_1,G)$. This informally means that the restriction of $X$ to the elements of $fpd(2)$ (that are, the free categories over $2$-pasting diagrams) corresponds exactly to the set of morphisms from $2$-pasting diagrams to a given $2$-globular set $G$.

\subsubsection{Arities for $n$-categories}

The reader may have remarked that all that was done in the previous section may be generalized to higher dimensions. Again, the idea is to use pasting diagrams - which are defined in any dimension, and still compose in these dimensions - as arities for the free $n$-category monad. The collection of all the pasting diagrams of finite dimension $pd$ even is an $\omega$-category, and probably provides a category of arities for $\omega$-categories.

\section{Monads and computer science}
\label{compsci}
\subsection{Usual monads for semantics}
In semantics, the study of the mathematical structure of computational effects is a very important subject. Monads give a very convenient way to model them: from a set of \emph{values} $A$, one can deduce a set of \emph{computations} over them - which will be $TA$, the computational steps being simulated by the multiplication of the monad. It is frequent to use finitary monads, or even \emph{countable} monads, that is, monads whose arity functor has domain $FinSet \sqcup \{[\omega]\}$ ($[\omega]$ being the usual ordinal corresponding to $\mathbb{N}$). Several easy examples of such monads may be given:
\begin{itemize}
\item The \emph{partiality} monad: since a computation might infinitely loop and thus never end, there is a notion of \emph{diverging computation}, denoted $\bot$, which is added to the set of potential values, so that given a set of potential values $A$, we have $TA\,=\,A \sqcup \{\bot\}$ defining the monad of partiality.
\item The \emph{nondeterminism} monad: in computer science, it is very common to consider nondeterministic models, where at some point several paths of computation may be taken; one then considers that they were all taken simultaneously, giving a set of results instead of just one result to the computation. Thus this monad is given by $TA\,=\,\mathcal{P}_{fin}(A)$ (this being the set of finite parts of the set $A$).
\end{itemize}

But more complex monads are used in semantics, and it is interesting to have a dual viewpoint on these monads and on the corresponding Lawvere theories\footnote{Most of the following is extracted from \cite{hyland07}.}:

\begin{itemize}
\item The \emph{state} monad: given a set of \emph{states} $S$, it is defined by $TA\,=\,(S \times A)^S$. This may seem a little barbarian, but is better understood if we consider the arrows of its Kleisli category: an arrow $A \rightarrow B$ in this category corresponds to an arrow $A \rightarrow (S\times B)^S$, which decurrifies to $S\times A \rightarrow S \times B$: so, the monad takes a value and a state to a new value together with a new state. We will see in the next subsection that there is a nice way to characterize the corresponding Lawvere theory when the set of states defines a global store, that is when $S$ is of the shape $V^L$ ($V$ being a countable set of \emph{values}, indexed by a finite set of \emph{locations}: this is the usual description of memory in computer science) - with some operations on the store.
\item The \emph{exceptions} monad: in some programming languages, one may use \emph{exceptions}, that are, special calls which interrupt the normal run of the program to be handled. For example, an exception may be launched if an attempt to divide by zero is made: it prevents this division to happen and asks an exception handler what to do now (usually, just abort with an error). The monad is then $TA\,=\,A + E$, with $E$ a set of exceptions\footnote{One may remark that the monad of partiality was just a monad of exceptions, endowed with only one exception corresponding to the diverging computation. However this lacks a bit of computational sense, since a real machine couldn't launch this exception by indecidability of the halting problem.}. The corresponding Lawvere theory is the one generated by as many operations as there are in $E$: for every exception $e \in E$, a nullary operation\\\mbox{$raise_e\,:\,0\rightarrow 1$} is provided, corresponding to the raise of this exception. Later in this paper we will see that this monad is dual to the state monad, after a generalization of its definition, allowing also the handling of the risen exceptions.
\item The monad of \emph{interactive input/output} is a bit frightening at first sight but should be understood more easily than it seems. Roughly, it is defined by $TA\,=\,\mu Y.(O \times Y + Y^I + A)$, where $I$ is a countable set of inputs and $O$ a countable set of outputs. The idea is that, as in modal $\mu$-calculus\footnote{Refer e.g. to \cite[Section 2.1]{grellois2010} for an informal introduction.}, $\mu$ represents the finite iteration operator. The easiest way to understand this monad is to see it as a combination of:
\begin{itemize}
\item The \emph{interactive input} monad $TA\,=\,\mu Y.(Y^I + A)$: this describes the set of finite trees whose branches are labelled by $I$ and leaves by $A$ (so that every path from the top of the tree to a leaf may be understood as a succession of inputs followed by "answer", which is the value the program returns). It is a bit tricky to understand this tree shape for the reader unaware of the recursion operator $\mu$: the idea is that this formula describes some kind of a computational process: at first, there is a choice to do between $Y^I$ and $A$ (the $+$ acts as a disjunction). If $A$ is chosen, since there is no instance of $Y$ in it and since $Y$ denotes the value we loop on, the process is over and our tree is just a one-node tree, so it's a leaf, and it is from $A$ as we claimed. Else, a set of values from $Y$ is chosen, and the recursion goes on every of them, building a subtree on every of these values, until $A$ is finally chosen in every taken computational path. It has to happen, since $\mu$ only allows finite recursion.
\item The \emph{interactive output} monad $TA\,=\,\mu Y.(O\times Y + A)$ which generates a language isomorphic to $O^* \times A$ - it should be quite clear: the recursion starts with a choice between $O\times Y$ and $A$, if $A$ is chosen, it is finished, else $O\times Y$ is picked and $Y$ is given by this finite recursion process. So, $TA$ is the set of words of finite length built out of letters of $O$ and then of one and only one letter of $A$: it describes a program which writes a given number of outputs and then returns a value at the end of its execution. 
\end{itemize}
The combination of both monads leads to a structure where inputs and outputs alternate until a final value is returned by the program at the end of its execution. Without surprise, the associated (countable, there is an additional arity $[\omega]$) Lawvere theory is the one generated by the operations $read\,:\,I\rightarrow 1$ and $write\,:\,1 \rightarrow O$.
\item The nondeterminism may also be accounted of by Lawvere theories: for example, when the branching structure of nondeterminism is binary, the corresponding Lawvere theory is the one freely generated by a binary operation $\vee\,:\,2 \rightarrow 1$ subject to the equations of a semilattice (associativity, commutativity, idempotence). In \cite[Example 6.6]{hyland07} one may found a probabilistic nondeterminism monad with its associated theory.
\end{itemize}

\subsection{The state monad}
\subsubsection{The state monad is finitary}
Remember Section \ref{finfact}: a characterization of finitary monads \emph{via} the factorization of arrows exists; here we have, for $T$ the state monad, that every function $h\,:\,i_0[n] \rightarrow (S\times A)^S$ can be seen after decurryfication as a function $S\times [n] \rightarrow S \times A$, and such a function factors as:
\[
S\times [n] \xrightarrow{e} S \times [p] \xrightarrow{S\times f} S \times A
\]
with $f$ defined as an injective enumeration of the finite image of $h$. This factorization is unique modulo the "zig-zag" equivalence relation we described before. The state monad may thus be presented by operations of finite arities and equations between them, as we do next.

\subsubsection{An algebraic presentation of the state monad}
In \cite{plotkinpower}, Plotkin and Power give an explicit description of the algebras of the state monad in the particular case where the set of states is of the shape $S\,=\,V^L$. It was reformulated by Melli\`es in \cite{lics2010}, and we will follow his description after giving the original one. The idea is that the category of algebras for the state monad, when $S\,=\,V^L$, is equivalent to a global store, this being:

\begin{definition}[Global store on a category]
Given a category $\mathcal{C}$ with countable products, a finite set $L$, and a countable set $V$, define the category $GS(\mathcal{C})$ as follows: 
\begin{itemize}
\item An object of $GS(\mathcal{C})$ is given by:
\begin{itemize}
\item an object $A$ of $\mathcal{C}$
\item a "lookup" map $lookup\,:\,A^V \rightarrow A^L$
\item an "update" map $update\,:\,A\rightarrow A^{L\times V}$
\end{itemize} 
satisfying a series of commutative diagrams we will give later

\item An arrow $(A,update,lookup)\rightarrow (A',update',lookup')$ is just a map\\\mbox{$f\,:\,A\rightarrow A'$} in $\mathcal{C}$ subject to commutativity of $f$ with $lookup$ and $lookup'$ and to commutativity of $f$ with $update$ and $update'$.
\end{itemize}
\end{definition}
The $update$ and $lookup$ operations behave as follows: given a $V$-indexed family of elements of A, $lookup$ takes a location $loc$, gets the corresponding value in $S$, and computes the element of $A$ corresponding to this value. Given an element of $A$ together with a location $loc$ and a value $val$, $update$ changes the state by giving value $v$ to location $loc$, then allows the computation to run.\\

Instead of giving Plotkin and Power's conditions of commutativity on the objects of $GS(\mathcal{C})$, we give their reformulation by Melli\`es in the case where we consider instead of just $update$ and $lookup$ two families of applications:
\begin{itemize}
\item $update_{loc,val}\,:\,A \rightarrow A$ - these are unary operations
\item $lookup_{loc}\,:\,A^V \rightarrow A$ - these are $V$-ary operations
\end{itemize}

\begin{remark}[A set-theoretic interpretation, for the finite case]
Recall from Section \ref{intuitn} that to every finitary monad corresponds a Lawvere theory, whose $n$-ary operations correspond to the elements of the free algebra over $n$ elements $Tn$. When $V$ is finite, $S\,=\,V^L$ is too, and the corresponding state monad is finitary. So, every $n$-ary operation of the associated Lawvere theory corresponds to an element of $Tn\,=\,(S\times n)^S$, that is, to a set-theoretic function $S \rightarrow S \times n$. So, what is the meaning of our two families of operations in the Lawvere theory associated to the state monad at the light of this remark ?
\begin{itemize}
\item $update_{loc,val}$ is unary: so it is a map $S\rightarrow S\times 1 \cong S$, namely:
\begin{center}
\begin{tabular}{rcl}
$update_{loc,val}$ & $:$ & 
\begin{tabular}{rcl}
$S$ & $\rightarrow$ & $S$ \\
$state$ & $\mapsto$ & $state[loc:=val]$\\ 
\end{tabular}\\
\end{tabular}
\end{center}
\item $lookup_{loc}$ is a $V$-ary operation, and thus corresponds to a $S \rightarrow S \times V$ operation:
\begin{center}
\begin{tabular}{rcl}
$lookup_{loc}$ & $:$ & 
\begin{tabular}{rcl}
$S$ & $\rightarrow$ & $S\times V$ \\
$state$ & $\mapsto$ & $(state,state(loc))$\\ 
\end{tabular}\\
\end{tabular}
\end{center}
\end{itemize}
which meets the intuition behind the notion of a global store. Using the generalization of arities we described before in the paper leads to a similar interpretation of the case where $V$ is countable, for an arity functor $FinSet \sqcup \{[\omega ]\}$.
\end{remark}

The diagrammatic conditions on these arrows then are the following\footnote{The interpretation may seem weird, as arrows may look in the wrong order: the point is that we are here in the description of algebras of the state monad, whereas in the corresponding Lawvere theory, the models have arrows in the other direction (they are presheaves), which is the one corresponding to intuition.}:

\begin{enumerate}
\item \emph{annihilation update-update}: reading the value of a location $loc$ and then updating $loc$ with the value just read is equivalent to doing nothing. Defining $A^{val}\,:\,A^V \rightarrow A$ as the $val$-th projection of $A^V$ over $A$, and $update_{loc,V}$ as the unique morphism making the following diagram for every value $val \in V$:
\begin{center}
\begin{tabular}{c}
\xymatrix @!0 @R=22mm @C=22mm {
\ & A^V \ar[dr]^{A^{val}} & \ \\
A \ar[ur]^{update_{loc,V}} \ar[rr]_{update_{loc,val}} & \ & A\\
}
\end{tabular}
\end{center}
this expresses as requiring that the following diagram is commutative:
\begin{center}
\begin{tabular}{c}
\xymatrix @!0 @R=22mm @C=22mm {
\ & A^V \ar[dr]^{lookup_{loc}} & \ \\
A \ar[ur]^{update_{loc,V}} \ar[rr]_{id} &\ & A\\
}
\end{tabular}
\end{center}

\item \emph{interaction lookup-lookup}: reading twice a same location $loc$ can be reduced to just one reading operation:
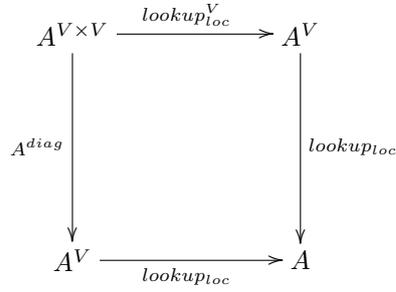
\begin{center}
\begin{tabular}{c}
\xymatrix @!0 @R=30mm @C=30mm {
A^{V\times V} \ar[d]_{A^{diag}} \ar[r]^{lookup^V_{loc}} & A^V \ar[d]^{lookup_{loc}}\\
A^V \ar[r]_{lookup_{loc}} & A\\
}
\end{tabular}
\end{center}

\item \emph{interaction update-update}: storing $val$ then $val'$ in $loc$ is the same as just storing $val'$:
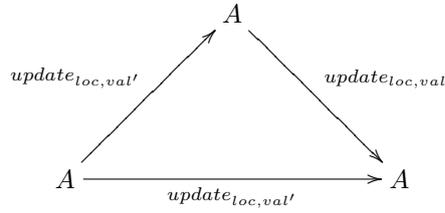
\begin{center}
\begin{tabular}{c}
\xymatrix @!0 @R=22mm @C=22mm {
\ & A \ar[dr]^{update_{loc,val}} & \ \\
A \ar[ur]^{update_{loc,val'}} \ar[rr]_{update_{loc,val'}} &\ & A\\
}
\end{tabular}
\end{center}

\item \emph{interaction update-lookup}: after storing $val$ in $loc$, reading $loc$ provides $val$:
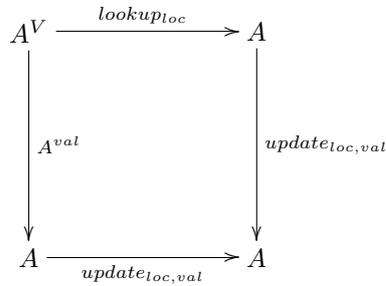
\begin{center}
\begin{tabular}{c}
\xymatrix @!0 @R=30mm @C=30mm {
A^V \ar[d]^{A^{val}} \ar[r]^{lookup_{loc}} & A \ar[d]^{update_{loc,val}} \\
A \ar[r]_{update_{loc,val}} & A\\
}
\end{tabular}
\end{center}

\item \emph{commutation lookup-lookup}: the order of looking at $loc$ and $loc'$ is irrelevant:
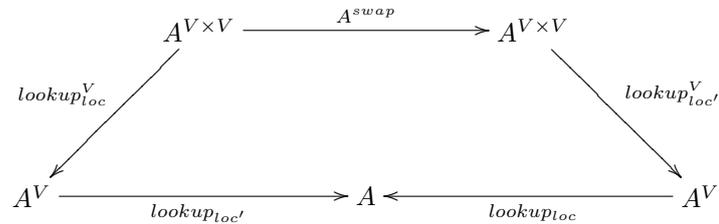
\begin{center}
\begin{tabular}{c}
\xymatrix @!0 @R=22mm @C=22mm {
\ & A^{V\times V} \ar[dl]_{lookup_{loc}^V} \ar[rr]^{A^{swap}} &\ & A^{V\times V} \ar[dr]^{lookup_{loc'}^V}   & \ \\
A^V  \ar[rr]_{lookup_{loc'}} & \ & A & \ & A^V \ar[ll]^{lookup_{loc}}\\
}
\end{tabular}
\end{center}

\item \emph{commutation update-update}: the order of storing values in different locations $loc \neq loc'$ is irrelevant:
\begin{center}
\begin{tabular}{c}
\xymatrix @!0 @R=30mm @C=30mm {
A \ar[r]^{update_{loc,val}} \ar[d]_{update_{loc',val'}} & A \ar[d]^{update_{loc',val'}}\\
A \ar[r]_{update_{loc,val}} & A\\
}
\end{tabular}
\end{center}

\item \emph{commutation update-lookup}: reading $loc$ and updating $loc' \neq loc$ can be performed in any order:
\begin{center}
\begin{tabular}{c}
\xymatrix @!0 @R=30mm @C=30mm {
A^V \ar[d]_{update^V_{loc,val}} \ar[r]^{lookup_{loc'}} & A \ar[d]^{update_{loc,val}}\\
A^V \ar[r]_{lookup_{loc'}} & A\\
}
\end{tabular}
\end{center}
\end{enumerate}

The interesting point is the following:

\begin{theorem}[Characterization of the algebras of the state monad]
The categories $GS(\mathcal{C})$ and $\mathcal{C}^T$ are equivalent.
\end{theorem}

In other words, every algebra of the state monad - that is, every \emph{concrete realization} of the computational process it describes - is an object $A$ of $\mathcal{C}$ together with the data of $lookup$ and $update$ operations over the set of states $S\,=\,V^L$.\\\ \\
There are (at least) two proofs of this theorem:
\begin{itemize}
\item In the original paper by Plotkin and Power \cite[Section 3]{plotkinpower}, the Beck theorem is used to show that the forgetful functor $GS(\mathcal{C})\rightarrow \mathcal{C}$ is \emph{monadic} for the state monad $TA\,=\,(S \times A)^S$: informally, this means that this functor is like the forgeful functor $\mathcal{C}^T\rightarrow \mathcal{C}$, and in particular that $GS(\mathcal{C})$ and $\mathcal{C}^T$ are equivalent categories: as claimed, global stores are the algebras of the state monad.
\item In the paper by Melli\`es \cite[Section VII]{lics2010}, the use of the Beck theorem is avoided by showing instead that the algebraic presentation of objects with global stores that was given by the diagrammatic conditions provides a presentation by generators and relations of the Lawvere theory associated to the state monad - and this by means of rewriting. 
\end{itemize}

\subsubsection{A remark on the exceptions monad}

Duval recently exposed \cite{duval1} with Dumas, Fousse and Reynaud \cite{duval2} an interesting generalization of the exceptions monad, which is dual to the state monad. The idea is to replace the usual formulation $TA\,=\,A + E$ (with $E$ a set of arities), which leads in the Kleisli category to arrows $A \rightarrow B + E$, by a formulation of the monad which leads to arrows $A + E \rightarrow B + E$ in the Kleisli category. The reader aware of exceptions may have remarked that this allows, in addition to the usual \emph{raise} of exceptions, their treatment by the use of \emph{catch} operations. So, the exceptions monad should be reformulated to $TA\,=\,E\cdot (A+E)$ (remember that $\cdot$ stands for the copower). Since $\cdot$ is dual to the usual power, and $+$ to $\times$, we have a duality of these two monads. A presentation by generators and relations of the algebras of the exceptions monad is provided; it would probably be interesting to use the technique Melli\`es used for the state monad to prove that the exceptions monad is really the one corresponding to the given equational theory.

\section{Conclusion}
We have investigated the notion of arity and its interest for monads starting from the usual correspondence between Lawvere theories and finitary monads, extending the situation to the more general case suggested by the dual notions of nerve and realization. The definition of monads with arities notably gives a definition of $n$-categories as $n$-globular sets satisfying a Segal condition, generalizing the $n=1$ case, but also covers the usual generalizations in the practice of semantics, such as the use of countable monads and theories. It would probably be interesting to have a look at how this generalization enriches semantics, by allowing the specification of theories over trees for example: is there a link with the notion of tree automata ?\\
Another interesting question would be to relate the work on the definition, using operads, of $n$-categories - including the weak case - by Leinster \cite{leinster} to a definition using monads with arities extending the work of Section \ref{higherarities}.

\appendix
\section{Further discussion: Operads and algebraic theories}
\subsection{Operads}
In a general category, arrows just have one object as both domain and codomain. Operads enriches this with arrows from any power of an object to the object itself - just like in linear algebra one can consider operations on vectors $\mathbb{R}^n \rightarrow \mathbb{R}$ instead of just morphisms $\mathbb{R}\rightarrow \mathbb{R}$. A similarity with Lawvere theories appear here, as we consider operations $A^n \rightarrow A$ for some object $A$. A reference on operads is \cite{leinster}.

\begin{definition}[Operad]
An operad $\mathcal{C}$ consists of:
\begin{itemize}
\item An object $A$,
\item For each $n\in\mathbb{N}$, a set of arrows $\mathcal{C}(n)$ - where for example an element of $\mathcal{C}(2)$ is depicted as follows:
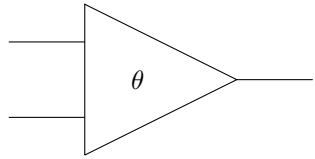
\begin{center}
\begin{tikzpicture}
\draw (-1,-1) -- (1, 0) -- (-1,1) -- (-1,-1);
\draw (1,0) -- (2,0);
\draw (-2,-0.5) -- (-1,-0.5);
\draw (-2,0.5) -- (-1,0.5);
\draw (-0.3,0) node {$\theta$};
\end{tikzpicture}
\end{center}
\item A notion of composition of these arrows: for $n,k_1,\cdots,k_n \in \mathbb{N}$, it is a function mapping an arrow of $P(n)$ together with an arrow from $P(k_i)$ per $i \in \{1,\cdots ,n\}$ to an arrow of $P(k_1 + \cdots + k_n)$; on a graphical example, it corresponds to building out of the three following arrows:
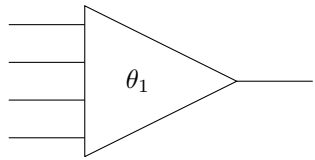
\begin{center}
\begin{tikzpicture}
\draw (-1,-1) -- (1, 0) -- (-1,1) -- (-1,-1);
\draw (1,0) -- (2,0);
\draw (-2,-0.75) -- (-1,-0.75);
\draw (-2,0.75) -- (-1,0.75);
\draw (-2,-0.25) -- (-1,-0.25);
\draw (-2,0.25) -- (-1,0.25);
\draw (-0.3,0) node {$\theta_1$};
\end{tikzpicture}
\end{center}
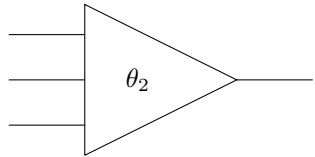
\begin{center}
\begin{tikzpicture}
\draw (-1,-1) -- (1, 0) -- (-1,1) -- (-1,-1);
\draw (1,0) -- (2,0);
\draw (-2,-0.6) -- (-1,-0.6);
\draw (-2,0.6) -- (-1,0.6);
\draw (-2,0) -- (-1,0);
\draw (-0.3,0) node {$\theta_2$};
\end{tikzpicture}
\end{center}
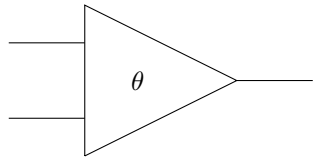
\begin{center}
\begin{tikzpicture}
\draw (-1,-1) -- (1, 0) -- (-1,1) -- (-1,-1);
\draw (1,0) -- (2,0);
\draw (-2,-0.5) -- (-1,-0.5);
\draw (-2,0.5) -- (-1,0.5);
\draw (-0.3,0) node {$\theta$};
\end{tikzpicture}
\end{center}
the following arrow:
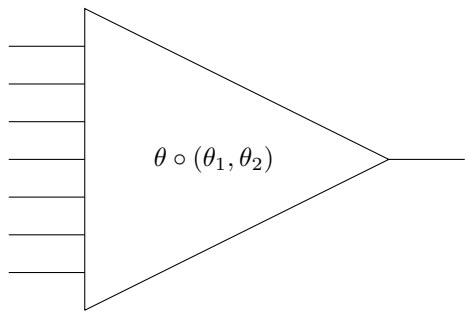
\begin{center}
\begin{tikzpicture}
\draw (-2,-2) -- (2, 0) -- (-2,2) -- (-2,-2);
\draw (2,0) -- (3,0);
\draw (-3,-1.5) -- (-2,-1.5);
\draw (-3,-1) -- (-2,-1);
\draw (-3,-0.5) -- (-2,-0.5);
\draw (-3,0) -- (-2,-0);
\draw (-3,0.5) -- (-2,0.5);
\draw (-3,1) -- (-2,1);
\draw (-3,1.5) -- (-2,1.5);
\draw (-0.3,0) node {$\theta \circ (\theta_1, \theta_2 )$};
\end{tikzpicture}
\end{center}
to be thought of as:
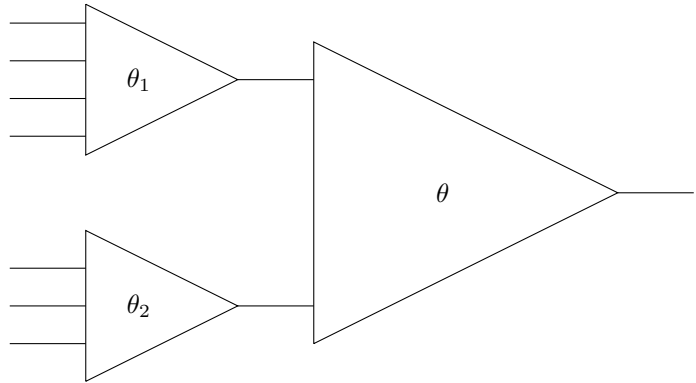
\begin{center}
\begin{tikzpicture}
\draw (0,-2) -- (4, 0) -- (0,2) -- (0,-2);
\draw (4,0) -- (5,0);
\draw (-1,1.5) -- (0,1.5);
\draw (-1,-1.5) -- (0,-1.5);
\draw (-3,-2.5) -- (-1, -1.5) -- (-3,-0.5) -- (-3,-2.5);
\draw (-4,-1.5) -- (-3,-1.5);
\draw (-4,-2) -- (-3,-2);
\draw (-4,-1) -- (-3,-1);
\draw (-3,2.5) -- (-1, 1.5) -- (-3,0.5) -- (-3,2.5);
\draw (-4,2.25) -- (-3,2.25);
\draw (-4,1.75) -- (-3,1.75);
\draw (-4,1.25) -- (-3,1.25);
\draw (-4,0.75) -- (-3,0.75);
\draw (1.7,0) node {$\theta$};
\draw (-2.3,-1.5) node {$\theta_2$};
\draw (-2.3,1.5) node {$\theta_1$};
\end{tikzpicture}
\end{center}
Moreover, this composition law has to have an identity $i \in \mathcal{C}(1)$ and to obey an associativity property.
\end{itemize}
\end{definition}

Beware that in the general definition there is no action of the symmetric group $S_n$ over $\mathcal{C}(n)$: when it is the case, the operad is called \emph{symmetric}.

\begin{remark}[Multicategories]
It is easy to define a version of operads where all wires do not necessarily correspond to the same object: this gives the notion of multicategory, see \cite[Section 2.1]{leinster} - operads being then a degenerated case of multicategories where all wires correspond to the same object. Multicategories are to operads what multisorted Lawvere theories are to usual Lawvere theories.
\end{remark}

Monads have \emph{algebras}, Lawvere theories have \emph{models}; operads have as well a notion of concrete realization, called algebra for the operad.

\begin{definition}[Algebra for an operad]
Given an operad $\mathcal{C}$, an algebra for it is a set $X$ together with a function $\widehat{\theta}\,:\,X^n\rightarrow X$ for every $n \in \mathbb{N}$ and $\theta \in \mathcal{C}(n)$, so that these functions respect identity and obey an associativity law.
\end{definition}

\begin{example}[Some operads and their algebras]\ 
\label{operadsalgeb}
\begin{itemize}
\item The terminal operad $1$ has for all $i \in \mathbb{N}$ exactly one arrow of arity $i$. An algebra for it is then a set $X$ together with a function:
\begin{center}
\begin{tabular}{rcl}
$X^n$ & $\rightarrow$ & $X$\\
$(x_1,\cdots,x_n)$ & $\mapsto$ & $(x_1 \cdot \ldots \cdot x_n)$\\
\end{tabular}
\end{center}
and these functions have to respect composition, associativity and identities, that is, to obey the two following requirements:
\begin{itemize}
\item $((x_1^{1} \cdot \ldots \cdot x_1^{k_1}) \cdot \ldots \cdot (x_n^{1} \cdot \ldots \cdot x_n^{k_n}))\,=\,(x_1^1 \cdot \ldots \cdot x_n^{k_n})$
\item $x\,=\,(x)$
\end{itemize}
so that the category of algebras for this operad is the category of monoids.
\item The operad which has only one unary operation and no operations of other arities has sets as its algebras.
\item The operad which has exactly one $n$-ary operation for $n \geq 1$ and no nullary operation has semigroups as its algebras - since this operad has no nullary operation, an algebra for it has no operation $1\rightarrow X$ picking an unit element.
\end{itemize}
\end{example}

\subsection{Operads and algebraic theories}
Example \ref{operadsalgeb} shows that some algebraic theories can be described by operads, but actually not all of them can be. In \cite{operadstopo}, an example of Lawvere theory that is not operadic (that is, such that there exists no operad whose category of algebras is equivalent to the category of models of the theory) is provided: the theory of Boolean  rings, where $x^2\,=\,x$, cannot be described as an algebra of an operad.\\
The algebraic theories which can be represented by operads are, according to \cite[Appendix C.1]{leinster}, the finitary ones which admit a strongly regular presentation, that is, the finitary algebraic theories whose relations have the same variables occuring without repetition and in the same order on each side of the equality predicate. For example, the following relations are strongly regular:
\[
x\cdot (y \cdot z)\,=\,(x\cdot y)\cdot z,\ \ x \cdot 1\,=\, x,\ \ (x^y)^z\,=\,x^{y\cdot z}
\]
whereas the following relations are not:
\[
x \cdot 0\,=\,0,\ \ x\cdot y\,=\,y \cdot x,\ \ x\cdot (y+z)\,=\,x\cdot y + x \cdot z
\]
So, the theory of monoids is strongly regular, and can thus be presented using an operad as we have seen it in Example \ref{operadsalgeb}. But where we had a Lawvere theory for groups in Example \ref{lawgroup}, there is no presentation of the theory of groups using operads, since $x^{-1}\cdot x\,=\,1$ is not strongly regular and since it turns out that this relation can not be replaced by a strongly regular equivalent one (see \cite[Example 4.1.6]{leinster} for indications towards a proof).\\
But operads are not just a special kind of algebraic theories: in \cite{leinster_algeb}, Leinster exhibits a pair of operads which are not isomorphic but induce isomorphic monads (corresponding to isomorphic algebraic theories). The idea is, given an operad $\mathcal{C}$, to define its \emph{reverse} $\mathcal{C}^{rev}$, where the only change is in the composition law: 
\[
\theta \circ_{rev} (\theta_1, \cdots, \theta_n)\,=\,\theta \circ (\theta_n,\cdots ,\theta_1)
\]
Basically, the only change is on the order of plugging of the inputs when composing - so we seek an operad, obviously nonsymmetric, whose reverse is not isomorphic to itself - since we have that the monad induced by an operad an its reverse are isomorphic, see \cite{leinster_algeb} - the monad induced by an operad being:

\begin{definition}[Monad induced by an operad] Given an operad $\mathcal{C}$, let $T_{\mathcal{C}}$ be the following endofunctor of $Set$:
\[
T_{\mathcal{C}}\,:\,X \mapsto \coprod_{n\in \mathbb{N}} \mathcal{C}(n) \times X^n
\]
$T_{\mathcal{C}}$ is a monad together with a multiplication induced by the composition law of the operad, unit being induced by the identity of the operad.
\end{definition}

\begin{proposition}
The algebras for $\mathcal{C}$ and $T_{\mathcal{C}}$ are the same.
\end{proposition}

An operad which is not isomorphic to its reverse is given by Leinster in \cite[Section 2]{leinster_algeb}. Finally, even if Lawvere theories and operads can both model some theories, they do not reduce to each other.

\subsection{Operads and generalized arities}
An interesting generalization of our first definition of operads is given by Leinster in \cite[Section 4.2]{leinster}: the point is to replace the usual notion of arity, given by integers, by a notion of arity given by a special class of monads - the cartesian ones. It is thus interesting to discuss a little about this widening of the notion of arity.

\begin{definition}[Cartesian monad]\ 
\begin{itemize}
\item A category is cartesian if it has all pullbacks
\item A functor whose domain is a cartesian category is cartesian if it preserves pullbacks
\item A natural transformation is cartesian if all its naturality squares are pullbacks
\item A monad $(T,\mu,\eta)$ on a category $\mathcal{C}$ is cartesian if $\mathcal{C}$, $T$, $\mu$ and $\eta$ are cartesian
\end{itemize}
\end{definition}

These monads will be used to describe arities, the preservation of pullbacks allowing to deal properly with spans as we will see soon. Here are some examples of our to-be arities:

\begin{example}[Some cartesian monads]\ 
\begin{itemize}
\item The identity monad on a cartesian category is cartesian - so, $FinSet$ being cartesian, the usual notion of arities will still correspond to a notion of arity in our generalization
\item The monad of trees with labelled leaves is cartesian (given a set $A$ of labellings, $TA$ is built inductively by the rules "if $a\in A$ then $a \in TA$" and "if $t_1,\cdots,t_n \in TA$ then $(t_1,\cdots ,t_n) \in TA$", the unity of the monad being obvious, and its multiplication corresponding to append trees)
\item The free strict $n$-category monad on a $n$-globular set, and the free strict $\omega$-category monad on a globular set are cartesian.
\end{itemize}
\end{example}

We can now define the notion of $T$-operad, for a cartesian monad $T$. It uses the notion of \emph{span}, whose idea is that a category can be defined as a set of arrows together with their source and target projections which implicitely give the objects.

\begin{definition}[$T$-operad]
Given a cartesian monad $T$ on a cartesian category $\mathcal{C}$, and a final object $C_0$ of $\mathcal{C}$, a $T$-operad is a diagram:
\begin{center}
\begin{tabular}{c}
  \xymatrix @R=2pc @C=1.5pc{
  \ & C_1 \ar[dl]_{dom} \ar[dr]^{cod} & \ \\
  TC_0 & \ & C_0\\
  }
\end{tabular}
\end{center}
together with maps $C_1 \circ C_1\,=\,C_1 \times_{TC_0} TC_1 \xrightarrow{comp} C_1$ and $C_0 \xrightarrow{ids} C_1$ satisfying identity and associativity axioms.
\end{definition}

The idea is that the composition map is induced by the cartesian structure of the monad, which makes the following diagram - defining $C_1 \circ C_1$ - commute:\\
\begin{center}
\begin{tabular}{c}
  \xymatrix @R=2pc @C=1.5pc{
\ &\ &\ &C_1 \circ C_1 \ar[dl] \ar[dr] &\ &\ \\
\ &\ & TC_1 \ar[dl]_{Tdom} \ar[dr]^{Tcod} &\ &C_1 \ar[dl]_{dom} \ar[dr]^{cod} &\ \\
\ &T^2C_0 \ar[dl]_{\mu_{C_0}} &\ &TC_0&\ &C_0\\
TC_0 & \ & \ & \ &\ &\ \\
  }
\end{tabular}
\end{center}

where the upper square, defining $C_1 \circ C_1$, is obtained as a pullback.

\paragraph{Globular operads and weak $\omega$-categories} If $T$ is the free strict $\omega$-category monad and $1$ the terminal globular set:
\begin{center}
\begin{tabular}{ccc}
$1$ & $\ =\ $ &
\xymatrix {
\cdots \ar@/^/[r] \ar@/_/[r] & 1 \ar@/^/[r] \ar@/_/[r] & \cdots \ar@/^/[r] \ar@/_/[r] & 1 \\
}\\
\end{tabular}
\end{center}
then $T1$ is the set $pd$ of pasting diagrams\footnote{As described at Section \ref{pd}.} of all dimensions - the $n$-pasting diagrams, whose set is denoted by $pd(n)$, being the elements for $(T1)(n)$.

\begin{definition}[Globular operad]
A globular operad is a $T$-operad.
\end{definition}

So, a globular operad $\mathcal{C}$ can be thought of as some kind of theory whose arities are the pasting diagrams. The notion of composition is the one we have seen before - paste the diagrams together according to the composition shape described by a pasting diagram. As mentioned before, the algebras of an operad are the ones of a monad; here it turns out to be the monad $T_{\mathcal{C}}$ defined on globular sets $X$ (presheaves over $\mathbb{G}$, that is $\mathbb{G}_n$ with $n = \omega$) by:
\[
(T_{\mathcal{C}} X)(n)\ \cong\ \coprod_{\pi \in pd(n)}{\mathcal{C}(\pi )\times [\mathbb{G}^{op},Set](\widehat{\pi},X)}
\]

The idea is that an algebra for this monad, and thus for the corresponding globular operad, just amounts to label the pasting diagrams by cells of the globular set $X$. In \cite[Section 9.2]{leinster}, a particular notion of globular operads, \emph{operads with contraction}, is used to give an operad whose algebras are the weak $\omega$-categories. The process is then extended to the definition of weak $n$-categories.

\newpage

\bibliographystyle{plain}
\bibliography{rapport}

\end{document}